\documentclass[11pt,letterpaper]{amsart}

\usepackage{amsmath}
\usepackage{amssymb}
\usepackage{mathrsfs}
\usepackage{graphicx}
\usepackage{etoolbox}

\makeatletter
\patchcmd{\@settitle}{\uppercasenonmath\@title}{}{}{}
\makeatother

\usepackage[table]{xcolor}
\usepackage{tikz}

\usepackage{enumitem}

\usepackage[
  letterpaper,
  textwidth=6in,
  textheight=8.5in,
  centering
]{geometry}

\numberwithin{equation}{section}
\allowdisplaybreaks

\title[Stability of composite waves for the zero heat-conductivity system]
{Nonlinear stability of a rarefaction-viscous shock composite wave
for the compressible Navier-Stokes equations with zero heat conductivity}

\author[Y. Guo]{Yilin Guo}
\address[Yilin Guo]{\newline School of Mathematical Sciences, Laboratory of Mathematics and Complex Systems, MOE, Beijing Normal University, Beijing 100875, P. R. China.}
\email{gyl@mail.bnu.edu.cn}

\author[L.-A. Li]{Lin-An Li}
\address[Lin-An Li]{\newline School of Mathematical Sciences, Laboratory of Mathematics and Complex Systems, MOE, Beijing Normal University, Beijing 100875, P. R. China.}
\email{linanli@amss.ac.cn}

\author[J. Wu]{Jiahong Wu}
\address[Jiahong Wu]{\newline Department of Mathematics, University of Notre Dame, Notre Dame, IN 46556, USA.}
\email{jwu29@nd.edu}

\author[X. Xu]{Xiaojing Xu}
\address[Xiaojing Xu]{\newline School of Mathematical Sciences, Laboratory of Mathematics and Complex Systems, MOE, Beijing Normal University, Beijing 100875, P. R. China.}
\email{xjxu@bnu.edu.cn}

\date{\today}

\theoremstyle{plain}

\newtheorem{theorem}{Theorem}[section]
\newtheorem{lemma}{Lemma}[section]
\newtheorem{proposition}{Proposition}[section]
\newtheorem{remark}{Remark}[section]
\newtheorem{corollary}{Corollary}[section]

\newcommand{\s}{\sigma}
\newcommand{\x}{\xi}
\newcommand{\eps}{\varepsilon}
\newcommand{\deltas}{\delta_S}
\newcommand{\deltar}{\delta_R}

\newcommand{\lam}{\lambda}
\newcommand{\bbr}{\mathbb {R}}
\newcommand{\mb}{\mathbf}

\begin{document}

\begin{abstract}
We prove the nonlinear stability of a composite wave consisting of a rarefaction wave and a viscous shock for the one-dimensional compressible Navier-Stokes equations with zero heat conductivity. This system is partially dissipative and violates the Shizuta-Kawashima condition: the entropy mode is undamped. To the best of our knowledge, this is the first stability result involving a viscous shock in this setting; the stability of a single viscous shock is included as a special case. We show that the
pressure and velocity converge uniformly to the composite wave, up to a time-dependent shift of the shock, whereas the entropy perturbation remains uniformly small, which cannot be improved to convergence in general. The proof combines the $a$-contraction method with shifts with two new ingredients: an effective viscous flux, which transfers the velocity dissipation to the pressure and allows us to close the estimates at the $H^1$ level, and weighted space-time estimates for the entropy derivative
along characteristics, which exploit the fact that the shock moves relative to the stationary entropy characteristics. This study is motivated by Liu's discussion of partially dissipative systems in \cite{LiuShockWaves}.
\end{abstract}


\maketitle 
\tableofcontents

\section{Introduction} 
We consider the one-dimensional compressible Navier-Stokes equations
with zero heat conductivity, which in Lagrangian coordinates are given
by
\begin{equation}\label{VUE-1}
    \begin{cases}
        v_t-u_x=0,
        & x\in\mathbb{R},\quad t>0,\\[1mm]
        u_t+p_x
        =\left(\mu\frac{u_x}{v}\right)_x,\\[1mm]
        E_t+(pu)_x
        =\left(\mu\frac{uu_x}{v}\right)_x.
    \end{cases}
\end{equation}
Here $v>0$, $u$, and $E$ denote the specific volume, velocity, and
total energy, respectively, and $\mu>0$ is the viscosity
coefficient. The total energy is given by $E=e+\frac{u^2}{2}$,
where $e$ denotes the internal energy. The pressure is given by $p=p(v,e)$ and satisfies $p_e(v,e)>0$ and $p_v(v,e)<0$.

System \eqref{VUE-1} with zero heat conductivity can be derived as a
first-order Chapman-Enskog approximation of a simple relaxation model
for gas dynamics in thermal non-equilibrium. For non-monatomic gases,
besides the translational energy, there are rotational, vibrational, and
other forms of internal energy. These energies may be approximately in
fixed proportions for small temperature variations, but this
approximation generally fails in high-temperature gas dynamics. A
simple non-equilibrium model assumes that the vibrational energy is not
in equilibrium, while the translational and rotational energies are in
equilibrium and remain in fixed proportions. The vibrational energy
then relaxes toward its equilibrium value within a relaxation time. As
the relaxation time tends to zero, one obtains the equilibrium
compressible Euler system, while retaining the first-order correction
through the Chapman-Enskog expansion gives system \eqref{VUE-1} with
positive viscosity and zero heat conductivity; see Liu
\cite{LiuShockWaves} and Vincenti-Kruger \cite{VK65}.

From the viewpoint of dissipation, system \eqref{VUE-1} is a partially
dissipative hyperbolic-parabolic system due to the absence of
temperature dissipation. Similar partially dissipative structures also
arise in several physically relevant models, including
magneto-hydrodynamic systems in which some dissipation parameters
vanish.

For simplicity, in this paper we consider the ideal polytropic gas, for
which the pressure, internal energy, and entropy satisfy
\begin{equation}\label{pes}
    p=\frac{R\theta}{v}
    =A v^{-\gamma}
      \exp\left(\frac{\gamma-1}{R}s\right),
    \quad
    e=\frac{R}{\gamma-1}\theta+\mathrm{const},
    \quad
    s
    =R\log v
    +\frac{R}{\gamma-1}\log\theta+\mathrm{const},
\end{equation}
where $A>0$, $R>0$, and $\gamma>1$ are constants. We then consider the initial values $(p_0,u_0,s_0)(x)$ satisfying
\begin{equation}\label{p0-s0}
    (p,u,s)(0,x)=(p_0,u_0,s_0)(x)
    \longrightarrow
    (p_\pm,u_\pm,s_\pm),
    \qquad x\to\pm\infty,
\end{equation}
where the far-field states $p_\pm>0$, $u_\pm\in\mathbb{R}$, and
$s_\pm\in\mathbb{R}$ are given constants. We study the long-time
behavior of solutions to \eqref{VUE-1} with initial values satisfying
\eqref{p0-s0} toward the composite wave determined by the Riemann solution of the corresponding Euler system
\begin{equation}\label{E}
\begin{cases}
v_t-u_x=0,&x\in\mathbb{R},\ t>0,\\[1mm]
u_t+p_x=0,\\[1mm]
E_t+\bigl(pu\bigr)_x=0,
\end{cases}
\end{equation}
where $v=v(p,s)$ and $E=e(p,s)+\frac{u^2}{2}$ are determined by the
thermodynamic relations \eqref{pes}, with the corresponding Riemann
data
\begin{equation}\label{Euler-pus-data}
    (p,u,s)(0,x)=
    \begin{cases}
        (p_-,u_-,s_-), & x<0,\\[1mm]
        (p_+,u_+,s_+), & x>0.
    \end{cases}
\end{equation}
The Riemann solution to \eqref{E}-\eqref{Euler-pus-data} under consideration consists of a rarefaction wave
and a shock wave.

To review the related stability results under different dissipative
mechanisms, we first recall the general one-dimensional compressible
Navier-Stokes-Fourier system in Lagrangian coordinates:
\begin{equation}\label{NSF}
    \begin{cases}
        v_t-u_x=0, & x\in\mathbb{R},\ t>0,\\[1mm]
        u_t+p(v,\theta)_x
        =\left(\mu\frac{u_x}{v}\right)_x,\\[1mm]
        E_t+\bigl(p(v,\theta)u\bigr)_x
        =\left(\mu\frac{uu_x}{v}
        +\kappa\frac{\theta_x}{v}\right)_x,
    \end{cases}
\end{equation}
where $\kappa\geq0$ denotes the heat conductivity coefficient.

For the classical compressible Navier-Stokes-Fourier system with $\mu>0$ and $\kappa>0$ and for related viscous models, the nonlinear stability of elementary waves and composite waves has been extensively studied. In particular, the stability of viscous shock waves, rarefaction waves, viscous contact waves, and various composite waves has been established in a large
number of works; see, for example, \cite{MN85,Liu85,Good,SX93,LZ15,Xin-c,MN86,LX88,MN92,NYZ04,HMX,HXY,HLM,HM,HW16,HKK,KVW23,HL25} and the references therein. Recently, Kang-Vasseur-Wang \cite{KVW25} employed the $a$-contraction method with shifts, together with relative entropy and energy estimates, to establish the time-asymptotic stability of generic Riemann solutions consisting of a viscous shock wave, a viscous contact wave, and a rarefaction wave.

When $\mu=0$ and $\kappa>0$, system \eqref{NSF} reduces to the system
of a heat-conductive ideal gas without viscosity, for which the momentum
equation has no direct viscous dissipation and the dissipative effect
is provided solely by heat conduction. For viscous contact waves,
Ma \cite{Ma13} and Ma-Wang \cite{MW} established nonlinear stability
results. For rarefaction waves, Hou \cite{Hou} proved their asymptotic
stability. For composite waves, Fan-Matsumura \cite{FM} established
the time-asymptotic stability of a composite wave composed of two
viscous shock waves. Fan-Gong-Tang \cite{FGT} further proved the
asymptotic stability of a composite wave consisting of a viscous
contact wave and rarefaction waves. Recently, Guo-Li-Wu-Xu
\cite{GLWX} established the time-asymptotic stability of generic
Riemann solutions composed of a viscous shock wave, a viscous contact
wave, and a rarefaction wave.

A common structural feature of the two regimes discussed above is
that both the classical Navier-Stokes-Fourier system with $\mu>0$
and $\kappa>0$ and the heat-conductive ideal gas system with $\mu=0$
and $\kappa>0$ satisfy the Shizuta-Kawashima condition \cite{SK1,SK2,ShK85}. It allows dissipation to be transferred through the hyperbolic coupling to components without direct diffusion, thereby providing a fundamental mechanism for closing energy estimates and establishing global existence and stability of small perturbations. In contrast, when $\mu>0$ and $\kappa=0$, this condition fails because the eigenvector associated with the entropy characteristic field lies in the kernel of the viscosity matrix. Although viscosity still dissipates the acoustic modes, the entropy mode remains nondissipative. Consequently, entropy perturbations need not decay, and the standard compensating energy estimates do not yield time-integrable control of their spatial derivatives. This lack
of dissipation creates substantial difficulties in controlling
nonlinear interactions and precludes a direct application of the
classical stability theory based on the Shizuta-Kawashima condition.
Determining whether, and in what sense, nonlinear stability persists
beyond this condition is therefore an important and challenging
problem in the theory of partially dissipative systems; see also
the discussion in Liu \cite{LiuShockWaves}.

For the zero heat-conductivity system considered here, several
stability results have nevertheless been established despite this
structural degeneracy. Liu-Zeng \cite{LZ99} studied the long-time
behavior of solutions near constant states and obtained detailed
pointwise estimates by combining the Green's function method with
energy estimates. For rarefaction waves, Hu-Yin \cite{HY18}
established their nonlinear stability. However, the stability
theory for viscous shock waves, contact discontinuities, and
composite waves in this regime remains much less developed,
particularly when interactions between dissipative acoustic modes
and the nondissipative entropy mode must be taken into account.

In this paper, we establish the nonlinear stability of a composite
wave consisting of a 1-rarefaction wave and a viscous 3-shock wave
for system \eqref{VUE-1}. Although the system fails to satisfy the
Shizuta-Kawashima condition, we can close the energy estimates by combining the $a$-contraction method with shifts, an effective viscous flux, the wave structure of pressure and weighted estimates for the entropy derivative. This yields uniform convergence of the pressure and velocity to the corresponding shifted composite wave and uniform stability of the entropy. To the best of our knowledge, this is the first nonlinear stability result involving a viscous shock wave for a system that violates the Shizuta-Kawashima condition.

For the subsequent stability analysis, it is useful to employ two
equivalent formulations of \eqref{VUE-1}. The system
\eqref{v-u-theta} is mainly used to establish the basic $L^2$-estimates,
while \eqref{p-u-s}  is more convenient for the estimates of the first-order derivatives. Using the thermodynamic relations in \eqref{pes}, we present the two equivalent systems below:
\begin{equation}\label{v-u-theta}
    \begin{cases}
        v_t-u_x=0, & x\in\mathbb{R},\ t>0,\\[1mm]
        u_t+p(v,\theta)_x
        =\left(\mu\frac{u_x}{v}\right)_x,\\[1mm]
        \frac{R}{\gamma-1}\theta_t+p(v,\theta)u_x
        =\mu\frac{u_x^2}{v},
    \end{cases}
\end{equation}
and
\begin{equation}\label{p-u-s}
    \begin{cases}
        \displaystyle
        p_t+\frac{\gamma p}{v}u_x
        =\mu(\gamma-1)\frac{u_x^2}{v^2},
        & x\in\mathbb{R},\ t>0,\\[1mm]
        u_t+p_x
        =\left(\mu\frac{u_x}{v}\right)_x,\\[1mm]
        s_t=\mu\frac{u_x^2}{v\theta}.
    \end{cases}
\end{equation}

\subsection{Riemann problem for the inviscid system}
\label{subsec:Riemann}

We now recall the Riemann problem for the corresponding inviscid system
and the elementary wave curves associated with the composite wave. We first rewrite the initial values \eqref{p0-s0}, given in the $(p,u,s)$-coordinates, in the $(v,u,\theta)$-coordinates. By the thermodynamic
relations, the corresponding initial specific volume $v_0$ and temperature $\theta_0$ are given by
\begin{equation*}
v_0
=\left(\frac{A}{p_0}e^{\frac{\gamma-1}{R}s_0}\right)^{1/\gamma},
\qquad
\theta_0=\frac{p_0v_0}{R}.
\end{equation*}
It follows that
\begin{equation}\label{initial-v}
(v,u,\theta)(0,x)=(v_0,u_0,\theta_0)(x)
\longrightarrow
(v_\pm,u_\pm,\theta_\pm),
\qquad x\to\pm\infty,
\end{equation}
where $v_\pm=\left(\frac{A}{p_\pm}e^{\frac{\gamma-1}{R}s_\pm}\right)^{1/\gamma}$ and $\theta_\pm=\frac{p_\pm v_\pm}{R}$.
Equivalently,
\begin{equation}\label{vtheta-pm}
p_\pm=\frac{R\theta_\pm}{v_\pm},
\qquad
s_\pm=R\log v_\pm+\frac{R}{\gamma-1}\log\theta_\pm+\mathrm{const}.
\end{equation}

We consider the corresponding inviscid Euler system \eqref{E} with the Riemann data
\begin{equation}\label{Ei-v}
(v,u,\theta)(0,x)=
\begin{cases}
(v_-,u_-,\theta_-),&x<0,\\[1mm]
(v_+,u_+,\theta_+),&x>0.
\end{cases}
\end{equation}
We recall below the associated Riemann problem; see, for instance,
\cite{Riemann,Lax}. Let $V=(v,u,E)^{t}$, and then the Euler system \eqref{E} can be written as
\[
V_t+\mb A(V)V_x=0,
\]
where
\[
\mb A(V)=
\begin{pmatrix}
0 & -1 & 0 \\[1mm]
-\dfrac{p}{v}
& -\dfrac{(\gamma-1)u}{v}
& \dfrac{\gamma-1}{v} \\[2mm]
-\dfrac{pu}{v}
& p-\dfrac{(\gamma-1)u^2}{v}
& \dfrac{\gamma-1}{v}u
\end{pmatrix}.
\]
Since $v>0$ and $\theta>0$, the system is strictly
hyperbolic and has three real eigenvalues
\[
\lambda_1(v,\theta)
=
-\sqrt{\frac{\gamma p}{v}}<0,
\qquad
\lambda_2=0,
\qquad
\lambda_3(v,\theta)
=
\sqrt{\frac{\gamma p}{v}}>0.
\]

As discussed in \cite{Smoller}, the characteristic fields corresponding to $\lambda_1$ and $\lambda_3$ are genuinely nonlinear, whereas the field corresponding to $\lambda_2$ is linearly degenerate. In this paper, we consider a composite wave consisting of a rarefaction wave and a shock wave, corresponding to the solution of the
Riemann problem for the Euler system \eqref{E} with initial values
\eqref{Ei-v}. Without loss of generality, we take the rarefaction wave
and the shock wave to be a $1$-rarefaction wave and a $3$-shock wave,
respectively. For convenience, we first use the Riemann data
\begin{equation}\label{ELR}
(v,u,\theta)
=
\begin{cases}
(v_L,u_L,\theta_L), & x<0,\\[1mm]
(v_R,u_R,\theta_R), & x>0
\end{cases}
\end{equation}
to construct the elementary waves of system \eqref{E}, and then introduce the Riemann solution consisting of a 1-rarefaction wave and a 3-shock wave.
For the 1-rarefaction wave, given any state $(v_L,u_L,\theta_L)\in\mathbb{R}_+\times\mathbb{R}\times\mathbb{R}_+$, there exists the 1-rarefaction curve $R_1$ defined by
\[
R_1(v_L,u_L,\theta_L):=\left\{(v,u,\theta)\ \middle|\ v>v_L,\ s(v,\theta)=s(v_L,\theta_L),\ u=u_L-\int_{v_L}^{v}\lambda_1(s_L,v')\,dv'\right\},
\]
where the entropy $s$ is defined in \eqref{pes} and $\lambda_1(s,v)=-\sqrt{\frac{\gamma p(v,s)}{v}}$. Here and below, with a slight abuse of notation, $\lambda_1$ denotes the same quantity expressed either in the variables $(v,\theta)$ or in the variables $(s,v)$. Hence, if the initial values \eqref{ELR} satisfy $(v_R,u_R,\theta_R)\in R_1(v_L,u_L,\theta_L)$, then the solution to the Riemann problem
\eqref{E}-\eqref{ELR} is the self-similar 1-rarefaction wave
\[
(v^r,u^r,\theta^r)(t,x)
=
(v^r,u^r,\theta^r)\left(\frac{x}{t}\right),
\]
which is defined by
\begin{equation*}
\lambda_1\bigl(v^r(t,x),\theta^r(t,x)\bigr)
=
\begin{cases}
\lambda_1(v_L,\theta_L),
&
x<\lambda_1(v_L,\theta_L)t,
\\[2mm]
\displaystyle\frac{x}{t},
&
\lambda_1(v_L,\theta_L)t
\le x
\le
\lambda_1(v_R,\theta_R)t,
\\[2mm]
\lambda_1(v_R,\theta_R),
&
x>\lambda_1(v_R,\theta_R)t,
\end{cases}
\end{equation*}
together with
$$s(v^r(t,x),\theta^r(t,x))=s(v_L,\theta_L)=s(v_R,\theta_R),$$
and
$$
z_1(v^r(t,x),u^r(t,x))
=
z_1(v_L,u_L)
=
z_1(v_R,u_R),
$$
where $s(v,\theta)$ and $z_1(v,u)=u+\int_{v_L}^{v}\lambda_1(s_L,v') dv'$ are the 1-Riemann invariants to the Euler equations \eqref{E}.

For the 3-shock wave, given any state
$(v_L,u_L,\theta_L)\in\mathbb{R}_+\times\mathbb{R}\times\mathbb{R}_+$, the 3-shock curve $S_3(v_L,u_L,\theta_L)$ is determined by the Rankine-Hugoniot
conditions
\begin{equation}\label{RH-general}
\begin{cases}
-\sigma(v-v_L)-(u-u_L)=0,\\[1mm]
-\sigma(u-u_L)
+p(v,\theta)-p(v_L,\theta_L)=0,\\[1mm]
-\sigma(E-E_L)
+p(v,\theta)u-p(v_L,\theta_L)u_L=0,
\end{cases}
\end{equation}
where $E_L=e(v_L,\theta_L)+\frac{u_L^2}{2}$ and $\sigma=\sqrt{-\frac{p(v,\theta)-p(v_L,\theta_L)}{v-v_L}}>0$ satisfies the Lax entropy condition
$$
\lambda_3(v_R,\theta_R)
<
\sigma
<
\lambda_3(v_L,\theta_L).
$$
Thus, if the initial values \eqref{ELR} satisfy $(v_R,u_R,\theta_R)\in S_3(v_L,u_L,\theta_L)$, then the corresponding Riemann solution is given by
\[
(v^s,u^s,\theta^s)(t,x)
=
\begin{cases}
(v_L,u_L,\theta_L), & x<\sigma t,\\[1mm]
(v_R,u_R,\theta_R), & x>\sigma t.
\end{cases}
\]

We consider the Riemann solution to the Euler system \eqref{E}
with initial data \eqref{Ei-v}, consisting of a 1-rarefaction wave
and a 3-shock wave, that is,
$(v_+,u_+,\theta_+)\in S_3R_1(v_-,u_-,\theta_-)$. More precisely, it can be proven that there exists a unique intermediate state
$(v_m,u_m,\theta_m)\in\bbr_+\times\bbr\times\bbr_+$
such that $(v_m,u_m,\theta_m)\in R_1(v_-,u_-,\theta_-)$ and $ (v_+,u_+,\theta_+)\in S_3(v_m,u_m,\theta_m)$. Then, the Riemann
solution of the Euler system \eqref{E}-\eqref{Ei-v} can be expressed as
\begin{equation}\label{ICW}
(v,u,\theta)(t,x)
=
(v^r,u^r,\theta^r)(t,x)
+
(v^s,u^s,\theta^s)(t,x)
-
(v_m,u_m,\theta_m).
\end{equation}

Since the transformation $(v,u,\theta)\mapsto(p,u,s)$ is smooth and
one-to-one, the above elementary waves can equivalently be represented in the $(p,u,s)$-coordinates. For the 1-rarefaction wave, we define
\[
(p^r,u^r,s^r)(t,x)
:=
\bigl(p(v^r,\theta^r),u^r,s(v^r,\theta^r)\bigr)(t,x).
\]
Similarly, the 3-shock wave is represented by
\[
(p^s,u^s,s^s)(t,x)
:=
\bigl(p(v^s,\theta^s),u^s,s(v^s,\theta^s)\bigr)(t,x).
\]

Therefore, for the states
$(p_-,u_-,s_-),(p_+,u_+,s_+)\in\mathbb{R}_+\times\mathbb{R}\times\mathbb{R}$ connected by a 1-rarefaction wave and a 3-shock wave, let $(v_-,u_-,\theta_-)$ and $(v_+,u_+,\theta_+)$ be the corresponding states determined by \eqref{vtheta-pm}. Then, there exists an intermediate
state $(v_m,u_m,\theta_m)\in R_1(v_-,u_-,\theta_-)$ such that
$(v_+,u_+,\theta_+)\in S_3(v_m,u_m,\theta_m)$. We denote the corresponding
intermediate state in the $(p,u,s)$-coordinates by $(p_m,u_m,s_m)$, where
\begin{equation}\label{eq:sp-pm}
    p_m=\frac{R\theta_m}{v_m},
    \qquad
    s_m=R\log v_m+\frac{R}{\gamma-1}\log\theta_m+\mathrm{const}.
\end{equation}
Consequently, in the $(p,u,s)$-coordinates, the Riemann solution of
the Euler system \eqref{E}-\eqref{Euler-pus-data} can be equivalently expressed as
\begin{equation*}
(p,u,s)(t,x)
=(p^r,u^r,s^r)(t,x)+(p^s,u^s,s^s)(t,x)-(p_m,u_m,s_m),
\end{equation*}
where $(p^r,u^r,s^r)$ is the 1-rarefaction wave connecting
$(p_-,u_-,s_-)$ to $(p_m,u_m,s_m)$, and
$(p^s,u^s,s^s)$ is the 3-shock wave connecting
$(p_m,u_m,s_m)$ to $(p_+,u_+,s_+)$.

Subsequently, we turn to the viscous wave associated with the
Riemann solution \eqref{ICW}. 
The viscous 3-shock wave, also called the traveling wave, is $(v^S, u^S, \theta^S)(\xi)(\xi=x-\sigma t)$ connecting the intermediate state
$(v_m,u_m,\theta_m)$ to the right state $(v_+,u_+,\theta_+)$. It satisfies the following system of ODEs:
\begin{equation}\label{VS}
\begin{cases}
-\sigma (v^S)'-(u^S)'=0,\\
-\sigma (u^S)' +(p^S)'
=\displaystyle \mu\left(\frac{(u^S)'}{v^S}\right)',\\
-\sigma (E^S)' +(p^S u^S)'
=\displaystyle \mu\left(\frac{u^S(u^S)'}{v^S}\right)',\\
(v^S,u^S,\theta^S)(-\infty)
=(v_m,u_m,\theta_m),\qquad
(v^S,u^S,\theta^S)(+\infty)
=(v_+,u_+,\theta_+),
\end{cases}
\end{equation}
where $E^S=\frac{R}{\gamma-1}\theta^S+\frac{(u^S)^2}{2}$, $p^S=\frac{R\theta^S}{v^S}$, and the shock speed $\sigma
=\sqrt{-\frac{p(v_+,\theta_+)-p(v_m,\theta_m)}{v_+-v_m}}>0$ is determined by the Rankine-Hugoniot conditions and satisfies the Lax entropy condition
\begin{equation*}
\lambda_3(v_+,\theta_+)
<
\sigma
<
\lambda_3(v_m,\theta_m).
\end{equation*}
Since the main stability theorem is formulated in the $(p,u,s)$-coordinates, we also represent the viscous shock wave by $(p^S,u^S,s^S)$, where $p^S:=p(v^S,\theta^S)$ and $
s^S:=s(v^S,\theta^S)$. The properties of the viscous shock wave will be introduced in the next section. The composite wave $(\widetilde p,\widetilde u,\widetilde s)$ is then defined by
\begin{equation}\label{SW}
\begin{aligned}
(\widetilde p,\widetilde u,\widetilde s)(t,x)
:={}\Big(
&p^r(t,x)+p^S(x-\sigma t-{\mb X}(t))-p_m,\\
&
u^r(t,x)+u^S(x-\sigma t-{\mb X}(t))-u_m,\\
&
s^r(t,x)+s^S(x-\sigma t-{\mb X}(t))-s_m
\Big),
\end{aligned}
\end{equation}
which is composed of the 1-rarefaction wave connecting
$(p_-,u_-,s_-)$ and $(p_m,u_m,s_m)$ and the viscous 3-shock wave shifted by $\mb X(t)$ (to be defined in \eqref{X(t)}) connecting $(p_m,u_m,s_m)$ and $(p_+,u_+,s_+)$.

    \subsection{The main results}
\begin{theorem}\label{thm:main}
Let the far-field states $(p_\pm,u_\pm,s_\pm)\in\bbr_+\times\bbr\times\bbr$ be
connected by a 1-rarefaction wave and a 3-shock wave with intermediate state
$(p_m,u_m,s_m)$, as in Subsection \ref{subsec:Riemann}, and let
$(\widetilde p,\widetilde u,\widetilde s)$ be the composite wave defined in
\eqref{SW}. There exist constants $\delta_0,\varepsilon_0>0$ such that, if the
wave strengths satisfy
\[
|u_m-u_-|+|u_+-u_m|\leq \delta_0
\]
and the initial values $(p_0,u_0,s_0)$ satisfy
\begin{equation}\label{i-p}
\sum_{\pm}\left\|(p_0-p_\pm,u_0-u_\pm,s_0-s_\pm)\right\|_{L^2(\mathbb{R}_\pm)}
+\left\|(p_{0x},u_{0x},s_{0x})\right\|_{L^2(\mathbb{R})}
<\varepsilon_0,
\end{equation}
where $\mathbb{R}_-=(-\infty,0)$ and $\mathbb{R}_+=(0,+\infty)$, then the
Cauchy problem \eqref{VUE-1}-\eqref{p0-s0} admits a unique global solution
$(p,u,s)(t,x)$, and there exists an absolutely continuous time-dependent
shift $\mathbf{X}(t)$ such that the following statements hold.
\begin{enumerate}[label=\textup{(\roman*)},leftmargin=*]
\item \textup{(Uniform stability.)} For all $t>0$,
\begin{equation*}
\|(p- \widetilde p, u-\widetilde{u}, s-\widetilde s)(t,\cdot)\|_{L^2(\bbr)} + \|(p, u, s)_{x}(t,\cdot)\|_{L^2(\bbr)}\leq C(\eps_0+\delta_0^{\frac{1}{8}}).
\end{equation*}
\item \textup{(Convergence of the pressure and velocity.)}
\begin{equation}\label{con}
\lim_{t\to+\infty}\sup_{x\in\mathbb{R}}
\left|(p-\widetilde p,u-\widetilde u)(t,x)\right|=0.
\end{equation}
\item \textup{(Uniform bound for the entropy.)}
\begin{equation}\label{as}
\sup_{t>0}\sup_{x\in\mathbb{R}}
\left|s(t,x)-\widetilde s(t,x)\right|\leq C(\eps_0+\delta_0^{\frac 18}).
\end{equation}
\item \textup{(Shift.)}
\begin{equation}\label{shift-bound}
\lim_{t\to+\infty}|\dot{\mathbf{X}}(t)|=0.
\end{equation}
In particular, $\mathbf{X}(t)=o(t)$ as $t\to+\infty$, so that the shifted
viscous shock wave retains the original traveling wave profile and speed
time-asymptotically.
\end{enumerate}
Moreover, the solution has the regularity
\begin{align*}
(p-\widetilde p,u-\widetilde u,s-\widetilde s)&\in C((0,+\infty);H^1(\mathbb{R})),\\
u_{xx}(t,x)- u_{xx}^S(x-\s t-\mb X(t))&\in L^2(0,+\infty;L^2(\mathbb{R})).
\end{align*}
\end{theorem}

\begin{remark}
    The absence of heat conductivity leaves the entropy mode nondissipative, and system \eqref{VUE-1} therefore fails to satisfy the Shizuta-Kawashima condition. Nevertheless, Theorem \ref{thm:main} establishes the nonlinear stability of a composite wave consisting of a 1-rarefaction wave and a viscous 3-shock wave.
\end{remark}

\begin{remark}\label{rem:vtheta-asymptotics}
    Although neither the pressure equation nor the entropy equation contains a direct diffusion term, their dissipative properties are fundamentally different. The coupling between pressure and velocity allows us to find the wave structure and recover a dissipative estimate for the pressure, which is essential for the uniform convergence of $(p,u)$ in \eqref{con}. However, no analogous dissipative mechanism is available for the entropy, for which \eqref{as} provides only a uniform-in-time bound. This distinction is structural; see Liu-Zeng \cite{LZ99}. Indeed, system \eqref{VUE-1} admits stationary solutions of the form
    \[
        (p,u,s)(t,x)=(p_*,u_*,s_*(x)),
    \]
    where $p_*>0$ and $u_*$ are constants, while $s_*(x)$ may be nonconstant. Such solutions show that convergence to a prescribed entropy profile cannot, in general, be expected, even for small perturbations of a constant state. Hence the uniform bound \eqref{as} cannot, in general, be improved to convergence.

    The thermodynamic relations \eqref{pes}, together with the smallness of the wave strengths, also yield
    \begin{equation}\label{as-vtheta}
        \sup_{t>0}\sup_{x\in\bbr}
        \bigl|(v-\widetilde v,\theta-\widetilde\theta)(t,x)\bigr|
        \leq C\bigl(\eps_0+\delta_0^{\frac18}\bigr),
    \end{equation}
    where
    \begin{align*}
        \widetilde v(t,x)
        &:=v^r(t,x)
          +v^S\bigl(x-\sigma t-{\mb X}(t)\bigr)-v_m,\\
        \widetilde\theta(t,x)
        &:=\theta^r(t,x)
          +\theta^S\bigl(x-\sigma t-{\mb X}(t)\bigr)-\theta_m.
    \end{align*}
    Thus, the theorem establishes uniform stability of the specific volume and temperature. Since both variables depend on the entropy as well as the pressure, their convergence to the corresponding composite wave does not follow from \eqref{con}. This conclusion is consistent with the behavior near constant states described in \cite{LZ99}.
\end{remark}

\begin{remark}
The exponent $\frac18$ in Theorem \ref{thm:main} is not claimed to be
optimal. It originates from the error term of order $\delta_R^{\frac14}$
generated by the approximate rarefaction wave in the (squared) basic energy
estimate; see Lemmas \ref{lem-sharp2} and \ref{Re-en}.
\end{remark}

\begin{corollary}[Stability of a single viscous shock]
Let $(p_+,u_+,s_+)\in\bbr_+\times\bbr\times\bbr$ be connected to
$(p_-,u_-,s_-)$ by a 3-shock wave with $|u_+-u_-|\le\delta_0$, and let the
initial values satisfy \eqref{i-p}. Then all the conclusions of Theorem
\ref{thm:main} hold with
$(\widetilde p,\widetilde u,\widetilde s)(t,x)
=(p^S,u^S,s^S)(x-\sigma t-\mathbf X(t))$.
\end{corollary}
\begin{proof}
This is the case of vanishing rarefaction strength in Theorem
\ref{thm:main}: the intermediate state equals the left state, the
rarefaction wave is constant, and the composite wave reduces to the shifted
viscous shock wave.
\end{proof}

\begin{remark}
    The approach developed here also applies to a composite wave consisting of a viscous 1-shock and a viscous 3-shock for the same zero heat-conductivity system. Since that case requires two independent shifts and additional interaction estimates, it is treated in a separate paper.
\end{remark}

\begin{remark}
    The stability of wave patterns involving a contact discontinuity remains open for the zero heat-conductivity system considered here. A central difficulty is that the contact mode lacks thermal diffusion, so the classical construction of a viscous contact wave and the associated stability analysis cannot be applied directly. From the viewpoint of the present method, an additional obstruction is that a contact wave is stationary in Lagrangian coordinates and therefore sits on the entropy characteristics, so that the mechanism behind the weighted estimates of Subsection \ref{subsec:weighted}, namely that the wave moves relative to the entropy characteristics, is not available.
    This issue will be left for our future research.
\end{remark}

\subsection{Main difficulties and key ideas of the proof}
To establish the nonlinear stability of the composite wave consisting of a rarefaction wave and a viscous shock wave for system \eqref{VUE-1}, we employ the $a$-contraction method with shifts developed in \cite{KVW25}, which allows viscous shock waves and rarefaction waves to be treated within a common energy framework. The absence of heat conductivity, however, creates additional difficulties in both the relative entropy estimate and the first-order derivative estimates. In particular, the entropy mode remains nondissipative, so the missing entropy dissipation cannot be recovered through coupling with the velocity. We address below the main difficulties and the mechanisms used to overcome them.

\medskip
\noindent\textbf{(i) Relative entropy estimates without temperature
dissipation.}

In the relative entropy framework of \cite{KVW25}, heat conduction provides space-time control of the temperature derivative, which is used to estimate nonlinear error terms. This control is unavailable for system \eqref{VUE-1}. To handle the terms involving the temperature, we exploit the expansive property of the rarefaction wave together with its temporal decay. More precisely, the expansion yields favorable contributions to the energy estimate, while the decay of the wave derivatives allows the remaining terms to be bounded by coefficients that are integrable in time. Combining these bounds with the dissipation supplied by viscosity and the wave profiles, we close the basic relative entropy estimate by Gronwall's inequality.

\medskip
\noindent\textbf{(ii) Pressure derivative estimates through the effective viscous flux.}

The first-derivative estimates in \cite{KVW25} also rely
essentially on temperature dissipation. For the present
system, we recover the required pressure estimate by exploiting the coupling between the momentum equation \eqref{p-u-s}$_2$, which contains the pressure gradient, and the pressure equation
\eqref{p-u-s}$_1$. This observation motivates us to work
with the $(p,u,s)$ variables and introduce an effective
viscous flux. Indeed, setting $F:=\mu\frac{u_x}{v}-p$, a direct
computation based on \eqref{p-u-s}$_1$, \eqref{p-u-s}$_2$ and $v_t=u_x$ shows
that $u_t=F_x$ and that $F$ satisfies the exact identity
\[
    F_t-\frac{\mu}{v}F_{xx}+\frac{\gamma u_x}{v}F=0,
\]
that is, a parabolic equation which involves neither the entropy nor
derivatives of the pressure. To illustrate the underlying mechanism, consider the linearization about the wave profile with specific volume $\bar v>0$, velocity $\bar u$, and pressure $\bar p>0$ defined in \eqref{per} and \eqref{per-v-theta}. Let $W$, $\psi$, and $\mathcal F$ denote the pressure, velocity, and effective viscous flux perturbations. Then the linearized equations for velocity and the effective viscous flux perturbations exhibit, to leading order, the following coupling structure:
\[
    \psi_t=\mathcal F_x,\qquad
    \mathcal F_t-\frac{\mu}{\bar v}\mathcal F_{xx}
    +\frac{\gamma\bar u_x}{\bar v}\mathcal F
    =\frac{\gamma\bar p}{\bar v}\psi_x.
\]
The associated energy estimate provides space-time estimate
of $\mathcal F_x$. Combining this estimate with the velocity
dissipation yields the required space-time estimate for $W_x$.

The same mechanism is reflected in the damped wave structure arising from the coupling between \eqref{p-u-s}$_1$ and \eqref{p-u-s}$_2$. We can also obtain the linearized wave equation
\[
    W_{tt}-\frac{\gamma\bar p}{\bar v}W_{xx}
    -\frac{\mu}{\bar v}W_{txx}=0.
\]
This structure shows how velocity dissipation is transferred
to the pressure through their coupling.

In comparison with the rarefaction-wave stability result of Hu-Yin
\cite{HY18}, which requires small $H^2$ initial perturbations, the effective viscous flux that reflects the pressure wave structure in our argument allows us to close the energy estimates at the $H^1$ level, requiring only a small $H^1$ initial perturbation.

\medskip
\noindent\textbf{(iii) Weighted estimates for the entropy derivative.}

Although the system of an ideal gas without viscosity lacks direct velocity dissipation in \cite{GLWX}, it still satisfies the Shizuta-Kawashima condition. By exploiting the wave structure, the dissipation of the temperature can be transferred to the density and velocity. For system \eqref{VUE-1}, however, the Shizuta-Kawashima condition fails. Therefore, the entropy equation \eqref{p-u-s}$_3$ contains neither dissipation nor damping, and no dissipative wave structure is available to transfer velocity dissipation to the entropy. This structural degeneracy accounts for the distinction between the asymptotic convergence of pressure and velocity and the uniform stability of entropy. It also prevents us from deducing convergence of the temperature and density to their prescribed wave profiles, consistently with the behavior near constant states described by Liu-Zeng \cite{LZ99}.

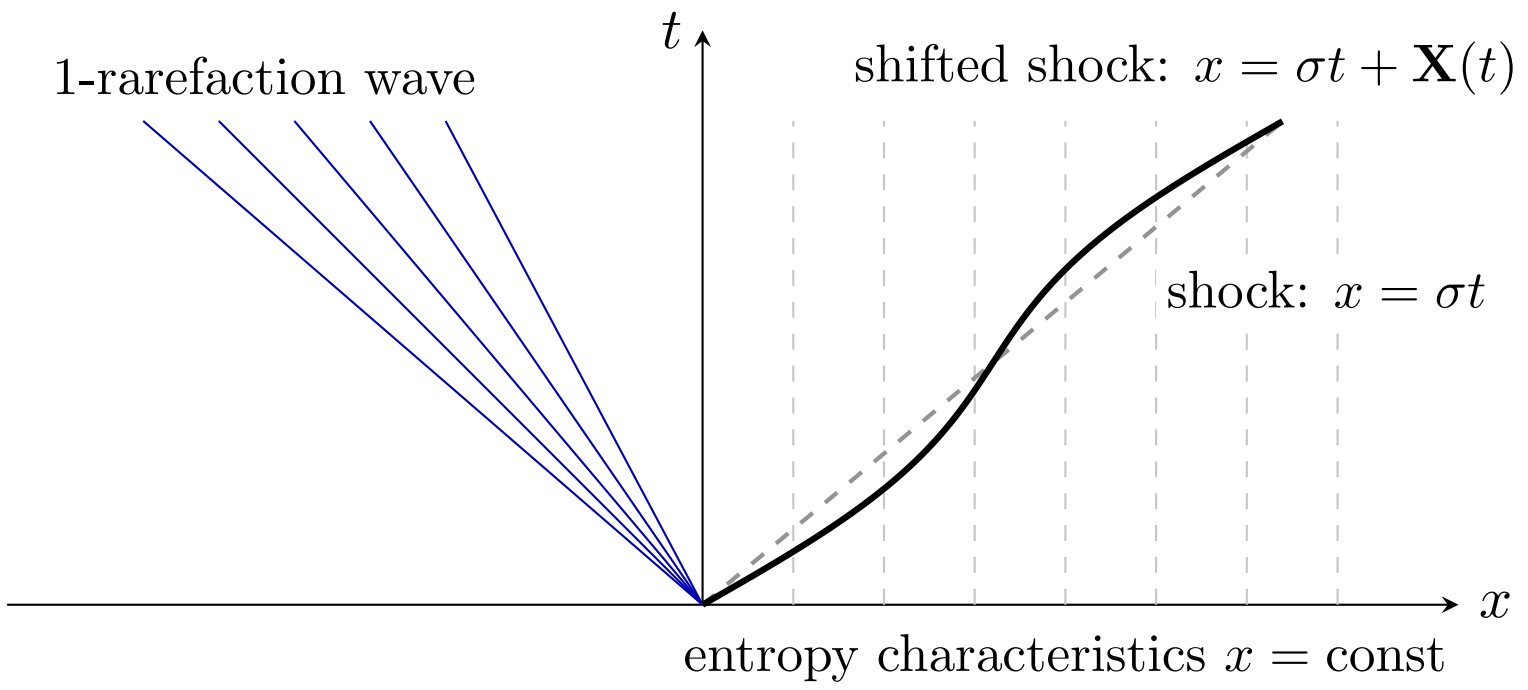
\begin{figure}[ht]
    \centering
    \begin{tikzpicture}[>=stealth]
        \draw[->] (-4.6,0) -- (5.0,0) node[right] {$x$};
        \draw[->] (0,0) -- (0,3.8) node[left] {$t$};
        \foreach \a in {-3.7,-3.2,-2.7,-2.2,-1.7}
        \draw[blue!70!black] (0,0) -- (\a,3.2);
        \node at (-2.9,3.5) {\small $1$-rarefaction wave};

        \foreach \x in {0.6,1.2,1.8,2.4,3.0,3.6,4.2}
            \draw[dashed,gray!45] (\x,0) -- (\x,3.2);

        \draw[thick,dashed,gray!85]
            (0,0) -- (3.84,3.2)
            node[pos=0.75,below right=3pt,fill=white,inner sep=2pt,text=black
            ] {\small shock: $x=\sigma t$};
            
        \draw[very thick,black]
            plot[domain=0:3.2,samples=160,variable=\t]
            ({1.2*\t + 0.28*sin(112.5*\t)},{\t});
        \node at (3.2,3.55)
        {\small shifted shock: $x=\sigma t+\mathbf X(t)$};
        \node at (2.4,-0.35)
        {\small entropy characteristics $x=\mathrm{const}$};
    \end{tikzpicture}
\caption{Schematic picture in Lagrangian coordinates. The entropy is transported along the vertical lines $x=\mathrm{const}$ (dashed, drawn only for $x>0$), which the shifted viscous shock crosses transversally. Each entropy characteristic therefore meets the shock layer only during a time interval of length $O(\delta_S^{-1})$, which yields \eqref{key-mech}.}
\label{fig:mechanism}
\end{figure}

The key observation is that the energy argument does not require an unweighted space-time estimate for the entropy derivative. Instead, the relevant terms contain weights arising from the shock profile or decaying factors associated with the rarefaction wave. For the terms weighted by the
shock profile, we integrate the differentiated entropy equation along its characteristics and then apply the Fubini-Tonelli theorem to change the order of integration. Since the shock propagates relative to the entropy
characteristics, the localized shock weight has a uniformly bounded time integral along these curves. 
More precisely, since the entropy equation \eqref{p-u-s}$_3$ contains no spatial derivative of $s$, the entropy characteristics are the vertical lines $x=\mathrm{const}$ in Lagrangian coordinates, whereas the shifted shock travels with speed $\sigma+\dot{\mathbf X}(t)\ge\frac{\sigma}{2}>0$. Hence, for every fixed $x\in\bbr$,
\begin{equation}\label{key-mech}
    \int_0^{\infty}\bigl|(v^S)'\bigl(x-\sigma t-\mathbf X(t)\bigr)\bigr|\,dt
    \le\frac{2}{\sigma}\int_{\bbr}|(v^S)'(y)|\,dy
    \le C\delta_S;
\end{equation}
see Figure \ref{fig:mechanism} and Lemma \ref{trans}.
This yields the required weighted space-time estimate for the entropy derivative. For the rarefaction-weighted terms, we use its expansive property and temporal decay to obtain coefficients that are integrable in time. Gronwall's inequality then gives bounds that are uniform in time. These weighted estimates, combined with the pressure
and velocity estimates, close the global energy estimates despite the absence of entropy dissipation.

The remainder of this paper is organized as follows.
In Section 2, we collect the properties of the approximate
rarefaction waves and viscous shock waves. In Section 3,
we introduce the time-dependent shift, state the main
{\it a priori} estimate, and show how it implies global
existence and nonlinear stability. Sections 4 and 5
provide the proof of Proposition \ref{prop:priori}. In Section 4,
we derive the basic $L^{2}$-estimates by the relative
entropy method and the $a$-contraction method with shifts. In Section 5, we establish the first-order derivative estimates and the weighted space-time estimates for the entropy derivative, which together close the global $H^{1}$ {\it a priori} estimate. Section 4 follows the framework of \cite{KVW25}, with the modifications required by the absence of temperature dissipation, whereas Section 5 contains the main new ingredients of this paper; a roadmap of the argument is given at the beginning of Subsection \ref{subsec:first}.

\noindent$\bullet$ {\bf Notation.} We use the following notation throughout the paper. \\
    1. $C,C_0,C_1,\ldots$ denote generic positive constants, which may vary from line to line and are independent of the small parameters $\delta_0,\eps_0,\varepsilon_1,\delta_S,\delta_R$, $\lambda$ (see \eqref{weight}), and the time $T$. For simplicity, such constants will be used throughout the paper without further specification. \\
	  2. For any function $f : \bbr_+\times \bbr\to \bbr$ and any time-dependent shift $\mb X(t)$, 
	\[
	  f^{\pm \mb X}(t, \xi):=f(t,\xi\pm \mb X(t)).
	  \]
	  To simplify the notation, we drop the explicit arguments of the waves whenever the context is clear: for instance,
	  \begin{align*}
	    u^S:=u^S(\xi),\qquad\quad\quad  &(u^S)^{-\mb X}:=u^S(\xi- \mb X(t)),\\
        u^R:=u^R(t,\xi+\s t),\quad  &(u^R)^{\mb X}:=u^R(t,\xi+\s t + \mb X(t)).
	  \end{align*}
	  3. For functional spaces, $H^s(\mathbb{R})$ denotes the $s$-th order Sobolev space with the norm  
	    \[  
	    \|f\|_{H^s(\mathbb{R})}^2 := \sum_{j=0}^{s} \|\partial_x^j f\|_{L^2(\mathbb{R})}^2.  
	    \]  

\section{Preliminaries}\label{sec-preli}

In this section, we collect the basic properties of the approximate
1-rarefaction wave and the viscous 3-shock wave that will be used in
the stability analysis. 

\subsection{Approximate rarefaction waves}
Following \cite{MN92}, we construct a smooth approximation of the
1-rarefaction wave by means of a smooth solution to the Burgers equation
\begin{equation}\label{ABE}
\begin{cases}
\displaystyle w_t+ww_x=0,\\[1mm]
\displaystyle
w(0,x)=w_0(x)
=\frac{w_m+w_-}{2}
+K\frac{w_m-w_-}{2}\int_0^{\beta x}(1+y^2)^{-2}\,dy,
\end{cases}
\end{equation}
where $K=\left(\int_0^{+\infty}(1+y^2)^{-2}\,dy\right)^{-1}$ and $\beta>0$ is a parameter to be specified later.

We define the smooth approximate 1-rarefaction wave
$(v^R,u^R,\theta^R)(t,x)$ associated with the self-similar rarefaction
wave $(v^r,u^r,\theta^r)(\frac xt)$ by
\begin{align}
\begin{aligned}\label{AR}
&\lambda_{1-}:=\lambda_1(v_-,\theta_-)=w_-,
\qquad
\lambda_{1m}:=\lambda_1(v_m,\theta_m)=w_m,\\
&\lambda_1(v^R,\theta^R)(t,x)=w(1+t,x),\\
&z_1(v^R,u^R)(t,x)
=z_1(v_-,u_-)=z_1(v_m,u_m),\\
&s(v^R,\theta^R)(t,x)
=s(v_-,\theta_-)=s(v_m,\theta_m),
\end{aligned}
\end{align}
where $w(t,x)$ denotes the smooth solution to \eqref{ABE}.

A direct calculation shows that $(v^R,u^R,\theta^R)$ satisfies the
Euler system, and the basic properties of the approximate rarefaction wave are summarized in the following lemma; see \cite{MN92}.

\begin{lemma}\label{lem:rare-v}
The smooth approximate 1-rarefaction wave
$(v^R,u^R,\theta^R)(t,x)$ defined in \eqref{AR} satisfies the following
properties. Let $\delta_R:=|v_m-v_-|\sim |u_m-u_-|\sim |\theta_m-\theta_-|$ denote the rarefaction wave strength, and set $\beta=\deltar$. Then, there exists a constant $C>0$ such that:

\noindent
(1)~It holds that $u^R_x=\frac{2v^R}{\gamma+1}w_x>0$,
$v^R_x=\frac{v^R}{\sqrt{R\gamma\theta^R}}u^R_x>0$ and
$\theta^R_x=-\frac{(\gamma-1)\theta^R}{v^R}v^R_x<0$,
for all $x\in\bbr$ and $t\geq0$.

\noindent
(2)~For all $t\geq0$ and $p\in[1,+\infty]$, it holds that
\begin{align*}
&\|(v^R_x,u^R_x,\theta^R_x)\|_{L^p}
\leq C\min\left\{
\delta_R^{2-\frac1p},
\,
\delta_R^{\frac1p}(1+t)^{-1+\frac1p}
\right\},\\
&\|(v^R_{xx},u^R_{xx},\theta^R_{xx})\|_{L^p}
\leq C\min\left\{
\delta_R^{3-\frac1p},
\,
\delta_R^{\frac12-\frac1{2p}}
(1+t)^{-\frac54+\frac1{4p}}
\right\},\\
&|u^R_{xx}|\leq C\deltar|u^R_x|,
\qquad
|\theta^R_{xx}|\leq C\deltar|\theta^R_x|,
\qquad \forall\,x\in\bbr.
\end{align*}

\noindent
(3)~For $x\geq\lambda_{1m}(1+t)$ and $t\geq0$, it holds that
\begin{align*}
&|(v^R,u^R,\theta^R)(t,x)-(v_m,u_m,\theta_m)|
\leq C\delta_R
\left(1+\delta_R|x-\lambda_{1m}(1+t)|\right)^{-3},\\
&|(v^R_x,u^R_x,\theta^R_x)(t,x)|
\leq C\delta_R^2
\left(1+\delta_R|x-\lambda_{1m}(1+t)|\right)^{-4}.
\end{align*}

\noindent
(4)~For $x\leq\lambda_{1-}(1+t)$ and $t\geq0$, it holds that
\begin{align*}
&|(v^R,u^R,\theta^R)(t,x)-(v_-,u_-,\theta_-)|
\leq C\delta_R
\left(1+\delta_R|x-\lambda_{1-}(1+t)|\right)^{-3},\\
&|(v^R_x,u^R_x,\theta^R_x)(t,x)|
\leq C\delta_R^2
\left(1+\delta_R|x-\lambda_{1-}(1+t)|\right)^{-4}.
\end{align*}
\end{lemma}

We define the smooth approximate 1-rarefaction wave in the
$(p,u,s)$-coordinates by
\begin{equation}\label{eq:pus-r}
(p^R,u^R,s^R)(t,x):=
\bigl(p(v^R,\theta^R),u^R,s(v^R,\theta^R)\bigr)(t,x).
\end{equation}

\begin{lemma}\label{lem:rare-p}
Let $(p^R,u^R,s^R)(t,x)$ be the smooth approximate
$1$-rarefaction wave defined in \eqref{eq:pus-r}.
Then there exists a constant $C>0$ such that the following estimates hold.

\noindent
(1) For all $t>0$,
\[
\left\|
(p^R,u^R,s^R)(t,\cdot)
-
(p^r,u^r,s^r)\left(\frac{\cdot}{t}\right)
\right\|_{L^2}
\leq
C\delta_R^{\frac12}
(1+\delta_R^2 t)^{-\frac12}.
\]

\noindent
(2) For all $t>0$,
\[
\sup_{x\in\bbr}
\left|
(p^R,u^R,s^R)(t,x)
-
(p^r,u^r,s^r)\left(\frac{x}{t}\right)
\right|
\leq
C\delta_R
(1+\delta_R^2 t)^{-\frac34}.
\]
\end{lemma}
Lemma \ref{lem:rare-p} can be derived from the characteristic representation of the Burgers solution $w$ in \eqref{ABE}, as in \cite{MN92}. Since the map $(v,\theta)\mapsto(p,s)$ is a smooth diffeomorphism, the same estimates hold with $(p^R,s^R)-(p^r,s^r)$ replaced by $(v^R,\theta^R)-(v^r,\theta^r)$; we shall refer to this as the $(v,\theta)$-form of Lemma \ref{lem:rare-p}.

\subsection{Viscous shock waves}

For the viscous 3-shock wave, integrating \eqref{VS} from the two end
states and using $E^S=\frac{R}{\gamma-1}\theta^S+\frac{(u^S)^2}{2}$,
we obtain
\begin{equation}\label{shock-v-u-theta-1}
\begin{cases}
\displaystyle
-\mu\sigma\frac{(v^S)'}{v^S}
=(p^S-p_+)+\sigma^2(v^S-v_+)
=(p^S-p_m)+\sigma^2(v^S-v_m),\\[3mm]
0=\displaystyle
\frac{R}{\gamma-1}(\theta^S-\theta_+)
+p_+(v^S-v_+)-\frac12\sigma^2(v^S-v_+)^2\\
\quad=\displaystyle\frac{R}{\gamma-1}(\theta^S-\theta_m)
+p_m(v^S-v_m)-\frac12\sigma^2(v^S-v_m)^2.
\end{cases}
\end{equation}
Thus, the existence of the viscous shock profile
$(v^S,u^S,\theta^S)(\xi)$ is equivalent to that of a solution to
\eqref{shock-v-u-theta-1}. Once $(v^S,u^S,\theta^S)$ is determined,
we define
\[
p^S(\xi):=p(v^S(\xi),\theta^S(\xi)),
\qquad
s^S(\xi):=s(v^S(\xi),\theta^S(\xi)),
\]
so that $(p^S,u^S,s^S)(\xi)$ is the corresponding viscous shock wave connecting $(p_m,u_m,s_m)$ and $(p_+,u_+,s_+)$. Based on \eqref{shock-v-u-theta-1} and the arguments in \cite{FM, KVW25, MN85}, we summarize below the basic properties of the viscous 3-shock waves $(v^S,u^S,\theta^S)$ and $(p^S,u^S,s^S)$.

\begin{lemma}\label{lem:shock-v}
For any state $(v_+,u_+,\theta_+)$, there exists a constant $C>0$ such
that the following holds. For any state $(v_m,u_m,\theta_m)$ such that
$(v_+,u_+,\theta_+)\in S_3(v_m,u_m,\theta_m)$, there exists a unique
(up to a constant shift) solution $(v^S,u^S,\theta^S)(\xi)$ to \eqref{shock-v-u-theta-1}. Let $\delta_S$ denote the shock strength as $\delta_S:=|v_+-v_m|\sim |u_+-u_m|\sim |\theta_+-\theta_m|$. Then, the following estimates hold 
\begin{align}
\begin{aligned}\label{shock-base}
&u^S_\xi<0,\qquad
v^S_\xi>0,\qquad
\theta^S_\xi<0,
\qquad \forall\,\xi\in\bbr,\\
&|(v^S(\xi)-v_m,u^S(\xi)-u_m,\theta^S(\xi)-\theta_m)|
\leq
C\delta_S e^{-C\delta_S|\xi|},
\qquad \xi<0,\\[1mm]
&|(v^S(\xi)-v_+,u^S(\xi)-u_+,\theta^S(\xi)-\theta_+)|
\leq
C\delta_S e^{-C\delta_S|\xi|},
\qquad \xi>0,\\[1mm]
&|(v^S_\xi,u^S_\xi,\theta^S_\xi)|
\leq
C\delta_S^2e^{-C\delta_S|\xi|},
\qquad \forall\,\xi\in\bbr,\\[1mm]
&|(v^S_{\xi\xi},u^S_{\xi\xi},\theta^S_{\xi\xi})|
\leq
C\delta_S|(v^S_\xi,u^S_\xi,\theta^S_\xi)|,
\qquad \forall\,\xi\in\bbr.
\end{aligned}
\end{align}
In particular,
$|v^S_\xi|\sim|u^S_\xi|\sim|\theta^S_\xi|$ for all $\xi\in\bbr$.
More explicitly,
\begin{equation}\label{shock-vu}
\left|(u^S)_\xi+\sigma_m(v^S)_\xi\right|
\leq
C\delta_S|(v^S)_\xi|,
\qquad
\forall\,\xi\in\bbr,
\end{equation}
and
\begin{equation}\label{theta-s}
\left|
(\theta^S)_\xi
+\frac{(\gamma-1)p_m}{R}(v^S)_\xi
\right|
\leq
C\delta_S|(v^S)_\xi|,
\qquad
\forall\,\xi\in\bbr,
\end{equation}
where $\sigma_m:=\sqrt{\frac{\gamma p_m}{v_m}}=\frac{\sqrt{\gamma R\theta_m}}{v_m}$, $p_m,s_m$ satisfy \eqref{eq:sp-pm}, and
\begin{equation}\label{sigma}
|\sigma-\sigma_m|\leq C\delta_S.
\end{equation}
\end{lemma}

We then prove the following refined relation for the
viscous shock profile, which will be used in the subsequent energy
estimates.

\begin{lemma}\label{sharp-D}
For the viscous shock profile $(v^S,u^S,\theta^S)(\xi)$ connecting
$(v_m,u_m,\theta_m)$ and $(v_+,u_+,\theta_+)$, it holds that
\[
\left|
\frac{p^S-p_+}{v^S-v_+}
-\frac{p^S-p_m}{v^S-v_m}
-\sigma_m\alpha_m(v_+-v_m)
\right|
\leq C\delta_S^2,
\]
where $\alpha_m
:=
\frac{\gamma(\gamma+1)p_m}
{2(v_m)^2\sigma_m}$.
\end{lemma}
\begin{proof}
By $\eqref{shock-v-u-theta-1}_2$, together with $R\theta^S=p^S v^S$ and $R\theta_m=p_m v_m$, we have
\[
p^S v^S-p_m v_m
=-(\gamma-1)p_m(v^S-v_m)
+\frac{\gamma-1}{2}\sigma^2(v^S-v_m)^2,
\]
which means 
\[
(p^S-p_m)v^S
=-\gamma p_m(v^S-v_m)
+\frac{\gamma-1}{2}\sigma^2(v^S-v_m)^2,
\]
Consequently,
\[
\frac{p^S-p_m}{v^S-v_m}
=\frac{1}{v^S}
\left[
-\gamma p_m
+\frac{\gamma-1}{2}\sigma^2(v^S-v_m)
\right].
\]
Similarly, we can also get
\[
\frac{p^S-p_+}{v^S-v_+}
=\frac{1}{v^S}
\left[
-\gamma p_+
+\frac{\gamma-1}{2}\sigma^2(v^S-v_+)
\right].
\]

On the other hand, the Rankine-Hugoniot relations imply
\[
p_+-p_m=-\sigma^2(v_+-v_m).
\]
Subtracting the above two identities and using this relation,
we find
\[
\begin{aligned}
\frac{p^S-p_+}{v^S-v_+}
-\frac{p^S-p_m}{v^S-v_m}
&=\frac{1}{v^S}
\left[
\gamma(p_m-p_+)
-\frac{\gamma-1}{2}\sigma^2(v_+-v_m)
\right]\\
&=\frac{\gamma+1}{2}
\frac{\sigma^2}{v^S}(v_+-v_m).
\end{aligned}
\]

Then, using \eqref{shock-base} and \eqref{sigma}, we obtain
\[
\begin{aligned}
\left|
\frac{\sigma^2}{v^S}
-\frac{\gamma p_m}{v_m^2}
\right|
&\le
\frac{1}{v^S}
\left|\sigma^2-\frac{\gamma p_m}{v_m}\right|
+\frac{\gamma p_m}{v_m}
\left|\frac{1}{v^S}-\frac{1}{v_m}\right|\\
&\le C\left(
|v_+-v_m|+|v^S-v_m|
\right)\\
&\le C|v_+-v_m|.
\end{aligned}
\]
Therefore,
\begin{align*}
&\left|
\frac{p^S-p_+}{v^S-v_+}
-\frac{p^S-p_m}{v^S-v_m}
-\sigma_m\alpha_m(v_+-v_m)
\right|\\
=&\frac{\gamma+1}{2}|v_+-v_m|
\left|
\frac{\sigma^2}{v^S}
-\frac{\gamma p_m}{v_m^2}
\right|
\le  C|v_+-v_m|^2
\le C\delta_S^2,
\end{align*}
where $\alpha_m=
\frac{\gamma(\gamma+1)p_m}
{2(v_m)^2\sigma_m}$.
\end{proof}
By \eqref{pes}, \eqref{eq:sp-pm}, and Lemma \ref{lem:shock-v}, we further
obtain the following estimates for $(p^S,s^S)$.

\begin{lemma}\label{lem:shock-p}
For the viscous 3-shock wave $(p^S,u^S,s^S)(\xi)$ formed by
$(v^S,u^S,\theta^S)(\xi)$, let $\delta_S$ denote the shock strength
defined in Lemma \ref{lem:shock-v}. Then, for $\delta_S>0$ sufficiently
small,
\[
|p_+-p_m|\sim\delta_S,
\qquad
|s_+-s_m|\sim\delta_S^3,
\]
and it holds that
\begin{equation*}
\begin{alignedat}{3}
|p^S_\xi|
&\leq C\delta_S^2 e^{-C\delta_S|\xi|},
&\qquad
|s^S_\xi|
&\leq C\delta_S^4 e^{-C\delta_S|\xi|},
&\qquad \forall\,\xi\in\bbr,
\\
|p^S_{\xi\xi}|
&\leq C\delta_S |p^S_\xi|,
&
|s^S_{\xi\xi}|
&\leq C\delta_S |s^S_\xi|,
&\qquad \forall\,\xi\in\bbr.
\end{alignedat}
\end{equation*}
In particular,
\[
|p^S_\xi|
\sim |v^S_\xi|
\sim |u^S_\xi|,
\qquad
|s^S_\xi|
\sim |u^S_\xi|^2,
\qquad
\forall\,\xi\in\bbr.
\]
\end{lemma}
  
  \section{Proof of the main results}

   For simplicity, we rewrite the system \eqref{p-u-s} in the following system, based on the change of variable associated to the speed of propagation of the shock $(t, x) \mapsto (t, \xi = x - \sigma t)$:            
         \begin{equation}\label{p-u-s-1}
        	\begin{cases}
        		p_t -\sigma p_\x+\frac{\gamma R \theta}{v^2} u_\x = \mu(\gamma-1)\frac{u^2_\x}{v^2},\\[1mm]
        		u_t - \sigma u_\xi + p(v, \theta)_\xi = \left( \mu \frac{u_\xi}{v} \right)_\xi, \\[1mm]
        		s_t  -\sigma s_\x = \mu\frac{u^2_\x}{v\theta}.
        	\end{cases}
        \end{equation}
In these coordinates, the corresponding composite wave is given by       
\begin{equation}\label{per}
\begin{aligned}
(\bar p,\bar u,\bar s)(t,\xi)
:={}\Big(
&p^R(t,\xi+\sigma t)+p^S(\xi-\mathbf X(t))-p_m,\\
&
u^R(t,\xi+\sigma t)+u^S(\xi-\mathbf X(t))-u_m,\\
&
s^R(t,\xi+\sigma t)+s^S(\xi-\mathbf X(t))-s_m
\Big),
\end{aligned}
\end{equation}
where $\mathbf X(t)$ will be defined in \eqref{X(t)}. Then, the composite wave $(\bar{p}, \bar{u}, \bar{s})(t, \xi)$ satisfies the system   
         \begin{equation}\label{shock-p-u-s}
        	\begin{cases}
        		\bar p_t -\sigma \bar p_\x+ \dot{\mathbf{X}}(t) (p^S)_\xi^{-\mathbf{X}} +\frac{\gamma \bar p}{\bar v} \bar u_\x = \mu(\gamma-1)\frac{\bar u^2_\x}{\bar v^2}+R_1,\\[1mm]
        		\bar{u}_t - \sigma \bar{u}_\xi + \dot{\mathbf{X}}(t) (u^S)_\xi^{-\mathbf{X}} + \bar{p}_\xi = \mu \left( \frac{\bar{u}_\xi}{\bar{v}} \right)_\x+R_2, \\[1mm]
        		\bar s_t  -\sigma \bar s_\x + \dot{\mathbf{X}}(t) (s^S)_\xi^{-\mathbf{X}} = \mu\frac{\bar u^2_\x}{\bar v\bar\theta}+R_3,
        	\end{cases}
        \end{equation}
        where the error terms satisfy
        	\begin{align*}
        		&R_j=R_j^I+R_j^R, \quad  j=1,2,3,
        	\end{align*}
        with the wave interaction terms 
        \begin{equation}\label{RI}
        \begin{aligned}
        	& R_1^I=\gamma\left(\frac{\bar p}{\bar v}\bar u_\xi-\frac{p^R}{v^R}u_\xi^R-\frac{(p^S)^{-\mathbf X}(u^S)_{\xi}^{-\mathbf X}}
            {(v^S)^{-\mathbf X}}\right)-(\gamma-1)\mu\left(\frac{\bar u_\xi^2}{\bar v^2}-\frac{(u_\xi^R)^2}{(v^R)^2}-\frac{\big((u^S)_{\xi}^{-\mathbf X}\big)^2}{\big((v^S)^{-\mathbf X}\big)^2}\right),\\
        	&R_2^I=-\mu\left(\frac{\bar u_\x}{\bar v}-\frac{u^R_\x}{v^R}-\frac{(u^S)^{-\mb X}_\x}{(v^S)^{-\mb X}}\right)_\x,
        	\qquad R_3^I=-\mu\left(\frac{\bar u^2_\x}{\bar v\bar \theta}-\frac{(u^R_\x)^2}{v^R \theta^R}-\frac{((u^S)^{-\mb X}_\x)^2}{(v^S)^{-\mb X}(\theta^S)^{-\mb X}}\right),
        \end{aligned}
        \end{equation}
        and the error terms due to the inviscid rarefaction wave
        \begin{equation*}
        	\begin{aligned}
        		&R^R_1=-(\gamma-1)\mu \frac{(u^R_\x)^2}{(v^R)^2},\qquad R^R_2=	-\mu \left(\frac{u^R_\x}{v^R}\right)_\x, \qquad R^R_3=-\mu\frac{(u^R_\x)^2}{v^R \theta^R}.
        	\end{aligned}
        \end{equation*} 	

       \subsection{Construction of the shift} 
       For notational convenience, we also denote
      \begin{equation}\label{per-v-theta}
          \begin{aligned}
          (\bar v,\bar\theta)(t,\xi)
          :={}\Big(
          &v^R(t,\xi+\sigma t)+v^S(\xi-\mathbf X(t))-v_m,\\
          &\theta^R(t,\xi+\sigma t)+\theta^S(\xi-\mathbf X(t))-\theta_m
          \Big).
         \end{aligned}
         \end{equation}
        For the main result, we use the continuation argument based on the {\it a priori} estimates of Proposition \ref{prop:priori}. These estimates depend on the shift function, and for this reason, we first give its definition. The definition depends on the weight function $a : \mathbb{R} \to \mathbb{R}$ defined in $\eqref{weight}$. For now, we will only use the fact that $\|a\|_{C^1(\mathbb{R})} \leq 2$. We then define the shift $\mathbf{X}(t)$ as a solution to the ODE: 
        for any given $x_0 \in \mathbb{R}$,
        
        \begin{equation}\label{X(t)} 
        	\begin{cases}
        		\dot{\mathbf{X}}(t) = -\frac{M}{\delta_S} \left\{ \displaystyle\int_{\mathbb{R}} a(\xi - \mathbf{X}) \left[  u_\xi^S(\xi - \mathbf{X})(u - \bar{u}) + \frac{R \theta_\xi^S(\xi - \mathbf{X})}{(\gamma - 1)\bar{\theta}} (\theta - \bar{\theta}) \right. \right. \\
        		\qquad\qquad\qquad\qquad\qquad\qquad
                \left. \left. + \frac{\bar p^{th} v_\xi^S(\xi - \mathbf{X})}{\bar{v}} (v - \bar{v}) \right] d\xi \right\}, \\
        		\mathbf{X}(0) = x_0,
        	\end{cases}
        \end{equation}
        
        where $M$ is the specific constant chosen as 
        \begin{equation}\label{M}
        M := \frac{3\gamma(\gamma+1)p_m}{4(v_m)^2(\sigma_m)^3} ,
        \end{equation}
        and 
        \begin{equation}\label{eq:p-th}
            \bar p^{th} = \frac{R \bar{\theta}}{\bar{v}}.
        \end{equation}

Note that the initial shock location $x_0\in\bbr$ can be chosen
arbitrarily, with or without the zero-mass conditions. For simplicity,
we set $x_0=0$ in the sequel.

By the standard Cauchy-Lipschitz theory, as in \cite{KVW25}, the shift equation \eqref{X(t)} admits a unique absolutely continuous
solution $\mb X(t)$ on $[0,T]$. In fact,
\[
    |\dot{\mb X}(t)|\le C,
    \qquad
    |\mb X(t)|\le Ct,
    \qquad 0\le t\le T.
\]
        
        \subsection{A priori estimates}
        We now state the key step for the proof of Theorem \ref{thm:main}. 
        
        \begin{proposition} \label{prop:priori}
        	
        	Let $(p_-,u_-,s_-)\in\mathbb{R}_+\times\mathbb{R}\times\mathbb{R}$ be a given state. Assume that $(p,u,s)$ is a solution to system \eqref{p-u-s} on $[0,T]$ for some $T>0$, and let $(\bar p,\bar u,\bar s)$ be the composite wave constructed in \eqref{per}. Then, there exist positive constants $C_0$, $\delta_0$(recall that $\delta_R\sim|u_m-u_-|$ and $\delta_S\sim|u_+-u_m|$, so that, up to a constant, this is the smallness assumption of Theorem \ref{thm:main}), and $\varepsilon_1$ such that, if the strengths of the shock and rarefaction waves satisfy $\delta_S+\delta_R<\delta_0$, and if the perturbation satisfies
        	\begin{align*}
        		\begin{aligned}
        			&(p -\bar p, s-\bar s)\in C([0,T];H^1(\mathbb{R})), \\
        			&u-\bar u\in C([0,T];H^1(\mathbb{R})) \cap  L^2(0,T;H^2(\mathbb{R})),
        		\end{aligned}
        	\end{align*}
        	and 
        	\begin{equation}\label{apri-ass}
        	\|(p-\bar p, u-\bar u, s-\bar s)\|_{L^\infty(0,T;H^1(\mathbb{R}))} +\|(u-\bar u)_{\x\x}\|_ {L^2(0,T;L^2(\mathbb{R}))}\leq \varepsilon_1,
        	\end{equation}
        	then, for all $t\in[0,T]$, the following estimate holds:
        	\begin{align*}
        			&\|(p-\bar p, u-\bar u, s-\bar s)(t,\cdot)\|^2_{H^1(\mathbb{R})}+\delta_S\int_0^t|\dot{\mathbf{X}}(\tau)|^2 d\tau+\int_0^t\mathcal G^R(U)d\tau \\
        			&\quad+\int_0^t\mathcal G^S(U)d\tau+\int_0^t\|(p-\bar p)_\x\|_{L^2(\mathbb{R})}^2 d\tau+\int_0^t\|(u-\bar u)_\x\|^2_{H^1(\mathbb{R})} d\tau\\
        			&\le C_0 \|\big(p_0-\bar p(0,\cdot),u_0 -\bar u(0,\cdot),s_0-\bar s(0,\cdot)\big)\|^2_{H^1(\mathbb{R})}+C_0\delta_0^{\frac 14},
        	\end{align*}
            where 
            \begin{equation}\label{Gr}
            	\mathcal{G}^R(U):= \int_\bbr|v^R_\x| |(v-\bar v, \theta-\bar \theta)|^2 d\x,
            \end{equation}
            and 
            \begin{equation}\label{Gs}
            	\mathcal{G}^{S}(U):=\int_{\mathbb{R}}|(v^{S})_{\xi}^{-\mb X}||(v-\bar{v},u-\bar{u},\theta-\bar{\theta})|^{2}d\xi.
            \end{equation}
            
        	In addition, by \eqref{X(t)},
        	\begin{equation}\label{xprop}
        	|\dot{\mathbf{X}}(t)|\leq C_0(\eps_1+\delta_0),\qquad t\leq T.
        	\end{equation}
        \end{proposition}

    \subsection{‌Global existence and nonlinear stability}
    Based on the local existence theory for the solution (see \cite{Nash} for the general framework; in the present $H^1$ setting, local existence in the class of Proposition \ref{prop:priori} can be obtained by a standard iteration argument, because in Lagrangian coordinates the equations for $p$ and $s$ in \eqref{p-u-s} contain no spatial derivatives of $p$ and $s$, and are therefore ordinary differential equations in $t$ for given $u_x$, while $u$ solves a uniformly parabolic equation) and the {\it a priori} estimate in Proposition \ref{prop:priori}, we apply a standard continuation argument to obtain the global existence of the solution stated in the following proposition; see also \cite{KVW25,GLWX} for related continuation arguments.
    
    \begin{proposition}\label{thm:global}
Let $(p_-,u_-,s_-)\in\mathbb{R}_+\times\mathbb{R}\times\mathbb{R}$ be a given state, and let $(\bar p,\bar u,\bar s)$ be the composite wave constructed in \eqref{per}. Then there exist positive constants $C_0$, $\delta_0$, and $\varepsilon_0$ such that, if 
\[
|u_m-u_-|+|u_+-u_m|\leq \delta_0
\]
and the initial values $(p_0,u_0,s_0)$ satisfy \eqref{i-p}, then the Cauchy problem for \eqref{p-u-s} admits a unique global solution $(p,u,s)$ satisfying
\begin{align*}
(p-\bar p,u-\bar u,s-\bar s)\in C([0,\infty);H^1(\mathbb{R})).
\end{align*}
Moreover,
\begin{align}
\begin{aligned}\label{global-est}
&\sup_{t\geq 0}
\|(p-\bar p,u-\bar u,s-\bar s)(t,\cdot)\|_{H^1(\mathbb{R})}^2
+\delta_S\int_0^\infty |\dot{\mathbf X}(t)|^2 dt+\int_0^\infty\mathcal G^R(U)dt\\
&\quad
+\int_0^\infty\mathcal G^S(U)dt
+\int_0^\infty\|(p-\bar p)_\x\|_{L^2(\mathbb{R})}^2 dt
+\int_0^\infty\|(u-\bar u)_\x\|_{H^1(\mathbb{R})}^2 dt\\
&\leq C_0
\big\|
\big(p_0(\cdot)-\bar p(0,\cdot),
u_0(\cdot)-\bar u(0,\cdot),
s_0(\cdot)-\bar s(0,\cdot)\big)
\big\|_{H^1(\mathbb{R})}^2
+C_0\delta_0^{\frac14},
\end{aligned}
\end{align}
where $\mathcal G^R(U)$ and $\mathcal G^S(U)$ are defined in
\eqref{Gr} and \eqref{Gs}.
\end{proposition}
    
    We now use Proposition \ref{thm:global} to prove \eqref{con}-\eqref{shift-bound} for the long-time behavior as follows. We first prove the long-time behavior of the pressure and velocity. Define
\[
    g(t)\triangleq
    \|(p-\bar p)_\xi(t,\cdot)\|_{L^2(\mathbb{R})}^2
    +
    \|(u-\bar u)_\xi(t,\cdot)\|_{L^2(\mathbb{R})}^2.
\]
By Proposition \ref{thm:global}, we have
\[
    \int_0^\infty g(t)d t<\infty.
\]
Then, we need to prove
\[
    \int_0^\infty |g'(t)|d t<\infty.
\]

Using \eqref{PUS}$_1$ and integration by parts(the computation is first carried out for smooth solutions and then extended by a standard mollification argument), we find that
    \begin{align*}
        g'(t)
        ={}&
        2\int_{\mathbb{R}}
        (p-\bar p)_\xi
        \Bigg[
            \dot{\mathbf X}(t)(p^S)_{\xi\xi}^{-\mathbf X}
            -\left(
                \frac{\gamma p}{v}(u-\bar u)_\xi
            \right)_\xi\\
        &
            -\left[
                \left(
                    \frac{\gamma p}{v}
                    -\frac{\gamma\bar p}{\bar v}
                \right)\bar u_\xi
            \right]_\xi
            +\mu(\gamma-1)
            \left(
                \frac{u_\xi^2}{v^2}
                -\frac{\bar u_\xi^2}{\bar v^2}
            \right)_\xi
            -R_{1\xi}
        \Bigg]d\xi\\
        &-
        2\int_{\mathbb{R}}
        (u-\bar u)_{\xi\xi}(u-\bar u)_{t}d\xi.
    \end{align*}
    
We then estimate the terms of $g'(t)$ one by one. By Young's inequality and Lemma \ref{lem:shock-p}, we have
    \begin{align*}
        &|\dot{\mathbf X}(t)|
        \left|
            \int_{\mathbb{R}}
            (p-\bar p)_\xi
            (p^S)_{\xi\xi}^{-\mathbf X}d\xi
        \right|\le
        C\|(p-\bar p)_\xi\|_{L^2(\mathbb{R})}^2
        +C\delta_S|\dot{\mathbf X}(t)|^2.
    \end{align*}
Then, using H\"older's inequality and Sobolev's inequality, together with Proposition \ref{thm:global}, we obtain
    \begin{align*}
        &\left|
            \int_{\mathbb{R}}
            (p-\bar p)_\xi
            \left(
                \frac{\gamma p}{v}(u-\bar u)_\xi
            \right)_\xi d\xi
        \right|\\
        =&
        \gamma\left|
            \int_{\mathbb{R}}
            (p-\bar p)_\xi
            \left[
                \frac{p}{v}(u-\bar u)_{\xi\xi}
                +\left(
                    \frac{p_\xi}{v}
                    -\frac{pv_\xi}{v^2}
                \right)(u-\bar u)_\xi
            \right]d\xi
        \right|\\
        \le &
        C\|(p-\bar p)_\xi\|_{L^2(\mathbb{R})}
        \|(u-\bar u)_{\xi\xi}\|_{L^2(\mathbb{R})}+
        C\|(p-\bar p)_\xi\|_{L^2(\mathbb{R})}
        \|(p_\xi,v_\xi)\|_{L^2(\mathbb{R})}
        \|(u-\bar u)_\xi\|_{L^\infty(\mathbb{R})}\\
        \le &
        C\|(p-\bar p)_\xi\|_{L^2(\mathbb{R})}
        \left(
            \|(u-\bar u)_{\xi\xi}\|_{L^2(\mathbb{R})}
            +\|(u-\bar u)_\xi\|_{L^2(\mathbb{R})}^{1/2}
            \|(u-\bar u)_{\xi\xi}\|_{L^2(\mathbb{R})}^{1/2}
        \right)\\
        \le &
        C\|(p-\bar p)_\xi\|_{L^2(\mathbb{R})}^2
        +C\|(u-\bar u)_\xi\|_{H^1(\mathbb{R})}^2.
    \end{align*}
Similarly, we can estimate
    \begin{align*}
        &\left|
            \int_{\mathbb{R}}
            (p-\bar p)_\xi
            \left[
                \left(
                    \frac{\gamma p}{v}
                    -\frac{\gamma\bar p}{\bar v}
                \right)\bar u_\xi
            \right]_\xi d\xi
        \right|\\
     \le&
        C\int_{\mathbb{R}}
        |(p-\bar p)_\xi||\bar u_\xi|
        \left(
            |(p-\bar p)_\xi|
            +|(v-\bar v)_\xi|
        \right)d\xi\\
        &+
        C\int_{\mathbb{R}}
        |(p-\bar p)_\xi|
        \left(
            |\bar u_{\xi\xi}|
            +|\bar u_\xi|
            \bigl(|\bar p_\xi|+|\bar v_\xi|\bigr)
        \right)
        |(p-\bar p,v-\bar v)|d\xi\\
        \le&
        C\|(p-\bar p)_\xi\|_{L^2(\mathbb{R})}^2
        +C\|\bar u_\xi(v-\bar v)_\xi\|_{L^2(\mathbb{R})}^2\\
        &+
        C\int_{\mathbb{R}}
        \left[
            |\bar u_{\xi\xi}|^2
            +|\bar u_\xi|^2
            \left(
                |\bar p_\xi|^2+|\bar v_\xi|^2
            \right)
        \right]
        |(p-\bar p,v-\bar v)|^2d\xi\\
        \le&
        C\|(p-\bar p)_\xi\|_{L^2(\mathbb{R})}^2
        +C\|\bar u_\xi(v-\bar v)_\xi\|_{L^2(\mathbb{R})}^2\\
        &+
        C\int_{\mathbb{R}}
        \left(
            |v_\xi^R|
            +|(v^S)_\xi^{-\mathbf X}|
        \right)
        |(p-\bar p,v-\bar v)|^2d\xi,
    \end{align*}
    and 
    \begin{align*}
        &\left|
            \mu(\gamma-1)
            \int_{\mathbb{R}}
            (p-\bar p)_\xi
            \left(
                \frac{u_\xi^2}{v^2}
                -\frac{\bar u_\xi^2}{\bar v^2}
            \right)_\xi d\xi
        \right|\\
       \le&
        C\|(p-\bar p)_\xi\|_{L^2(\mathbb{R})}^2
        +C\|(u-\bar u)_\xi\|_{H^1(\mathbb{R})}^2
        +C\|\bar u_\xi(v-\bar v)_\xi\|_{L^2(\mathbb{R})}^2\\
        &+
        C\int_{\mathbb{R}}
        \left(
            |v_\xi^R|
            +|(v^S)_\xi^{-\mathbf X}|
        \right)|v-\bar v|^2d\xi.
    \end{align*}
Then, using Young's inequality, we have
    \begin{align*}
        &\left|
            \int_{\mathbb{R}}
            (p-\bar p)_\xi R_{1,\xi}d\xi
        \right|
        \le
        C\|(p-\bar p)_\xi\|_{L^2(\mathbb{R})}^2
        +C\|R_{1,\xi}\|_{L^2(\mathbb{R})}^2.
    \end{align*}

Combining the above estimates with Lemma \ref{lem:variables-weighted}, we obtain
\[
    \begin{aligned}
        \int_0^\infty |g'(t)|d t
        \le{}& C
        +C\int_0^\infty
        \Big(
            \|(p-\bar p)_\xi\|_{L^2(\mathbb{R})}^2
            +\|(u-\bar u)_\xi\|_{H^1(\mathbb{R})}^2\\
        &+\delta_S|\dot{\mathbf X}(t)|^2
            +\mathcal G^R(U)+\mathcal G^S(U)
        \Big)d t\\
        &+
        C\int_0^\infty
        \int_{\mathbb{R}}
        |(v^S)_\xi^{-\mathbf X}|
        |(s-\bar s)_\xi|^2d\xi d t\\
        &+
        C\int_0^\infty
        \|u_\xi^R\|_{L^\infty(\mathbb{R})}^2
        \|(s-\bar s)_\xi\|_{L^2(\mathbb{R})}^2d t\\
        &+
        C\int_0^\infty
            \|R_{1,\xi}\|_{L^2(\mathbb{R})}^2
            d t+
        C\int_0^\infty
            \|(u-\bar u)_{t}\|_{L^2(\mathbb{R})}^2
            d t.
    \end{aligned}
\]

By Proposition \ref{thm:global}, Lemmas \ref{lem:rare-v}, \ref{weight-s}, \ref{lem:wave-inter} and \eqref{psi-t}, we conclude that
\[
    \int_0^\infty
    \left(
        g(t)+|g'(t)|
    \right)d t<\infty.
\]
Thus,
 \begin{equation}\label{longtime}
     \lim_{t\to\infty}\|(p-\bar p, u-\bar u)_\x(t,\cdot)\|^2_{L^2(\bbr)}=0.
 \end{equation}
 By Sobolev's inequality, we have
  \begin{equation*}
     \|(p-\bar p, u-\bar u)\|^2_{L^{\infty}(\bbr)}\leq C \|(p-\bar p, u-\bar u)\|_{L^{2}(\bbr)}\|(p-\bar p, u-\bar u)_\x\|_{L^{2}(\bbr)},
 \end{equation*}
 which, together with \eqref{longtime}, implies
  \begin{equation*}
     \lim_{t\to\infty}\|(p-\bar p, u-\bar u)(t,\cdot)\|_{L^\infty(\bbr)}=0.
 \end{equation*}
 Then, using Lemma \ref{lem:rare-p} and the above limit, we can prove \eqref{con}.
        As for the entropy, using Sobolev's inequality and Proposition \ref{thm:global}, we can only obtain 
        $$
        \|s-\bar s\|^2_{L^{\infty}(\mathbb{R})} \leq C\| s-\bar s\|_{L^2(\mathbb{R})} \| (s-\bar s)_\x\|_{L^2(\mathbb{R})} ,   
        $$
        which means 
        $$
        \sup_{t>0}\sup_{\x\in\bbr} |(s-\bar s)(t,\x) | \le C\| s-\bar s\|_{H^1(\mathbb{R})}\leq C(\eps_0+\delta_0^{\frac 18}).
        $$
Also, using Lemma \ref{lem:rare-p} and the above inequality, we can prove \eqref{as}. 

As for the specific volume and temperature, by Lemma \ref{lem:variables-equivalence} and Proposition \ref{thm:global}, we obtain
\[
    \sup_{\x\in\bbr}
    |(v-\bar v,\theta-\bar\theta)(t,\x)|
    \le C\sup_{\x\in\bbr}
    |(p-\bar p,s-\bar s)(t,\x)|+C\delta_0.
\]
Therefore, using the smallness of $\delta_0$, we have
\[
    \sup_{t>0}
    \sup_{x\in\bbr}
    |(v-\bar v,\theta-\bar\theta)(t,x)|
    \le C\bigl(\eps_0+\delta_0^{\frac18}\bigr).
\]
Together with the $(v,\theta)$-form of Lemma \ref{lem:rare-p}, we obtain \eqref{as-vtheta} in Remark \ref{rem:vtheta-asymptotics}. 

Then, we prove \eqref{shift-bound}. Indeed, by an argument similar to that used to establish the long-time
behavior of the pressure and velocity, we obtain
\begin{equation*}
\int_0^{\infty}
\left(
\mathcal{G}^S(U)
+
\left|
\frac{d\mathcal{G}^S(U)}{dt}
\right|
\right)dt
<\infty,
\end{equation*}
which implies that
\[
\lim_{t\to\infty}\mathcal{G}^S(U)=0.
\]
By H\"older's inequality and Lemma \ref{lem:shock-v}, we have
\begin{align*}
|\dot{\mathbf X}(t)|
&\leq
\frac{C}{\delta_S}
\int_{\mathbb R}
|(v^S)_\xi^{-\mb X}|
|(v-\bar v,u-\bar u,\theta-\bar\theta)|\,d\xi\\
&\leq
\frac{C}{\delta_S}
\left(
\int_{\mathbb R}
|(v^S)_\xi^{-\mb X}|\,d\xi
\right)^{\frac12}
\left(
\int_{\mathbb R}
|(v^S)_\xi^{-\mb X}|
|(v-\bar v,u-\bar u,\theta-\bar\theta)|^2\,d\xi
\right)^{\frac12}\\
&\leq
\frac{C}{\sqrt{\deltas}}\sqrt{\mathcal G^S(U)},
\end{align*}
which proves \eqref{shift-bound}. Therefore, we complete the proof of Theorem \ref{thm:main}. It remains to prove Proposition \ref{prop:priori},
which will be the focus of the rest of the paper.
        
 \section{Relative entropy estimates}
        
  In this section, we use the relative entropy method together with the $a$-contraction method with shifts to establish the following lemma, which provides the basic $L^2$-estimates.
  The overall structure follows \cite{KVW25}. The main differences are the following: (a) no space-time control of the temperature derivative is available, so the terms involving the temperature ($\mathbf B_4$, $\mathbf B_6$ and $\mathbf S_2$ below) are closed using only the velocity dissipation, the sharp quadratic structure leading to the good terms $\mathcal G^S$ and $\mathcal G^R$, and the time decay of the rarefaction wave; (b) since the composite wave is defined in the $(p,u,s)$-coordinates, the passage between $(p,u,s)$ and $(v,u,\theta)$ produces the compatibility defects $\mathcal Q_p,\mathcal Q_s$, which are controlled in Lemmas \ref{lem:variables-boundedness} and \ref{lem:variables-equivalence}.
   \begin{lemma}\label{Re-en} 
    Under the hypotheses of Proposition \ref{prop:priori}, there  exists  $C > 0$ such that for all $t \in [0, T]$,
    \begin{equation*}
     \begin{aligned}
        & \left\| (p - \bar{p}, u - \bar{u}, s- \bar{s})(t, \cdot) \right\|_{L^2(\mathbb{R})}^2 + \delta_S \int_0^t | \dot{\mathbf{X}}(\tau) |^2 d\tau \\
        & \quad +\int_0^t\mathcal G^R(U)d\tau+\int_0^t\mathcal G^S(U)d\tau + \int_0^t \left\| (u - \bar{u})_{\xi} \right\|_{L^2(\mathbb{R})}^2 d\tau \\
        & \leq C \left\| (p_0 - \bar{p}(0, \cdot), u_0 - \bar{u}(0, \cdot), s_0 - \bar{s}(0, \cdot)) \right\|_{L^2(\mathbb{R})}^2  +C\delta_0^{\frac{1}{4}},
    \end{aligned}
   \end{equation*}
   where $\mathcal G^R(U)$ and $\mathcal G^S(U)$ are defined in
\eqref{Gr} and \eqref{Gs}.
       \end{lemma}
       The proof of Lemma \ref{Re-en} is mainly carried out in the  $(v,u,\theta)$-coordinates, and the estimates thus obtained are finally converted into corresponding estimates in the $(p,u,s)$-coordinates. For this purpose, we rewrite system \eqref{v-u-theta} as follows, based on the coordinate transformation $(t, x) \mapsto (t, \xi = x - \sigma t)$:            
        \begin{equation}\label{v-u-theta-1}
        	\begin{cases}
        		v_t - \sigma v_\xi - u_\xi = 0, \\[1mm]
        		u_t - \sigma u_\xi + p(v, \theta)_\xi = \left( \mu \frac{u_\xi}{v} \right)_\xi, \\[1mm]
        		\frac{R}{\gamma-1}\theta_t - \frac{R\sigma}{\gamma-1}\theta_\xi + p(v, \theta)u_\xi =  \mu \frac{u_\xi^2}{v}.
        	\end{cases}
        \end{equation}
        Then, the composite wave $(\bar{v}, \bar{u}, \bar{\theta})(t, \xi)$  satisfies the system   
        \begin{equation}\label{shock-v-u-theta}
        	\begin{cases}
        		\bar{v}_t - \sigma \bar{v}_\xi + \dot{\mathbf{X}}(t) (v^S)_\xi^{-\mathbf{X}} - \bar{u}_\xi = 0, \\[1mm]
        		\bar{u}_t - \sigma \bar{u}_\xi + \dot{\mathbf{X}}(t) (u^S)_\xi^{-\mathbf{X}} + \bar p^{th}_\xi = \mu \left( \frac{\bar{u}_\xi}{\bar{v}} \right)_\xi+Q_1,  \\[1mm]
        		\frac{R}{\gamma - 1} \bar{\theta}_t - \frac{R\sigma}{\gamma - 1} \bar{\theta}_\xi + \frac{R}{\gamma - 1} \dot{\mathbf{X}}(t) (\theta^S)_\xi^{-\mathbf{X}} + \bar p^{th} \bar{u}_\xi =\mu \frac{\bar{u}_\xi^2}{\bar{v}}+Q_2,
        	\end{cases}
        \end{equation}
        where $\bar p^{th}$ is defined in \eqref{eq:p-th} and the error terms are
        	\begin{align*}
        		Q_i=Q_i^I+Q_i^R, \quad i=1,2,
        	\end{align*}
        with the wave interaction terms 
        \begin{equation}\label{QI}
        \begin{aligned}
        	& Q_1^I= \left(\bar p^{th}-p^R-(p^S)^{-\mb X}\right)_\xi-\mu\left(\frac{\bar u_\x}{\bar v}-\frac{u^R_\x}{v^R}-\frac{(u^S)^{-\mb X}_\x}{(v^S)^{-\mb X}}\right)_\x,\\
        	& Q_2^I= \left(\bar p^{th}\bar u_\x-p^Ru^R_\xi-(p^S)^{-\mb X}(u^S)^{-\mb X}_\x\right)-\mu\left(\frac{\bar u^2_\x}{\bar v}-\frac{(u^R_\x)^2}{v^R}-\frac{((u^S)^{-\mb X}_\x)^2}{(v^S)^{-\mb X}}\right),
        \end{aligned}
        \end{equation}
        and the error terms due to the inviscid rarefaction wave
        \begin{equation}\label{QR}
        	\begin{aligned}
        		&Q^R_1=-\mu\left(\frac{u^R_\x}{v^R}\right)_\x,\qquad
        		Q^R_2= -\mu \frac{(u^R_\x)^2}{v^R}.
        	\end{aligned}
        \end{equation}
       
       As the method of $a$-contraction with shifts uses the Poincar\'e-type inequality, we first recall the following lemma.       		
       \begin{lemma}[Poincar\'e-type inequality] (see \cite[Lemma 2.9]{KVJEMS21}).\label{poincare}
       	For any $f : [0, 1] \to \mathbb{R}$ satisfying $\int_0^1 y(1 - y)|f'|^2 dy < \infty$,
       	\begin{equation*}
       		\int_0^1 \left| f - \int_0^1 f \, dy \right|^2 dy \leq \frac{1}{2} \int_0^1 y(1 - y)|f'|^2 dy. 
       	\end{equation*}
       \end{lemma}

        We will also use the estimates from Lemma \ref{lem:shock-v} and Lemma \ref{sharp-D} to complete the proof of Lemma \ref{Re-en}.

    \subsection{Wave interaction estimates}
     In this part, we present useful estimates for the wave interaction terms such as $Q_i^I(i=1,2)$ and $R^I_j(j=1,2,3)$ in \eqref{QI} and \eqref{RI}. Before proving Lemma \ref{lem:wave-inter}, we first establish the following lemma, which will be used in the proof of Lemmas \ref{lem:wave-inter} and \ref{lem:variables-equivalence}.

    \begin{lemma}\label{lem:variables-boundedness}
Recall that, for $\delta_0$ sufficiently small, the composite wave satisfies
\[
0<\underline c
\leq
\bar v(t,\x),\bar\theta(t,\x)
\leq \overline C,
\qquad
|\bar u(t,\x)|\leq \overline C.
\]
Then, under the same hypotheses as in Proposition \ref{prop:priori},
there exist constants $c_1,C_1>0$ such that, for all $t\in[0,T]$ and
$\x\in\bbr$,
\begin{equation}
\label{eq:physical-bounds}
0<c_1
\leq
p(t,\x),v(t,\x),\theta(t,\x)
\leq
C_1,
\end{equation}
and
\begin{equation}
\label{eq:u-s-bounds}
|s(t,\x)|\leq C_1.
\end{equation}
Furthermore,
\begin{equation}
\label{eq:vutheta-Linfty-smallness}
\sup_{0\leq t\leq T}
\|(v-\bar v,\theta-\bar\theta)(t,\cdot)\|_{L^\infty(\bbr)}
\leq
C(\varepsilon_1+\delta_0).
\end{equation}
\end{lemma}

\begin{proof}
Define the thermodynamic compatibility defects by
\begin{align}
\label{eq:Qp-expansion-proof}
\mathcal Q_p
={}&
p(v^R,\theta^R)
+
p\bigl(
(v^S)^{-\mathbf X},
(\theta^S)^{-\mathbf X}
\bigr)
-
p(v_m,\theta_m)
-
p(\bar v,\bar\theta),
\end{align}
and
\begin{align}
\label{eq:Qs-expansion-proof}
\mathcal Q_s
={}&
s(v^R,\theta^R)
+
s\bigl(
(v^S)^{-\mathbf X},
(\theta^S)^{-\mathbf X}
\bigr)
-
s(v_m,\theta_m)
-
s(\bar v,\bar\theta).
\end{align}
Applying Taylor's formula to second order to the mixed differences
\eqref{eq:Qp-expansion-proof}-\eqref{eq:Qs-expansion-proof}, we obtain
\begin{align}
\label{eq:Q-product-proof}
|\mathcal Q_p|+|\mathcal Q_s|
\leq
C
|(v^R-v_m,\theta^R-\theta_m)|
\left|
\bigl(
(v^S)^{-\mathbf X}-v_m,
(\theta^S)^{-\mathbf X}-\theta_m
\bigr)
\right|.
\end{align}
Then, from the estimates of the rarefaction waves and viscous shock waves,
we have
\begin{equation}
\label{eq:Q-Linfty-proof}
\|(\mathcal Q_p,\mathcal Q_s)(t,\cdot)\|_{L^\infty}
\leq
C\delta_0.
\end{equation}

By the {\it a priori} assumption \eqref{apri-ass} and Sobolev's
inequality, we have
\begin{equation}
\label{eq:pus-Linfty-proof}
\sup_{0\leq t\leq T}
\|(p-\bar p,u-\bar u,s-\bar s)(t,\cdot)\|_{L^\infty}
\leq
C\varepsilon_1.
\end{equation}

Since
\[
\bar p^{th}
=
\frac{R\bar\theta}{\bar v}
\]
and
\[
0<\underline c
\leq
\bar v,\bar\theta
\leq
\overline C,
\]
we have
\begin{equation}
\label{eq:bar-pth-bounds-proof}
\frac{R\underline c}{\overline C}
\leq
\bar p^{th}
\leq
\frac{R\overline C}{\underline c}.
\end{equation}
Since
\[
p
=
\bar p^{th}
+
\mathcal Q_p
+
(p-\bar p),
\]
it follows from \eqref{eq:Q-Linfty-proof}-\eqref{eq:bar-pth-bounds-proof} that
\begin{equation}
\label{eq:p-bounds}
0<c\leq p(t,\x)\leq C,
\qquad
(t,\x)\in[0,T]\times\bbr,
\end{equation}
where we have used the smallness of $\varepsilon_1$ and $\delta_0$.

Similarly,
\begin{equation}
\label{eq:s-bounds}
|s(t,\x)|\leq C.
\end{equation}

Using \eqref{pes}, together with \eqref{eq:p-bounds} and \eqref{eq:s-bounds}, we obtain
\[
0<c_1
\leq
v(t,\x),\theta(t,\x)
\leq
C_1.
\]
This proves \eqref{eq:physical-bounds} and
\eqref{eq:u-s-bounds}.

By the definitions of $\mathcal Q_p$ and $\mathcal Q_s$ in
\eqref{eq:Qp-expansion-proof} and \eqref{eq:Qs-expansion-proof}, we have
\begin{align*}
(p-\bar p+\mathcal Q_p,
s-\bar s+\mathcal Q_s)
=
\bigl(
p(v,\theta)-p(\bar v,\bar\theta),
s(v,\theta)-s(\bar v,\bar\theta)
\bigr).
\end{align*}
Using \eqref{eq:physical-bounds}, \eqref{eq:u-s-bounds}, and the mean
value theorem, we obtain
\begin{equation}\label{v-p}
\begin{aligned}
|(v-\bar v,\theta-\bar\theta)|
\leq
C\Bigl(
|(p-\bar p,s-\bar s)|
+
|\mathcal Q_p|+|\mathcal Q_s|
\Bigr).
\end{aligned}
\end{equation}
Then, it follows from \eqref{eq:Q-Linfty-proof},
\eqref{eq:pus-Linfty-proof}, and \eqref{v-p} that
\[
\begin{aligned}
\sup_{0\leq t\leq T}
\|(v-\bar v,\theta-\bar\theta)(t,\cdot)\|_{L^\infty}
\leq
C(\varepsilon_1+\delta_0),
\end{aligned}
\]
which proves \eqref{eq:vutheta-Linfty-smallness}.
\end{proof}
    
        Therefore, the ODE \eqref{X(t)}, together with Lemmas \ref{lem:shock-v} and \ref{lem:variables-boundedness}, gives
        \begin{equation}\label{dxbound}
        \begin{aligned}
        	|\dot{\mb X}(t)| &\le \frac{C}{\deltas}  \|(v-\bar v, u-\bar u, \theta-\bar \theta)(t,\cdot) \|_{L^\infty(\bbr)} \int_\bbr ( v^S)^{-\mb X}_\xi d\xi \\
        	&\le C \|(v-\bar v, u-\bar u, \theta-\bar \theta)(t,\cdot) \|_{L^\infty(\bbr)} \le C(\eps_1+\delta_0).
        \end{aligned}
        \end{equation}
        This proves \eqref{xprop}, and 
        will be used to get the wave interaction estimates in Lemma \ref{lem:wave-inter}. 
        \begin{lemma}\label{lem:wave-inter}
        	Let $\mb X$ be the shift defined by \eqref{X(t)}, and let $Q_i^I(i=1,2)$ and $R_j^I(j=1,2,3)$ be defined by \eqref{QI} and \eqref{RI}, respectively. Under the same hypotheses of Proposition \ref{prop:priori}, there exists a constant $C>0$ such that, for $i=1,2$, $j=1,2,3$ and all $t\le T$,
        	\begin{equation}\label{wave-interactions}
        	\begin{aligned}
        		&\|(Q_i^I,R_j^I,R^I_{1,\x},R^I_{3,\x})\|_{L^2(\bbr)} \le C \deltas \deltar \left(e^{-C \deltas t}+|(1+\deltar t)|^{-\frac{5}{2}}\right),\\[3mm]
        		& \| |(v^S)^{-\mb X}_\xi| |\big(v^R -v_m,\theta^R -\theta_m\big) |\|_{L^2(\bbr)}\leq C\deltas^{3/2} \deltar \left( e^{-C \deltas t}+|(1+\deltar t)|^{-\frac{5}{2}}\right).
        	\end{aligned}
        	\end{equation}
        \end{lemma}
        \begin{proof} For brevity, we only estimate $\|R_2^I\|_{L^2(\bbr)}$, since the estimates of $\|Q_i^I\|_{L^2(\bbr)}(i=1,2)$, $\|R_j^I\|_{L^2(\bbr)}(j=1,3)$ and $\|R_{j,\x}^I\|_{L^2(\bbr)}(j=1,3)$ are obtained similarly. 
        	Recall
        	\[
        	R_2^I:= -\mu\left(\frac{\bar u_\x}{\bar v}-\frac{u^R_\x}{v^R}-\frac{(u^S)^{-\mb X}_\x}{(v^S)^{-\mb X}}\right)_\x.
        	\]
        	Since $\bar{u}_\xi = u^R_\xi+(u^S)^{-\mb X}_\xi$ and $\bar{v}_\xi = v^R_\xi+(v^S)^{-\mb X}_\xi$, $R_2^I$ can be estimated as follows
        		\begin{align*} 
        			&\left|\mu\left(\frac{\bar u_\x}{\bar v}-\frac{u^R_\x}{v^R}-\frac{(u^S)^{-\mb X}_\x}{(v^S)^{-\mb X}}\right)_\x\right|\\[3mm]
        			&\le C\Big[|(u^R_{\xi\xi}, u^R_\x v^R_\xi)||(v^S)^{-\mb X}-v_m| +|((u^S)^{-\mb X}_{\xi\xi}, (u^S)^{-\mb X}_{\xi} (v^S)^{-\mb X}_{\xi})||v^R-v_m| \\
        			&\quad  +|u^R_\xi||(v^S)^{-\mb X}_\xi|+|(u^S)^{-\mb X}_\xi|| v^R_\xi| \Big]\\
        			&\le C\Big[|u^R_{\xi}||(v^S)^{-\mb X}-v_m| +|(u^S)^{-\mb X}_{\xi}||v^R-v_m|  +|u^R_\xi||(v^S)^{-\mb X}_\xi|+|(u^S)^{-\mb X}_\xi|| v^R_\xi| \Big]\\
        			&\le C(A_1+A_2+A_3),
        	\end{align*}
        	where
          \begin{align*}
           A_1&:=|u^R_{\xi}||(v^S)^{-\mb X}-v_m|,\\
           A_2&:=|(u^S)^{-\mb X}_{\xi}||v^R-v_m|,\\
           A_3&:=|u^R_{\xi}||(v^S)^{-\mb X}_{\xi}|
                +|(u^S)^{-\mb X}_{\xi}||v^R_{\xi}|.
        \end{align*}
        	First, for fixed $t\in[0,T]$ and $\lam_{1m}:=\lam_1(v_m, \theta_m)<0$, we have
        	\begin{align*}
        		\|A_1 \|_{L^2(\bbr)} &\le  \bigg[\bigg(\int_{\big\{\xi+\s t \leq \frac{\lam_{1m}}{2}(1+t)\big\}}+\int_{\big\{\xi+\s t > \frac{\lam_{1m}}{2}(1+t)\big\}}\bigg) \big|u^R_\xi\big|^2\big|(v^S)^{-\mb X}-v_m \big|^2 d\xi\bigg]^{\frac12}\\
        		& =: \sqrt{A_{11}+A_{12}}.
        	\end{align*}
        	By \eqref{dxbound}, 
        	it holds that
        	$$
        	|\dot{\mb{X}}(t)|\leq C(\eps_1+\delta_0),\qquad 0\le t\le T,
        	$$
        	which together with $\mb X(0)=0$ yields
        	$$
        	|\mb{X}(t)|\leq C(\eps_1+\delta_0) t,\qquad 0\le t\le T.
        	$$
        	Let us take $\eps_1$ and $\delta_0$ sufficiently small that the above bound is less than $ \frac{\s t}{4}$, that is,
        	\begin{equation*}
        	|\mb{X}(t)|\leq C(\eps_1+\delta_0) t \leq \frac{\s t}{4}.
        	\end{equation*}
        	Since 
        	\begin{equation*}
        		\begin{aligned}
        	\xi+\s t \leq  \frac{\lam_{1m}}{2}(1+t)\quad\Rightarrow\quad\xi- \mb{X}(t)&\leq \frac{\lam_{1m}}{2}(1+t)-\s t-\mb{X}(t)\\
        	&\leq  \frac{\lam_{1m}}{2}(1+t)-\frac 34\s t \leq -\frac 34\s t \leq0,
        	     \end{aligned}
        	\end{equation*}
        	and so $|\xi- \mb{X}(t)|\ge \frac 34\s t $, it follows from Lemma \ref{lem:shock-v} that for all $\xi$ with $\xi+\s t \leq  \frac{\lam_{1m}}{2}(1+t)< 0$,
        	\begin{align*}
        		\big|(v^S)^{-\mb X}-v_m \big| \le C\deltas e^{-C\deltas |\xi -\mb X(t)|} \le C\deltas e^{-\frac{C\deltas |\xi -\mb X(t)|}{2}} e^{-C\deltas t}.
        	\end{align*}
        	Thus,
        	\begin{align*}
        		A_{11} \le C(\deltas)^2 e^{-2C\deltas t} \int_\bbr \big|u^R_\xi\big|^2 d\xi,
        	\end{align*}
        	which together with Lemma \ref{lem:rare-v} yields
        	\[
        	A_{11} \le C \deltas^2 \deltar^3 e^{-2C\deltas t}.
        	\]
        	On the other hand, since 
        	\begin{align*}
        	\xi+\s t > \frac{\lam_{1m}}{2}(1+t)  \quad\Rightarrow\quad \xi+\s t-\lam_{1m}(1+t)&> -\frac{\lam_{1m}}{2}(1+t)\\
        	&=\frac{|\lam_{1m}|}{2}(1+t)>0,
        	\end{align*}
        	it follows from Lemma \ref{lem:rare-v} that  for all $\xi$ with $\xi+\s t > \frac{\lam_{1m}}{2}(1+t) $,
        	\begin{equation*}
         |(v_\xi^R,u^R_\x,\theta_\x^R)| \le C\deltar^2{\big(1+ \deltar|\xi +\s t-\lambda_{1m}(1+t)|\big)^{-4}}.
        	\end{equation*}
        	This together with  Lemma \ref{lem:shock-v} and Lemma \ref{lem:rare-v} implies
        	\begin{align*}
        		A_{12} &\le C \deltas^2\deltar^4  \int_{\big\{\xi+\s t > \frac{\lam_{1m}}{2}(1+t)\big\}}  (1+\deltar |\xi +\s t-\lambda_{1m}(1+t) |)^{-8}d\xi \\
        		&\le  C \deltas^2 \deltar^3 |(1+\deltar t)|^{-7}.
        	\end{align*}
        	Thus,
        	\[
        	\|A_1 \|_{L^2(\bbr)} \le C \deltas\deltar^{3/2} e^{-C \deltas t}+C \deltas\deltar^{3/2} |(1+\deltar t)|^{-\frac{7}{2}}.
        	\]
        	To estimate $ \|A_2\|_{L^2(\bbr)}$ and $ \|A_3\|_{L^2(\bbr)}$, we consider
        	\begin{equation*}
        	\begin{aligned}
        		\|A_2 \|_{L^2(\bbr)}
        		\le C \bigg[\bigg(\int_{\big\{\xi+\s t \leq \frac{\lam_{1m}}{2}(1+t)\big\}}+\int_{\big\{\xi+\s t > \frac{\lam_{1m}}{2}(1+t)\big\}}\bigg)    |(u^S)^{-\mb X}_{\xi}|^2 |v^R-v_m|^2 d\x \bigg]^{\frac{1}{2}}
        	\end{aligned}
        	\end{equation*}
        	and 
        	\begin{equation*}
        	\begin{aligned}
        		\|A_3 \|_{L^2(\bbr)}
        		&\le C \bigg[\bigg(\int_{\big\{\xi+\s t \leq \frac{\lam_{1m}}{2}(1+t)\big\}}+\int_{\big\{\xi+\s t > \frac{\lam_{1m}}{2}(1+t)\big\}}\bigg)    |(v^S)^{-\mb X}_{\xi}|^2 |u_\x^R|^2   \\
        		&\qquad +|(u^S)^{-\mb X}_\xi|^2| v^R_\xi|^2 d\x \bigg]^{\frac{1}{2}}.
        	\end{aligned}
        	\end{equation*}
        	Arguing as for $A_1$, we have
        	\[
        	\|A_2 \|_{L^2(\bbr)}+ \|A_3 \|_{L^2(\bbr)}\le C \deltas\deltar \left( e^{-C \deltas t}+|(1+\deltar t)|^{-\frac{5}{2}}\right).
        	\]
        	Combining the above estimates, we have
        	\begin{equation*}
        	\big\|R_2^I\big\|_{L^2(\bbr)}\leq C  \deltas \deltar\left(e^{-C \deltas t}+|(1+\deltar t)|^{-\frac{5}{2}}\right).
        	\end{equation*}
        	Therefore, the wave interaction estimate $\eqref{wave-interactions}_1$ holds true. Moreover, using Lemmas~\ref{lem:rare-v} and \ref{lem:shock-v} and arguing as in the estimate of $A_1$, we obtain
        	$$
        	 \| |(v^S)^{-\mb X}_\xi| |\big(v^R -v_m,\theta^R -\theta_m\big) |\|_{L^2(\bbr)}\leq C\deltas^{3/2}  \left(\deltar e^{-C \deltas t}+\deltar|(1+\deltar t)|^{-\frac{5}{2}}\right),
        	 $$
        	which proves the wave interaction estimate $\eqref{wave-interactions}_2$.
        \end{proof}  
        
\subsection{‌Equivalence of variables}

Throughout the subsequent proof, we shall alternate between the $(p,u,s)$ and $(v,u,\theta)$ formulations. For completeness, the following lemma provides their equivalence.

\begin{lemma}\label{lem:variables-equivalence}
Under the same hypotheses as in Proposition \ref{prop:priori}, there exists a constant $C>0$ such that, for all $t\in[0,T]$ and $\x\in\bbr$,
\begin{equation}
\label{eq:corrected-pointwise-equivalence}
\begin{aligned}
    C^{-1}|(v-\bar v,\theta-\bar\theta)|
    \leq
    |(p-\bar p+\mathcal Q_p,
       s-\bar s+\mathcal Q_s)|
    \leq
    C|(v-\bar v,\theta-\bar\theta)|,
\end{aligned}
\end{equation}
where $\mathcal Q_p$ and $\mathcal Q_s$ are defined in \eqref{eq:Qp-expansion-proof} and \eqref{eq:Qs-expansion-proof}.
Consequently,
\begin{equation}
\label{eq:raw-pointwise-forward}
    |(v-\bar v,\theta-\bar\theta)|
    \leq
    C\Bigl(
        |(p-\bar p,s-\bar s)|
        +
        |\mathcal Q_p|+|\mathcal Q_s|
    \Bigr),
\end{equation}
and
\begin{equation}
\label{eq:raw-pointwise-backward}
    |(p-\bar p,s-\bar s)|
    \leq
    C\Bigl(
        |(v-\bar v,\theta-\bar\theta)|
        +
        |\mathcal Q_p|+|\mathcal Q_s|
    \Bigr).
\end{equation}
Furthermore,
\begin{equation}
\label{eq:H1-equivalence-forward}
\begin{aligned}
    \|(v-\bar v,\theta-\bar\theta)\|_{H^1(\bbr)}
    \leq
    C\|(p-\bar p,s-\bar s)\|_{H^1(\bbr)}
    +
    C\delta_0,
\end{aligned}
\end{equation}
and
\begin{equation}
\label{eq:H1-equivalence-backward}
\begin{aligned}
    \|(p-\bar p,s-\bar s)\|_{L^2(\bbr)}
    \leq
    C\|(v-\bar v,\theta-\bar\theta)\|_{L^2(\bbr)}
    +
    C\delta_0.
\end{aligned}
\end{equation}
In particular,
\begin{equation}
\label{eq:vutheta-smallness}
\begin{aligned}
    \sup_{0\leq t\leq T}
    \|(v-\bar v,\theta-\bar\theta)\|_{H^1(\bbr)}
    \leq
    C(\varepsilon_1+\delta_0).
\end{aligned}
\end{equation}
\end{lemma}

\begin{proof}
We first establish the estimates on the thermodynamic compatibility
defects. 
By \eqref{eq:Q-product-proof} and an argument similar to the proof of Lemma \ref{lem:wave-inter}, we obtain
\begin{equation}
\label{eq:Q-L2-proof}
\begin{aligned}
    \|(\mathcal Q_p,\mathcal Q_s)(t,\cdot)\|_{L^2}
    \leq
    C\delta_R\delta_S^{\frac12}
    e^{-c\delta_St}+
    C\delta_S\delta_R^{\frac12}
    (1+\delta_Rt)^{-\frac52}.
\end{aligned}
\end{equation}

Differentiating \eqref{eq:Qp-expansion-proof} and
\eqref{eq:Qs-expansion-proof}, and applying the Taylor formula once
again, we obtain
\begin{align*}
    |(\mathcal Q_p,\mathcal Q_s)_\x|
    \leq
    &C
    |(v^R, \theta^R)_\x|
    \left|
        \bigl(
            (v^S)^{-\mathbf X}-v_m,
            (\theta^S)^{-\mathbf X}-\theta_m
        \bigr)
    \right|\\
    &+C
    |(v^R-v_m,\theta^R-\theta_m)|
    \left|
        \bigl(
            (v^S)^{-\mathbf X},
            (\theta^S)^{-\mathbf X}
        \bigr)_\x
    \right|.
\end{align*}
Similarly, we have
\begin{equation}
\label{eq:Qx-sup-proof}
    \|(\mathcal Q_p,\mathcal Q_s)_\x(t,\cdot)\|_{L^2}
    \leq
    C\delta_R\delta_S^{\frac12}
    e^{-c\delta_St}+
    C\delta_S\delta_R^{\frac12}
    (1+\delta_Rt)^{-\frac52}.
\end{equation}
Together with \eqref{eq:Q-L2-proof}, this yields
\begin{equation}
\label{eq:Q-H1-sup-proof}
    \sup_{0\leq t\leq T}
    \|(\mathcal Q_p,\mathcal Q_s)(t,\cdot)\|_{H^1}
    \leq
    C\delta_0.
\end{equation}

We next prove the pointwise equivalence between the two sets of
thermodynamic perturbations. By the definitions of
$\mathcal Q_p$ and $\mathcal Q_s$ in \eqref{eq:Qp-expansion-proof} and \eqref{eq:Qs-expansion-proof}, we have
\begin{align}
\label{eq:corrected-identity-proof}
    (p-\bar p+\mathcal Q_p,
      s-\bar s+\mathcal Q_s)=
    \bigl(
        p(v,\theta)-p(\bar v,\bar\theta),
        s(v,\theta)-s(\bar v,\bar\theta)
    \bigr).
\end{align}
Applying the mean value theorem to \eqref{eq:corrected-identity-proof}, together with \eqref{eq:p-bounds}-\eqref{eq:s-bounds}, we obtain
\[
\begin{aligned}
    C^{-1}|(v-\bar v,\theta-\bar\theta)|
    \leq
    |(p-\bar p+\mathcal Q_p,
       s-\bar s+\mathcal Q_s)|
       \leq
    C|(v-\bar v,\theta-\bar\theta)|,
\end{aligned}
\]
which proves \eqref{eq:corrected-pointwise-equivalence}.
The triangle inequality then yields
\eqref{eq:raw-pointwise-forward} and
\eqref{eq:raw-pointwise-backward}.

To prove \eqref{eq:H1-equivalence-forward}, we differentiate
\eqref{eq:corrected-identity-proof}. Using the smoothness of the
thermodynamic transformation, we obtain
\begin{align*}
    &|(p-\bar p+\mathcal Q_p,
       s-\bar s+\mathcal Q_s)_\x|\\
    \leq &
    C|(v-\bar v,\theta-\bar\theta)_\x|
    +
    C|(\bar v,\bar\theta)_\x|
      |(v-\bar v,\theta-\bar\theta)|.
\end{align*}
Applying the same argument to the inverse thermodynamic
transformation gives
\begin{align}
\label{eq:corrected-derivative-backward-proof}
    |(v-\bar v,\theta-\bar\theta)_\x|
    &\leq
    C|(p-\bar p+\mathcal Q_p,
       s-\bar s+\mathcal Q_s)_\x|
    \nonumber\\
    &\quad+
    C|(\bar v,\bar\theta)_\x|
      |(p-\bar p+\mathcal Q_p,
         s-\bar s+\mathcal Q_s)|.
\end{align}
Using \eqref{eq:corrected-derivative-backward-proof}, we have
\begin{align*}
    \|(v-\bar v,\theta-\bar\theta)_\x\|_{L^2}
    \leq{}&
    C\|(p-\bar p,s-\bar s)_\x\|_{L^2}+
    C\|(\mathcal Q_p,\mathcal Q_s)_\x\|_{L^2}
    \\
    &+
    C\|(p-\bar p,s-\bar s)\|_{L^2}+
    C\|(\mathcal Q_p,\mathcal Q_s)\|_{L^2}.
\end{align*}
Together with \eqref{eq:raw-pointwise-forward} and \eqref{eq:Q-H1-sup-proof}, this gives
\[
\begin{aligned}
    \|(v-\bar v,\theta-\bar\theta)\|_{H^1(\bbr)}
    &\leq
    C\|(p-\bar p,s-\bar s)\|_{H^1(\bbr)}+
    C\|(\mathcal Q_p,\mathcal Q_s)\|_{H^1(\bbr)}\\
    &\leq
    C\|(p-\bar p,s-\bar s)\|_{H^1(\bbr)}+
    C\delta_0,
\end{aligned}
\]
which proves \eqref{eq:H1-equivalence-forward}.
Similarly, using \eqref{eq:raw-pointwise-backward} and \eqref{eq:Q-H1-sup-proof}, we obtain
\[
\begin{aligned}
    \|(p-\bar p,s-\bar s)\|_{L^2(\bbr)}
    &\leq
    C\|(v-\bar v,\theta-\bar\theta)\|_{L^2(\bbr)}+
    C\|(\mathcal Q_p,\mathcal Q_s)\|_{L^2(\bbr)}\\
    &\leq
    C\|(v-\bar v,\theta-\bar\theta)\|_{L^2(\bbr)}+
    C\delta_0,
\end{aligned}
\]
which proves \eqref{eq:H1-equivalence-backward}.

Finally, by the {\it a priori} assumption \eqref{apri-ass} and \eqref{eq:Q-H1-sup-proof}, we conclude that
\begin{align*}
    \sup_{0\leq t\leq T}
    \|(v-\bar v,\theta-\bar\theta)(t,\cdot)\|_{H^1}
    \leq
    C(\varepsilon_1+\delta_0).
\end{align*}
This proves \eqref{eq:vutheta-smallness}.
\end{proof}
        
        \subsection{Construction of the weight function}
        We define the weight function $a$ by
        \begin{equation}\label{weight}
        	a(\xi):=1+\frac{\lam}{\deltas}(v^S(\xi)-v_m),
        \end{equation}
        where the constant $\lam$ is chosen to be so small but much larger than $\deltas=v_+-v_m>0$ such that 
        \begin{equation}\label{lamsmall}
        \deltas\ll \lam \le C\sqrt{\deltas}.
        \end{equation}
        Thus, we have
        \begin{equation*}
        	1<a(\xi)<1+\lam,
        \end{equation*}
        and 
        \begin{equation}\label{a-prime}
        	a^\prime(\xi)=\frac{\lam}{\deltas}v^S_\xi>0,
        \end{equation}
        and so,
        \begin{equation}\label{d-weight}
        |a'|\sim\frac{\lam}{\deltas} |u^S_\xi|\sim \frac{\lam}{\deltas} |\theta^S_\xi|.
        \end{equation}
        
        \subsection{Relative entropy method}
        
        To simplify the relative entropy calculation, we will use the non-conserved quantities \(U=(v, u, \theta)^{t}\) and \(\bar{U}=(\bar{v}, \bar{u}, \bar{\theta})^{t}\), where \(U=(v, u, \theta)^{t}\) is the solution to the reformulated system $\eqref{v-u-theta-1}$, and \(\bar{U}=(\bar{v}, \bar{u}, \bar{\theta})^{t}\) is defined in \eqref{per-v-theta}.
        
        First of all, using the Gibbs relation and an argument similar to those in \cite{KVW25,GLWX} (see also \cite{Dafermos1} for the relative entropy method), we obtain the relative entropy of $U$ with respect to $\bar{U}$:
       $$
        \eta (U|\bar{U}) = R\left( \frac{v}{\bar{v}} - 1 - \log \frac{v}{\bar{v}} \right) + \frac{R}{\gamma - 1}\left( \frac{\theta}{\bar{\theta}} - 1 - \log \frac{\theta}{\bar{\theta}} \right) + \frac{(u - \bar{u})^2}{2\bar{\theta}}.
    $$
        Using the convex function \(\Phi(z):=z - 1 - \log z\), we have 
        \[
        \eta(U | \bar{U}) = R \Phi\left(\frac{v}{\bar{v}}\right) + \frac{R}{\gamma - 1} \Phi\left(\frac{\theta}{\bar{\theta}}\right) + \frac{(u - \bar{u})^2}{2 \bar{\theta}}.
        \]
        
        So, the relative entropy weighted by \(\bar{\theta}\) is given by 
        \begin{equation}\label{rela-entro}
        	\bar{\theta} \eta(U | \bar{U}) = R \bar{\theta} \Phi\left(\frac{v}{\bar{v}}\right) + \frac{R \bar{\theta}}{\gamma - 1} \Phi\left(\frac{\theta}{\bar{\theta}}\right) + \frac{(u - \bar{u})^2}{2}.
        \end{equation}
        Below, we will compute the evolution of the relative entropy of \(U\) and \(\bar{U}\) weighted by \(a(\xi - \mb X) \bar{\theta}(t, \xi)\): 
        \[
        \int_{\mathbb{R}} a^{-\mb X}(\xi) \bar{\theta}(t, \xi) \eta(U(t, \xi) | \bar{U}(t, \xi)) d\xi.
        \]
        Using the systems \eqref{v-u-theta-1} and \eqref{shock-v-u-theta} and following a computation similar to that in \cite{KVW25}, we obtain the following lemma without proof.
        \begin{lemma}
        	\label{lem3.2}
        	Let \(a\) be the weight function defined by \eqref{weight}. Then, 
        	\begin{equation*}
        	\frac{d}{dt} \int_{\mathbb{R}} a^{-\mb X}(\xi) \bar{\theta}(t, \xi) \eta(U(t, \xi) | \bar{U}(t, \xi)) d\xi = \dot{\mb X}(t) \mb Y(U) + \mathcal{J}^{bad}(U) - \mathcal{J}^{good}(U),
        	\end{equation*}
        	where 
        	\begin{align*}
        		\mb Y(U) &:= -\int_{\mathbb{R}} a_{\xi}^{-\mb X} \bar{\theta} \eta(U | \bar{U}) d\xi \\
        		& \quad + \int_{\mathbb{R}} a^{-\mb X}\left[ -R(\theta^S)_{\xi}^{-\mb X} \Phi\left(\frac{v}{\bar{v}}\right) - \frac{R}{\gamma - 1}(\theta^S)_{\xi}^{-\mb X} \Phi\left(\frac{\theta}{\bar{\theta}}\right) \right] d\xi \\
        		& \quad + \int_{\mathbb{R}} a^{-\mb X}\left[ (u^S)_{\xi}^{-\mb X}(u - \bar{u}) + \frac{(v^S)_{\xi}^{-\mb X} \bar p^{th}}{\bar{v}}(v - \bar{v}) + \frac{R}{\gamma - 1} \frac{(\theta^S)_{\xi}^{-\mb X}}{\bar{\theta}}(\theta - \bar{\theta}) \right] d\xi,
        	\end{align*}
        	\begin{align*}
        		\mathcal{J}^{bad}(U) &:= \int_{\mathbb{R}} a_\xi^{-\mathbf{X}} (u - \bar{u})(p - \bar p^{th}) d\xi \\
        		&+ \int_{\mathbb{R}} a^{-\mathbf{X}} \left[ R\left( (\theta^R_t-\sigma \theta_\x^R) - \sigma (\theta^S)_\xi^{-\mathbf{X}} \right) \Phi\left( \frac{v}{\bar{v}} \right) - \frac{\bar p^{th} \bar{u}_\xi}{v \bar{v}} (v - \bar{v})^2 \right] d\xi \\
        		&+ \int_{\mathbb{R}} a^{-\mathbf{X}} \left[ \frac{R}{\gamma - 1}\left( (\theta^R_t-\sigma \theta_\x^R) - \sigma (\theta^S)_\xi^{-\mathbf{X}} \right) \Phi\left( \frac{\theta}{\bar{\theta}} \right) - \frac{\bar{u}_\xi}{\theta} (\theta - \bar{\theta})(p - \bar p^{th})\right. \\
        		&\quad\quad\left.+ \bar p^{th} \bar{u}_\xi \frac{(\theta - \bar{\theta})^2}{\theta \bar{\theta}} \right] d\xi + \int_{\mathbb{R}} a^{-\mathbf{X}} \left[ \mu (u - \bar{u}) \left( \frac{u_\xi}{v} - \frac{\bar{u}_\xi}{\bar{v}} \right) \right]_\xi d\xi \\
        		&+ \int_{\mathbb{R}} a^{-\mathbf{X}} \left[ -\mu \bar{u}_\xi \left( \frac{1}{v} - \frac{1}{\bar{v}} \right) (u - \bar{u})_\xi + \mu \frac{\theta - \bar{\theta}}{\theta} \left( \frac{u_\xi^2}{v} - \frac{\bar{u}_\xi^2}{\bar{v}} \right) -  \frac{(\theta - \bar{\theta})^2}{\theta \bar{\theta}}  \right.\\
                &\qquad\qquad\cdot\left.\left(\mu\frac{\bar{u}_\xi^2}{\bar{v}} +Q_2\right)\right]d\x- \int_{\mathbb{R}} a^{-\mathbf{X}} \left[Q_1(u-\bar u)+Q_2(1-\frac{\bar\theta}{\theta})\right] d\xi,
           \end{align*}
        	and 
        	\begin{align*}
        	\mathcal{J}^{good}(U) := \sigma \int_{\mathbb{R}} a_{\xi}^{-\mb X} \bar{\theta} \eta(U | \bar{U}) d\xi + \int_{\mathbb{R}} a^{-\mb X}\left( \frac{\mu}{v} |(u - \bar{u})_{\xi}|^2 \right) d\xi,
        	\end{align*}
        where $\bar p^{th}$ is defined in \eqref{eq:p-th}.
        \end{lemma}

        \begin{remark}
        	\label{rem:4.1}
        	Since \(\sigma a'(\xi) > 0\), \(\mathcal{J}^{good}(U)\) consists of nonnegative terms.
        \end{remark}

       	\subsection{Decompositions}
       	
       	First of all, we will decompose the second and third terms of $\mathcal{J}^{bad}$ into the main term $\mathbf{B}_1$ of leading order, and some lower order terms
       	
       	\begin{equation}\label{B1}
       		\begin{split}
       			\mathbf{B}_1(U) &:= \int_{\mathbb{R}} a^{-\mb X}(v^S)^ {-\mathbf{X}}_\xi \left[ \frac{(\gamma+1) p_m \sigma_m}{2 (v_m)^2} (v - \bar{v})^2 + \frac{R \sigma_m}{2 v_m \theta_m} (\theta - \bar{\theta})^2 \right. \\
       			&\qquad\quad \left. - \frac{p_m \sigma_m}{v_m \theta_m} (v - \bar{v})(\theta - \bar{\theta}) \right] d\xi,
       		\end{split}
       	\end{equation}

       	where the constants $p_m$ and $\sigma_m$ are 	
       	\begin{equation}\label{p-sigma-}
       		p_m := p(v_m, \theta_m) = \frac{R \theta_m}{v_m} \quad \text{and} \quad \sigma_m := \sqrt{\frac{\gamma p_m}{v_m}} = \frac{\sqrt{\gamma R \theta_m}}{v_m}.
       	\end{equation}
            	
        Therefore, it follows from Lemma \ref{lem3.2} that 
        \begin{equation}\label{d-a-1}
            \begin{aligned}
            \frac{d}{d t} & \int_{\mathbb{R}} a^{-\mb X}(\xi) \bar{\theta}(t, \xi) \eta(U(t, \xi) | \bar{U}(t, \xi)) d \xi \\
            & = \dot{\mathbf X}(t) \mathbf Y(U)+\sum_{i=1}^6 \mathbf B_{i}(U)-\mathbf G(U)-\mathbf D(U)+\mb S_1(U)+\mb S_2(U),
            \end{aligned}
        \end{equation}
    where $\mathbf {B}_{1}$ is as in $\eqref{B1}$, and 
    \begin{align*}
        &\mathbf B_{2}(U):=\int_{\mathbb{R}} a_{\xi}^{-\mb X}(u-\bar{u})(p-\bar p^{th}) d \xi,\\
        &\mathbf B_{3}(U):=-\int_{\mathbb{R}} a_{\xi}^{-\mb X}\mu(u-\bar{u})\left(\frac{u_{\xi}}{v}-\frac{\bar{u}_{\xi}}{\bar{v}}\right) d \xi,\\
        &\mathbf B_{4}(U):= \int_{\mathbb{R}} a^{-\mb X}\left\{-\mu \bar{u}_{\xi}\left(\frac{1}{v}-\frac{1}{\bar{v}}\right)(u-\bar{u})_{\xi}+\mu \frac{\theta-\bar{\theta}}{\theta}\left(\frac{u_{\xi}^{2}}{v}-\frac{\bar{u}_{\xi}^{2}}{\bar{v}}\right)-\frac{(\theta-\bar{\theta})^{2}}{\theta \bar{\theta}}\mu \frac{\bar{u}_{\xi}^{2}}{\bar{v}}\right\} d \xi,\\
        &\mb B_5(U):=\int_\bbr \!a^{-\mb X} \Big[R\Big((\theta_t^R-\sigma \theta_\x^R)-\sigma (\theta^S)^{-\mb X}_\xi \Big)\Phi\left(\frac{v}{\bar v}\right)-\frac{\bar p^{th}\bar u_\xi}{v\bar v}(v-\bar v)^2\\
        &\qquad\qquad\qquad-(v^S)_\xi^{-\mb X} \frac{(\gamma+1)p_m\sigma_m}{2(v_m)^2} (v-\bar v)^2\Big]d\x,\\
        &\mb B_6(U) := \int_\bbr \!a^{-\mb X} \left\{\frac{R}{\gamma-1}\Big((\theta_t^R-\sigma \theta_\x^R)-\sigma (\theta^S)^{-\mb X}_\xi \Big)\Phi\left(\frac{\theta}{\bar \theta}\right) - \frac{\bar u_\xi}{\theta}(\theta-\bar \theta)(p-\bar p^{th}) \right.\\
        &\left.\qquad\qquad\qquad+\bar p^{th} \bar u_\xi \frac{(\theta-\bar\theta)^2}{\theta\bar\theta} - (v^S)_\xi^{-\mb X} \Big[\frac{R\sigma_m}{2v_m\theta_m}(\theta-\bar\theta)^2-\frac{p_m\sigma_m}{v_m\theta_m}(v-\bar v)(\theta-\bar\theta)\Big]\right\}d\x,
    \end{align*}
    and 
    \begin{align*}
        \mathbf G(U) & :=\sigma \int_{\mathbb{R}} a_{\xi}^{-\mb X} \bar{\theta} \eta(U | \bar{U}) d \xi, \\
        \mathbf D(U) & :=\int_{\mathbb{R}} a^{-\mb X}\frac{\mu}{v}\left|(u-\bar{u})_{\xi}\right|^{2} d \xi,\\
        \mb S_1(U)&:=-\int_{\mathbb{R}} a^{-\mathbf{X}} Q_1(u-\bar u)d\xi,  \\
        \mb S_2(U)&:=- \int_{\mathbb{R}} a^{-\mathbf{X}} Q_2 \left[ (1-\frac{\bar\theta}{\theta})+\frac{(\theta - \bar{\theta})^2}{\theta \bar{\theta}} \right]  d\xi.
    \end{align*}
    We decompose the functional $\mathbf{Y}(U)$ as 
    \begin{equation*}
        \mathbf{Y}(U):=\sum_{i=1}^{6} \mathbf{Y}_{i}(U),
    \end{equation*}
    where 
    \begin{align*}
        \mathbf{Y}_{1}(U) &  :=\int _{\mathbb{R}}a^{-\mb X}\left(u^{S}\right)_{\xi}^{-\mathbf X}(u-\bar{u})d\xi,\\
        \mathbf{Y}_{2}(U) & :=\int_{\mathbb{R}} a^{-\mathbf X} \frac{\left(v^{S}\right)_{\xi}^{-\mb X} \bar p^{th}}{\bar{v}}(v-\bar{v}) d \xi,\\
        \mathbf{Y}_{3}(U) & :=\int_{\mathbb{R}} a^{-\mathbf X} \frac{R}{\gamma-1} \frac{\left(\theta^{S}\right)_{\xi}^{-\mathbf X}}{\bar{\theta}}(\theta-\bar{\theta}) d \xi,\\
        \mathbf{Y}_{4}(U) & :=-\int _{\mathbb{R}} a^{-\mathbf X} R\left(\theta^{S}\right)_{\xi}^{-\mb X} \Phi\left(\frac{v}{\bar{v}}\right) d \xi,\\
        \mathbf{Y}_{5}(U) & :=-\int_{\mathbb{R}} a^{-\mathbf X} \frac{R}{\gamma-1}\left(\theta^{S}\right)_{\xi}^{-\mathbf X} \Phi\left(\frac{\theta}{\bar\theta}\right) d \xi,\\
        \mathbf{Y}_6(U) & := -\int_{\mathbb{R}} a_{\xi}^{-\mathbf X} \bar{\theta} \eta(U| \bar{U}) d\xi.            		
    \end{align*}           
	Notice from $\eqref{X(t)}$ that
    \begin{equation}\label{dot-X}
        \dot{\mathbf{X}}(t) = -\frac{M}{\delta_S} \left( \mathbf{Y}_1 + \mathbf{Y}_2 + \mathbf{Y}_3 \right),
     \end{equation}
    and so,
    \begin{equation}\label{de-X}
      \dot{\mathbf{X}}(t) \mathbf{Y} = -\frac{\delta_S}{M} | \dot{\mathbf{X}}(t) |^2 + \dot{\mathbf{X}}(t) \sum_{i=4}^6 \mathbf{Y}_i.
    \end{equation}
            	
     Then, for convenience of estimation, by \eqref{de-X}, we split the right-hand side of \eqref{d-a-1} into two parts as follows
    \begin{equation}\label{d-a-2}
        \begin{aligned} 
        \frac{d}{d t} \int_{\mathbb{R}} a^{-\mathbf{X}} \bar{\theta} \eta(U | \bar{U}) d \xi =
        & -\frac{\delta_{S}}{2 M}|\dot{\mathbf{X}}(t)|^{2}+\mathbf B_{1}+\mathbf B_{2}-\mathbf G-\frac{3}{4} \mathbf D \\ 
        & -\frac{\delta_{S}}{2 M}|\dot{\mathbf{X}}(t)|^{2}+\dot{\mathbf{X}}(t) \sum_{i=4}^{6} \mathbf Y_{i}  +\sum_{i=3}^{6} \mathbf B_{i}+\mb S_1+\mb S_2-\frac{1}{4} \mathbf D . 
        \end{aligned}
    \end{equation}
    \subsection{Leading order estimates}  		
    The following lemma estimates the first part of the right-hand side of \eqref{d-a-2}.
    \begin{lemma}\label{leading-order}
    There exists \(C>0\) such that 
    \begin{equation*}
       -\frac{\delta_{S}}{2 M}|\dot{\mathbf X}(t)|^{2}+\mathbf B_{1}+\mathbf B_{2}-\mathbf G-\frac{3}{4} \mathbf D  \leq-C \mathcal{G}^{S}(U)+ C \deltas^{4/3} \deltar^{4/3}  \left( e^{-C \deltas t}+|(1+\deltar t)|^{-10/3}\right),
    \end{equation*}
    where $\mathcal G^S(U)$ is defined in \eqref{Gs}.
	\end{lemma}
            		
    \begin{proof}
        We first rewrite the main terms in terms of the new variables \(y\) and \(w\): 
    \begin{equation}\label{omega}
        w:=u^{\mb X}-\bar{u}^{\mb X},
    \end{equation}
    and 
    \begin{equation*}
        y:=\frac{v^{S}(\xi)-v_m}{\delta_{S}}.
    \end{equation*}
    Since \(\mb X(t)\) is bounded on \([0, T]\) , it follows from \(\delta_{S}=v_+-v_m>0\) and \((v^{S})_\x>0\) that for any fixed \(t\), the change of variables \(\xi \in \mathbb{R} \mapsto y \in(0,1)\) is well defined, together with 
    \begin{equation}\label{dy}
        \frac{d y}{d \xi}=\frac{\left(v^{S}\right)_{\xi}}{\delta_{S}}>0.
    \end{equation}
    Note also that \(a(\xi)=1+\lambda y\) and \(a'(\xi)=\lambda(d y / d \xi)>0\). To perform the sharp estimates, we will use the \(O(1)\)-constants \(p_m\), \(\sigma_m\) defined in \eqref{p-sigma-}, which are indeed independent of the small constants $\delta_0$, $\deltar$, \(\delta_{S}\), \(\varepsilon_{1}\), since \(\frac{v_{+}}{2} \leq v_- \leq 2v_{+}\) and \(\frac{\theta_{+}}{2} \leq\theta_- \leq2\theta_{+}\).
            			
    \noindent$\bullet$ {\bf Estimates on \(\mathbf B_{2}-\mathbf G\):}

    First, we estimate the terms in the integrand‌ by $\eqref{rela-entro}$, 
    \begin{align*}
        & (u-\bar{u})(p-\bar p^{th})-\sigma \bar{\theta} \eta(U(t, \xi) | \bar{U}(t, \xi)) \\
        & \quad=(u-\bar{u})(p-\bar p^{th})-\sigma\left[R \bar{\theta} \Phi\left(\frac{v}{\bar{v}}\right)+\frac{R \bar{\theta}}{\gamma-1} \Phi\left(\frac{\theta}{\bar{\theta}}\right)+\frac{(u-\bar{u})^{2}}{2}\right].
    \end{align*}

    Since $|(v - \bar{v},\theta-\bar{\theta})| \leq C(\varepsilon_1+\delta_0)$ and $\Phi(1) = \Phi'(1) = 0, \Phi''(1) = 1$, we have
    \begin{equation}\label{phi-v}
        \Phi\left(\frac{v}{\bar{v}}\right) = \frac{(v - \bar{v})^2}{2\bar{v}^2} + O(|v - \bar{v}|^3)
    \end{equation}
    and 
    \begin{equation*}
        \Phi\left(\frac{\theta}{\bar{\theta}}\right) = \frac{(\theta - \bar{\theta})^2}{2\bar{\theta}^2} + O(|\theta - \bar{\theta}|^3).
    \end{equation*}
    Using \eqref{shock-base}-\eqref{theta-s} and
    \begin{equation}\label{bar-p}
        \begin{aligned}
        |\bar p^{th} - p_m| &\leq C(|\bar{v} - v_m| + |\bar{\theta} - \theta_m|)\\
        & \le C(|v^R - v_m|+ |(v^S)^{-\mb X} - v_m| + |\theta^R - \theta_m|+ |(\theta^S)^{-\mb X} - \theta_m| )\\
        &\leq C(\delta_S+\deltar)\le C\delta_0, 
        \end{aligned}
    \end{equation}
    together with \(p-\bar p^{th}=\frac{R}{v}(\theta-\bar{\theta})-\frac{\bar p^{th}}{v}(v-\bar{v})\), we have
    \begin{align*}
    	(u- & \bar{u})(p-\bar p^{th})-\sigma \bar{\theta} \eta(U(t, \xi) | \bar{U}(t, \xi)) \\
        = & (u-\bar{u})\left(\frac{R}{v}(\theta-\bar{\theta})-\frac{\bar p^{th}}{v}(v-\bar{v})\right)-\sigma\left[R \bar{\theta} \Phi\left(\frac{v}{\bar{v}}\right)+\frac{R \bar \theta}{\gamma-1} \Phi\left(\frac{\theta}{\bar{\theta}}\right)+\frac{(u-\bar{u})^{2}}{2}\right] \\
        \leq & (u-\bar{u})\left[\frac{R}{v_m}(\theta-\bar{\theta})-\frac{p_m}{v_m}(v-\bar{v})\right]-\frac{R \sigma_m \theta_m}{2\left(v_m\right)^{2}}(v-\bar{v})^{2}-\frac{R \sigma_m}{2(\gamma-1) \theta_m}(\theta-\bar{\theta})^{2} \\
        & -\frac{\sigma_m}{2}(u-\bar{u})^{2}+C\left(\delta_S+|(v-\bar{v},\theta-\bar\theta)|+\left|\bar{v}-v_m\right|+| \bar{\theta}-\theta_m|\right)|(v-\bar{v}, u-\bar{u}, \theta-\bar{\theta})|^{2} \\
        \leq & -\frac{R \sigma_m \theta_m}{2\left(v_m\right)^{2}}\left[(v-\bar{v})+\frac{u-\bar{u}}{\sigma_m}\right]^{2}-\frac{R \sigma_m}{2(\gamma-1) \theta_m}\left[(\theta-\bar{\theta})-\frac{(\gamma-1) \theta_m}{v_m \sigma_m}(u-\bar{u})\right]^{2} \\
        & +C\left(|(v-\bar{v},\theta-\bar\theta)|+\delta_{S}+|(v^R-v_m,\theta^R-\theta_m) |\right)|(v-\bar{v}, u-\bar{u}, \theta-\bar{\theta})|^{2},
    \end{align*}
    where the last inequality is obtained by using 
	\begin{equation*}
        \frac{R \theta_m}{2(v_m)^{2} \sigma_m}+\frac{R(\gamma-1) \theta_m}{2(v_m)^{2} \sigma_m}=\frac{R \gamma \theta_m}{2\left(v_m\right)^{2} \sigma_m}=\frac{\sigma_m}{2} \quad \text{by} \quad \sigma_m=\frac{\sqrt{\gamma R \theta_m}}{v_m}.
    \end{equation*}
            			
    Therefore, we have 
    \begin{equation}\label{B2-G}
        \mathbf B_{2}(U)-\mathbf G(U)\leq -\mathbf G_{1}(U)-\mathbf G_{2}(U)+\mathbf B_{\text{new}}(U),
    \end{equation}            			
    where             			
    \begin{align*}
        \mathbf G_{1}(U) & :=\frac{R\sigma _m\theta _m}{2(v_m)^{2}}\int _{\mathbb{R}}a_{\xi }^{-\mathbf X}\left[(v-\bar{v})+\frac{u-\bar{u}}{\sigma _m}\right]^{2}d\xi,\\
        \mathbf G_{2}(U) & :=\frac{R \sigma_m}{2(\gamma-1) \theta_m} \int_{\mathbb{R}} a_{\xi}^{-\mathbf X}\left[(\theta-\bar{\theta})-\frac{(\gamma-1) \theta_m}{v_m \sigma_m}(u-\bar{u})\right]^{2} d \xi,\\
        \mathbf B_{\text{new}}(U) & :=  C \delta_{S} \int_{\mathbb{R}} a_{\xi}^{-\mathbf X}|(v-\bar{v}, u-\bar{u}, \theta-\bar{\theta})|^{2} d \xi  +C \int_{\mathbb{R}} a_{\xi}^{-\mathbf X}|(v-\bar{v}, u-\bar{u}, \theta-\bar{\theta})|^{3} d \xi\\
        &\quad  + C  \int_\bbr a_\xi^{-\mb X}  |(v^R - v_m ,\theta^R - \theta_m)|  |(v-\bar v, u-\bar u,\theta-\bar\theta)|^2 d\x.
    \end{align*}            		
            			
    As in \eqref{B2-G}, a careful estimate gives two additional good terms, $\mathbf G_1$ and $\mathbf G_2$. The good term $\mathbf G_1$ is used to separate the term involving only $u$ from the terms involving both $v$ and $u$, while $\mathbf G_2$ allows us to make the same separation for the coupled terms in $\theta$ and $u$. The terms depending only on $u$ can then be controlled by the diffusion of $u$. Moreover, the new bad term $\mathbf B_{\mathrm{new}}$ can be controlled by $\mathbf G_1$, $\mathbf G_2$, $\mathbf D$, and $\mathcal G^{S}$, as follows.
            			
    First, using $\eqref{d-weight}$, we get
    \begin{equation*}
         \delta_{S} \int_{\mathbb{R}} a_{\xi}^{-\mathbf X}|(v-\bar{v}, u-\bar{u}, \theta-\bar{\theta})|^{2} d \xi \leq C \lambda \mathcal{G}^{S}.
    \end{equation*}
    Using $\eqref{d-weight}$, Sobolev's inequality, and \(\lambda \leq C \sqrt{\delta_{S}}\) from $\eqref{lamsmall}$, we have 
    \begin{align*}
        & \int_{\mathbb{R}} a_{\xi}^{-\mathbf X}\left|(v-\bar{v}, u-\bar{u}, \theta-\bar{\theta})\right|^{3} d \xi \\
    	& \leq C \int_{\mathbb{R}}|a_{\xi}^{-\mathbf X}|\left|(v-\bar{v})+\frac{u-\bar{u}}{\sigma_m}\right|^{3} d \xi+C \int_{\mathbb{R}}|a_{\xi}^{-\mathbf X}||u-\bar u|^{3} d \xi \\
        & \quad+C \int_{\mathbb{R}}|a_{\xi}^{-\mathbf X}|\left|(\theta-\bar{\theta})-\frac{(\gamma-1) \theta_m}{v_m \sigma_m}(u-\bar{u})\right|^{3} d \xi \\
        & \leq C  (\eps_1+\delta_0)\left(\mathbf G_{1}+\mathbf G_{2}\right)+C \frac{\lambda}{\delta_{S}} \int_{\mathbb{R}}|\left(v^{S}\right)_{\xi}^{-\mathbf X}|\|u-\bar u\|_{L^{\infty}(\mathbb{R})}^{2}|u-\bar u| d \xi \\
        & \leq C (\eps_1+\delta_0)\left(\mathbf G_{1}+\mathbf G_{2}\right)\\
        &\quad+C \frac{\lambda}{\delta_{S}}\|u-\bar u\|_{L^2}\|(u-\bar u)_\xi\|_{L^2} \sqrt{\int_{\mathbb{R}}|\left(v^{S}\right)_{\xi}^{-\mathbf X}| |u-\bar u|^{2} d \xi} \sqrt{\int_{\mathbb{R}}|\left(v^{S}\right)_{\xi}^{-\mathbf X}| d \xi}  \\
        & \leq C (\eps_1+\delta_0)\left(\mathbf G_{1}+\mathbf G_{2}\right)+C \varepsilon_{1} \frac{\lambda}{\sqrt{\delta_{S}}} \sqrt{\mathbf D} \sqrt{\mathcal{G}^{S}} \\
        & \leq C (\eps_1+\delta_0)\left(\mathbf G_{1}+\mathbf G_{2}+\mathbf D+\mathcal{G}^{S}\right),
    \end{align*}
    Likewise, using the interpolation inequality and Lemma \ref{lem:wave-inter}, we have
    \begin{align*}
        & \int_\bbr a_\xi^{-\mb X}  |(v^R - v_m, \theta^R - \theta_m)|  |(v-\bar v, u-\bar u,\theta-\bar\theta)|^2 d\x \\
    	&\quad \le C\delta_R  (\mb{G}_1+\mb{G}_2) +C \frac{\lam}{\deltas} \int_\bbr  |(v^S)_\xi^{-\mb X} |  |(v^R - v_m ,\theta^R - \theta_m)|   |u-\bar u|^2 d\xi  \\
        &\quad\le C\delta_R  (\mb{G}_1+\mb{G}_2)+C\frac{\lam}{\deltas} \|u-\bar u\|_{L^4(\bbr)}^2 \| |(v^S)_\xi^{-\mb X}| | (v^R - v_m ,\theta^R - \theta_m)|\|_{L^2(\bbr)}\\
        &\quad \le C \delta_0(\mb{G}_1+\mb{G}_2)+C\frac{1}{\sqrt\deltas} \|(u-\bar u)_\xi\|_{L^2(\bbr)}^{1/2}  \|u-\bar u\|_{L^2(\bbr)}^{3/2}  \deltas^{3/2}\deltar  \left(e^{-C \deltas t}+|(1+\deltar t)|^{-5/2}\right)\\
        &\quad\leq  C \delta_0(\mb{G}_1+\mb{G}_2)+C\mb D^{\frac14} \eps_1^{\frac32}\deltas \left(\deltar e^{-C \deltas t}+\deltar |(1+\deltar t)|^{-5/2}\right)\\
        &\quad\leq  C \delta_0(\mb{G}_1+\mb{G}_2)+C\eps_1 \mb D +C \deltas^{4/3} \deltar^{4/3}   \left(e^{-C \deltas t}+|(1+\deltar t)|^{-10/3}\right).
    \end{align*}
    Therefore, 
    \begin{equation}\label{Bnew}
        \mb B_{\text {new}} \le C (\lam+\eps_1+\delta_0) ( \mb{G}_1+\mb{G}_2 + \mb D +\mathcal{G}^S)+C \deltas^{4/3} \deltar^{4/3}   \left(e^{-C \deltas t}+|(1+\deltar t)|^{-10/3}\right).
    \end{equation}
            			
    \noindent$\bullet$ {\bf Estimates on \(-\frac{\delta_{S}}{2 M}|\dot{\mathbf X}(t)|^{2}\):}
    To estimate the term \(-\frac{\delta_{S}}{2 M}|\dot{\mathbf X}(t)|^{2}\), we will first estimate \(\mathbf Y_{1}\), \(\mathbf Y_{2}\), \(\mathbf Y_{3}\) from the definition of $\dot{\mathbf{X}}(t)$ by \eqref{dot-X}.
            			
    Using \eqref{omega}, \eqref{dy}, and the change of variables $\x\mapsto \x+\mb X(t)$, we have 
    \begin{equation*}
        \mathbf Y_{1}=\int_{\mathbb{R}} a^{-\mathbf X}\left(u^{S}\right)_{\xi}^{-\mathbf X}(u-\bar{u}) d \xi=\int_{\mathbb{R}} a\left(u^{S}\right)_{\xi}(u^{\mb X}-\bar{u}^{\mb X}) d \xi=-\delta_{S} \sigma \int_{0}^{1} a w d y.
    \end{equation*}    			
    Using \eqref{sigma} and \(|a-1| \leq\lambda\), we have 
    \begin{equation}\label{Y1}
        \left| \mathbf Y_{1}+\delta_{S} \sigma_m \int_{0}^{1} w d y \right| \leq C \delta_{S} \lambda \int_{0}^{1}|w| d y.
    \end{equation}
	For 
    \begin{align*}
        \mathbf Y_{2} & =\int_{\mathbb{R}} a^{-\mathbf X} \frac{\left(v^{S}\right)_{\xi}^{-\mathbf X} \bar p^{th}}{\bar{v}}(v-\bar{v}) d \xi \\
        & =-\int_{\mathbb{R}} a^{-\mathbf X} \frac{\left(v^{S}\right)_{\xi}^{-\mathbf X} \bar p^{th}}{\bar{v} \sigma_m}(u-\bar{u}) d \xi+\int_{\mathbb{R}} a^{-\mathbf X} \frac{\left(v^{S}\right)_{\xi}^{-\mathbf X} \bar p^{th}}{\bar{v}}\left[(v-\bar{v})+\frac{u-\bar{u}}{\sigma_m}\right] d \xi,\\
        & =-\int_{\mathbb{R}} a \frac{\left(v^{S}\right)_{\xi} \left(\bar p^{th}\right)^{\mb X}}{\bar{v}^{\mb X} \sigma_m}(u^{\mb X}-\bar{u}^{\mb X}) d \xi+\int_{\mathbb{R}} a^{-\mathbf X} \frac{\left(v^{S}\right)_{\xi}^{-\mathbf X} \bar p^{th}}{\bar{v}}\left[(v-\bar{v})+\frac{u-\bar{u}}{\sigma_m}\right] d \xi,
    \end{align*}
    using $\eqref{d-weight}$ and $\eqref{bar-p}$ , we have 
    \begin{equation}\label{Y2}
        \left| \mathbf Y_{2}+\frac{p_m \delta_{S}}{v_m \sigma_m} \int_{0}^{1} w d y \right| \leq C \delta_{S}\left(\lambda+\delta_{0}\right) \int_{0}^{1}|w| d y+C \frac{\delta_{S}}{\lambda} \int_{\mathbb{R}} a_{\xi}^{-\mathbf X}\left|(v-\bar{v})+\frac{u-\bar{u}}{\sigma_m}\right| d \xi.
    \end{equation}
    Likewise, since
    \begin{align*}
        \mathbf Y_{3}= & \int_{\mathbb{R}} a^{-\mathbf X} \frac{R}{\gamma-1} \frac{\left(\theta^{S}\right)_{\xi}^{-\mathbf X}}{\bar{\theta}}(\theta-\bar{\theta}) d \xi \\
        = & \frac{R \theta_m}{v_m \sigma_m} \int_{\mathbb{R}} a^{-\mathbf X} \frac{\left(\theta^{S}\right)_{\xi}^{-\mathbf X}}{\bar{\theta}}\left(u-\bar{u}\right) d \xi \\
        & +\int_{\mathbb{R}} a^{-\mathbf X} \frac{R}{\gamma-1} \frac{\left(\theta^{S}\right)_{\xi}^{-\mathbf X}}{\bar{\theta}}\left[(\theta-\bar{\theta})-\frac{(\gamma-1) \theta_m}{v_m \sigma_m}(u-\bar{u})\right] d \xi,\\
        = & \frac{R \theta_m}{v_m \sigma_m} \int_{\mathbb{R}} a \frac{\left(\theta^{S}\right)_{\xi}}{\bar{\theta}^{\mb X}}\left(u^{\mb X}-\bar{u}^{\mb X}\right) d \xi \\
        & +\int_{\mathbb{R}} a^{-\mathbf X} \frac{R}{\gamma-1} \frac{\left(\theta^{S}\right)_{\xi}^{-\mathbf X}}{\bar{\theta}}\left[(\theta-\bar{\theta})-\frac{(\gamma-1) \theta_m}{v_m \sigma_m}(u-\bar{u})\right] d \xi,
    \end{align*}
    using \(|\bar{\theta}-\theta_m| \leq C \delta_{0}\) and \eqref{theta-s}, $\eqref{d-weight}$, we have 
    \begin{equation}\label{Y3}
    	\begin{aligned}
        \left| \mathbf Y_{3}+\frac{(\gamma-1) p_m \delta_{S}}{v_m \sigma_m} \int_{0}^{1} w d y \right| \leq& C \delta_{S}\left(\lambda+\delta_{0}\right) \int_{0}^{1}|w| d y \\
        & +C \frac{\delta_{S}}{\lambda} \int_{\mathbb{R}} a_{\xi}^{-\mathbf X}\left|(\theta-\bar{\theta})-\frac{(\gamma-1) \theta_m}{v_m \sigma_m}(u-\bar{u})\right| d \xi.
        \end{aligned}
    \end{equation}
            			
    Therefore, using \eqref{dot-X}, \eqref{Y1}, \eqref{Y2} and  $\eqref{Y3}$ together with \(\sigma_m=\sqrt{\frac{\gamma p_m}{v_m}}\), we have 
    \begin{align*}
        \left| \dot{\mathbf{X}}(t)-2 \sigma_m M \int_{0}^{1} w d y \right| 
        &\leq C\left(\lambda+\delta_{0}+\varepsilon_{1}\right) \int_{0}^{1}|w| d y+\frac{C}{\lambda} \int_{\mathbb{R}} a_{\xi}^{-\mathbf{X}}\left|(v-\bar{v})+\frac{u-\bar{u}}{\sigma_m}\right| d \xi \\
        & \quad+\frac{C}{\lambda} \int_{\mathbb{R}} a_{\xi}^{-\mathbf{X}}\left|(\theta-\bar{\theta})-\frac{(\gamma-1) \theta_m}{v_m \sigma_m}(u-\bar{u})\right| d \xi,
    \end{align*}
    which yields 
    \begin{align*}
        \left( \left|2 \sigma_m M \int_{0}^{1} w d y\right|-|\dot{\mathbf{X}}(t)|\right)^{2} 
        &\leq C\left(\lambda+\delta_{0}+\varepsilon_{1}\right)^{2} \int_{0}^{1}|w|^{2} d y \\
        & \quad+\frac{C}{\lambda^{2}}\left(\mathbf G_{1}+\mathbf G_{2}\right) \int_{\mathbb{R}} a_{\xi}^{-\mathbf{X}} d \xi.
    \end{align*}
    This estimate and the algebraic inequality \(\frac{p^{2}}{2}-q^{2} \leq(p-q)^{2}\) for all \(p, q \geq0\) imply 
    \begin{equation*}
        2\left(\sigma_m\right)^{2} M^{2}\left(\int_{0}^{1} w d y\right)^{2}-|\dot{\mathbf{X}}(t)|^{2} \leq C\left(\lambda+\delta_{0}+\varepsilon_{1}\right)^{2} \int_{0}^{1}|w|^{2} d y+\frac{C}{\lambda}\left(\mathbf G_{1}+\mathbf G_{2}\right).
    \end{equation*}
            			
    Thus, 
    \begin{equation}\label{es-dotX}
        -\frac{\delta_{S}}{2 M}|\dot{\mathbf{X}}(t)|^{2} \leq-\left(\sigma_m\right)^{2} M \delta_{S}\left(\int_{0}^{1} w d y\right)^{2}+C \delta_{S}\left(\lambda+\delta_{0}+\varepsilon_{1}\right)^{2} \int_{0}^{1} w^{2} d y+\frac{C \delta_{S}}{\lambda}\left(\mathbf G_{1}+\mathbf G_{2}\right).
    \end{equation}
            			
    \noindent$\bullet$ {\bf Estimates on \(\mathbf B_{1}\):}
    Set 
    \begin{equation}\label{de-B1}
        \begin{aligned}
        \mathbf B_{1}= & \int_{\mathbb{R}} a^{-\mathbf{X}}\left(v^{S}\right)_{\xi}^{-\mathbf{X}}\left[\frac{(\gamma+1) p_m \sigma_m}{2\left(v_m\right)^{2}}(v-\bar{v})^{2}+\frac{R \sigma_m}{2 v_m \theta_m}(\theta-\bar{\theta})^{2}-\frac{p_m \sigma_m}{v_m \theta_m}(v-\bar{v})(\theta-\bar{\theta})\right] d \xi \\
        := & \mathbf B_{11}(U)+\mathbf B_{12}(U)+\mathbf B_{13}(U).
        \end{aligned}
    \end{equation}
    First, using Young's inequality, we have
    \begin{align*}
        \mathbf B_{11}= & \frac{(\gamma+1) p_m \sigma_m}{2\left(v_m\right)^{2}} \int_{\mathbb{R}} a^{-\mathbf{X}}\left(v^{S}\right)_{\xi}^{-\mathbf{X}}\left|\left((v-\bar{v})+\frac{u-\bar{u}}{\sigma_m}\right)-\frac{u-\bar{u}}{\sigma_m}\right|^{2} d \xi \\
        \leq & \frac{(\gamma+1) p_m \sigma_m}{2\left(v_m\right)^{2}}\left(1+\lambda+\left(\frac{\delta_{S}}{\lambda}\right)^{\frac{1}{2}}\right) \int_{\mathbb{R}}\left(v^{S}\right)_{\xi}^{-\mathbf{X}}\left|\frac{u-\bar{u}}{\sigma_m}\right|^{2} d \xi \\
        & +C\left(\frac{\delta_{S}}{\lambda}\right)^{-\frac{1}{2}} \int_{\mathbb{R}}\left(v^{S}\right)_{\xi}^{-\mathbf{X}}\left|(v-\bar{v})+\frac{u-\bar{u}}{\sigma_m}\right|^{2} d \xi,
    \end{align*}
    which, together with $\eqref{d-weight}$ and the change of variables $\x\mapsto\x+\mb X(t)$, yields 
    \begin{equation}\label{B11}
        \mathbf B_{11} \leq \frac{(\gamma+1) p_m}{2\left(v_m\right)^{2} \sigma_m}\left(1+\lambda+\left(\frac{\delta_{S}}{\lambda}\right)^{\frac{1}{2}}\right) \delta_{S} \int_{0}^{1} w^{2} d y+C\left(\frac{\delta_{S}}{\lambda}\right)^{\frac{1}{2}} \mathbf G_{1}.
	\end{equation}
    Similarly, using \eqref{eq:sp-pm} and the change of variables $\x\mapsto\x+\mb X(t)$, we estimate that 
    \begin{equation}\label{B12}
        \begin{aligned}
        \mathbf B_{12} & =\frac{R \sigma_m}{2 v_m \theta_m} \int_{\mathbb{R}} a^{-\mathbf{X}}\left(v^{S}\right)_{\xi}^{-\mathbf{X}}\left|\left((\theta-\bar{\theta})-\frac{(\gamma-1) \theta_m}{v_m \sigma_m}(u-\bar{u})\right)+\frac{(\gamma-1) \theta_m}{v_m \sigma_m}(u-\bar{u})\right|^{2} d \xi \\
        & \leq \frac{(\gamma-1)^{2} p_m}{2\left(v_m\right)^{2} \sigma_m}\left(1+\lambda+\left(\frac{\delta_{S}}{\lambda}\right)^{\frac{1}{2}}\right) \delta_{S} \int_{0}^{1} w^{2} d y+C\left(\frac{\delta_{S}}{\lambda}\right)^{\frac{1}{2}} \mathbf G_{2}
        \end{aligned}
    \end{equation}
    and 
    \begin{equation}\label{B13}
        \begin{aligned}
    	\mathbf B_{13} \leq & \frac{p_m \sigma_m}{v_m \theta_m} \int_{\mathbb{R}} a^{-\mathbf{X}}\left(v^{S}\right)_{\xi}^{-\mathbf{X}}\left|\left((v-\bar{v})+\frac{u-\bar{u}}{\sigma_m}\right)-\frac{u-\bar{u}}{\sigma_m}\right| \\
        & \cdot\left|\left((\theta-\bar{\theta})-\frac{(\gamma-1) \theta_m}{v_m \sigma_m}(u-\bar{u})\right)+\frac{(\gamma-1) \theta_m}{v_m \sigma_m}(u-\bar{u})\right| d \xi \\
        & \leq \frac{(\gamma-1) p_m}{\left(v_m\right)^{2} \sigma_m}\left(1+\lambda+\left(\frac{\delta_{S}}{\lambda}\right)^{\frac{1}{2}}\right) \delta_{S} \int_{0}^{1} w^{2} d y+C\left(\frac{\delta_{S}}{\lambda}\right)^{\frac{1}{2}}\left(\mathbf G_{1}+\mathbf G_{2}\right).
        \end{aligned}
    \end{equation}
            			
    Let \(\alpha_m\) be the \(O(1)\)-constant defined by 
    \begin{equation*}
        \alpha_m:=\frac{\gamma(\gamma+1) p_m}{2\left(v_m\right)^{2} \sigma_m}.
    \end{equation*}        			
    Substituting $\eqref{B11}$, $\eqref{B12}$ and $\eqref{B13}$ into $\eqref{de-B1}$, we have 
    \begin{equation}\label{es-B1}
        \mathbf B_{1} \leq \alpha_m\left(1+\lambda+\left(\frac{\delta_{S}}{\lambda}\right)^{\frac{1}{2}}\right) \delta_{S} \int_{0}^{1} w^{2} d y+C\left(\frac{\delta_{S}}{\lambda}\right)^{\frac{1}{2}}\left(\mathbf G_{1}+\mathbf G_{2}\right).
    \end{equation}
            			
    \noindent$\bullet$ {\bf Estimates on \(\mathbf D\):}
    Set 
    \begin{equation*}
        \mathbf D(U)=\int_{\mathbb{R}} a^{-\mathbf{X}}\left(\frac{\mu}{v}\left|(u-\bar{u})_{\xi}\right|^{2}\right) d \xi.
    \end{equation*}            			
    First, using \(a \geq1\) and the change of variable, we have 
    \begin{equation*}
        \mathbf D \geq \int_{\mathbb{R}} \frac{\mu}{v}\left|(u-\bar{u})_{\xi}\right|^{2} d \xi=\int_{0}^{1}\left|\partial_{y} w\right|^{2} \frac{\mu}{v^{\mb X}}\left(\frac{d y}{d \xi}\right) d y.
    \end{equation*}
    To estimate the right-hand side, set \(p_{+}:=p(v_{+}, \theta_{+})\) for simplicity. Then, using $\eqref{shock-v-u-theta-1}_1$, we have 
    \begin{equation*}
        \mu \frac{v_{\xi}^{S}}{v^{S}}=-\sigma\left(v^{S}-v_m\right)-\frac{p^{S}-p_m}{\sigma},
    \end{equation*}
    and so, 
    \begin{align*}
        \delta_{S} \frac{\mu}{v^{S}} \frac{d y}{d \xi} & =\mu \frac{v_{\xi}^{S}}{v^{S}}=-\sigma\left(v^{S}-v_m\right)-\frac{p^{S}-p_m}{\sigma} \\
        & =-\frac{1}{\sigma}\left[\sigma^{2}\left(v^{S}-v_m\right)+p^{S}-p_m\right],
    \end{align*}
    which together with \(\sigma^{2}=\frac{p_m-p_{+}}{v_{+}-v_m}\) yields 
    \begin{align*}
        \delta_{S} \frac{\mu}{v^{S}} \frac{d y}{d \xi}= & -\frac{1}{\sigma\left(v_{+}-v_m\right)}\left[\left(p_m-p_{+}\right)\left(v^{S}-v_m\right)+\left(v_{+}-v_m\right)\left(p^{S}-p_m\right)\right] \\
        = & -\frac{1}{\sigma\left(v_{+}-v_m\right)}\left[\left(p^{S}-p_{+}\right)\left(v^{S}-v_m\right)+\left(v^{S}-v_m\right)\left(p_m-p^{S}\right)\right. \\
        & \left.+\left(v^{S}-v_m\right)\left(p^{S}-p_m\right)+\left(v_{+}-v^{S}\right)\left(p^{S}-p_m\right)\right] \\
        = & -\frac{1}{\sigma\left(v_{+}-v_m\right)}\left[\left(p^{S}-p_{+}\right)\left(v^{S}-v_m\right)+\left(v_{+}-v^{S}\right)\left(p^{S}-p_m\right)\right].
    \end{align*}
            			
    Since \(y=\frac{v^{S}-v_m}{\delta_{S}}\) and \(1-y=\frac{v_{+}-v^{S}}{\delta_{S}}\), 
    \begin{equation*}
        \frac{1}{y(1-y)} \frac{\mu}{v^{S}} \frac{d y}{d \xi}=\frac{1}{\sigma}\left(\frac{p^{S}-p_{+}}{v^{S}-v_{+}}-\frac{p^{S}-p_m}{v^{S}-v_m}\right).
    \end{equation*}
            			
    By Lemma \ref{sharp-D}, we have 
    \begin{equation*}
        \left| \frac{1}{y(1-y)} \frac{\mu}{v^{S}} \frac{d y}{d \xi}-\delta_{S} \alpha_m\right| \leq C \delta_{S}^{2}, 
    \end{equation*}
            			
    This together with \(|v^{\mb X}-v^{S}| \leq C(\varepsilon_{1}+\delta_{0})\) implies 
    \begin{align*}
        \mathbf D & \geq \int_{0}^{1}\left|\partial_{y} w\right|^{2} \frac{\mu}{v^{S}}\left(\frac{d y}{d \xi}\right) d y+\int_{0}^{1}\left|\partial_{y} w\right|^{2}\left(\frac{v^{S}}{v^{\mb X}}-1\right) \frac{\mu}{v^{S}}\left(\frac{d y}{d \xi}\right) d y \\
        & \geq \alpha_m \left(1-C\left(\delta_0+\varepsilon_{1}\right)\right) \delta_{S} \int_{0}^{1} y(1-y)\left|\partial_{y} w\right|^{2} d y.
    \end{align*}

    Our intention is to use the Poincar\'e-type inequality of Lemma \ref{poincare} to absorb the main bad term \(\mathbf B_{1}\) by the diffusion term \(\mathbf D\) with $\eqref{es-dotX}$. 
    First, using Lemma \ref{poincare} and 
    \begin{equation*}
        \int_{0}^{1}|w-\bar{w}|^{2} d y=\int_{0}^{1} w^{2} d y-\bar{w}^{2}, \quad \bar{w}:=\int_{0}^{1} w d y,
    \end{equation*}
    we have 
    \begin{equation}\label{D}
        \begin{aligned}
        \mathbf D \geq & 2 \alpha_m \left(1-C\left(\delta_0+\varepsilon_{1}\right)\right) \delta_{S}\left[\int_{0}^{1} w^{2} d y-\bar{w}^{2}\right].
        \end{aligned}
    \end{equation}
            			
    \noindent$\bullet$ {\bf  Conclusion:}
    First, combining $\eqref{B2-G}$, $\eqref{Bnew}$, $\eqref{es-dotX}$, $\eqref{es-B1}$, $\eqref{D}$ and using the smallness of \(\frac{\delta_{S}}{\lambda}\), \(\lambda\) , and \(\delta_0\), \(\varepsilon_{1}\), we have 
    \begin{align*}
        & -\frac{\delta_{S}}{2 M}|\dot{\mathbf{X}}(t)|^{2}+\mathbf B_{1}+\mathbf B_{2}-\mathbf G-\frac{3}{4} \mathbf D \\
        & \leq C\left(\lambda+\varepsilon_{1}+\delta_{0}\right) \mathcal{G}^{S} -\frac{\alpha_m}{4} \delta_{S} \int_{0}^{1} w^{2} d y-\frac{1}{2}\left(\mathbf G_{1}+\mathbf G_{2}\right) \\
        &\quad -\left(\sigma_m\right)^{2} M \delta_{S} \bar{w}^{2}+\frac{3}{2} \alpha_m \delta_{S} \bar{w}^{2}+C \deltas^{4/3} \deltar^{4/3}  \left( e^{-C \deltas t}+|(1+\deltar t)|^{-10/3}\right).
    \end{align*}
    Choosing the specific \(O(1)\)-constant \(M\) as in $\eqref{M}$, that is, 
    \begin{equation*}
        M=\frac{3}{2\left(\sigma_m\right)^{2}} \alpha_m ,
    \end{equation*}
    we have 
    \begin{align*}
        & -\frac{\delta_{S}}{2 M}|\dot{\mathbf{X}}(t)|^{2}+\mathbf B_{1}+\mathbf B_{2}-\mathbf G-\frac{3}{4} \mathbf D \\
        & \leq C\left(\lambda+\varepsilon_{1}+\delta_{0}\right) \mathcal{G}^{S} -\frac{\alpha_m}{4} \delta_{S} \int_{0}^{1} w^{2} d y-\frac{1}{2}\left(\mathbf G_{1}+\mathbf G_{2}\right)\\
        &\quad+C \deltas^{4/3} \deltar^{4/3}  \left( e^{-C \deltas t}+|(1+\deltar t)|^{-10/3}\right).
    \end{align*}
            			
    Finally, using
    \begin{align*}
        &\int_{\mathbb{R}}|(v^{S})_{\xi}^{-\mathbf{X}}|(\theta-\bar{\theta})^{2} d \xi \\
        & \leq 2 \int_{\mathbb{R}}|(v^{S})_{\xi}^{-\mathbf{X}}|\left|(\theta-\bar{\theta})-\frac{(\gamma-1) \theta_m}{v_m \sigma_m}(u-\bar{u})\right|^{2} d \xi +C \int_{\mathbb{R}}|(v^{S})_{\xi}^{-\mathbf{X}}|(u-\bar{u})^{2} d \xi \\
        & \leq C \frac{\delta_{S}}{\lambda} \mathbf G_{2}+C \int_{\mathbb{R}}|(v^{S})_{\xi}^{-\mathbf{X}}|(u-\bar{u})^{2} d \xi,
    \end{align*}
    and
    \begin{align*}
        &\int_{\mathbb{R}}|(v^{S})_{\xi}^{-\mathbf{X}}|(v-\bar{v})^{2} d \xi \\
        & \leq 2 \int_{\mathbb{R}}\left|(v^{S})_{\xi}^{-\mathbf{X}}\right|\left|(v-\bar{v})+\frac{u-\bar{u}}{\sigma_m}\right|^{2} d \xi  +2 \int_{\mathbb{R}}|(v^{S})_{\xi}^{-\mathbf{X}}|\left(\frac{u-\bar{u}}{\sigma_m}\right)^{2} d \xi \\
        & \leq C \frac{\delta_{S}}{\lambda} \mathbf G_{1}+C \int_{\mathbb{R}}|(v^{S})_{\xi}^{-\mathbf{X}}|(u-\bar{u})^{2} d \xi,
    \end{align*}
    and
    \begin{equation*}
        \delta_{S} \int_{0}^{1} w^{2} d y=\int_{\mathbb{R}}|(v^{S})_{\xi}||u^{\mb X}-\bar{u}^{\mb X}|^{2} d \xi=\int_{\mathbb{R}}|(v^{S})_{\xi}^{-\mathbf{X}}||u-\bar{u}|^{2} d \xi,
    \end{equation*}            		
    we have
    \begin{align*}
        & -\frac{\delta_{S}}{2 M}|\dot{\mathbf{X}}(t)|^{2}+\mathbf B_{1}+\mathbf B_{2}-\mathbf G-\frac{3}{4} \mathbf D \\ 
        & \leq-C \int_{\mathbb{R}}|(v^{S})_{\xi}^{-\mathbf{X}}||(v-\bar{v}, u-\bar{u}, \theta-\bar{\theta})|^{2} d \xi+C \deltas^{4/3} \deltar^{4/3}  \left( e^{-C \deltas t}+|(1+\deltar t)|^{-10/3}\right), 
    \end{align*} 			            			
    which implies the desired estimate in the lemma.
    \end{proof}   			
            			
    \subsection{Proof of Lemma $\ref{Re-en}$}
            			
    In order to prove Lemma \ref{Re-en}, we have to estimate the right-hand side of \eqref{d-a-2} and we already get the leading order estimates in Lemma \ref{leading-order}, so we estimate the remaining terms of \eqref{d-a-2} as follows.
    \begin{lemma}\label{lem-sharp2}
    There exists $C>0$ such that 
    \begin{align*}
        &\frac{C}{\deltas}\sum_{i=4}^6 |\mb{Y}_i|^2  +\sum_{i=3}^6 \mb B_i +\mb S_1 +\mb S_2 -\frac{1}{8} \mb D \\
        &\leq -C \mathcal{G}^R(U)+ C(\delta_0+\eps_1)  \mathcal{G}^S(U)\\
        &\quad+C \deltas\deltar \left(e^{-C \deltas t}+|(1+\deltar t)|^{-\frac{5}{2}}\right)\big(\|u-\bar u\|_{L^2(\bbr)}+\|\theta-\bar \theta\|_{L^2(\bbr)}\big)\\
        &\quad + C(1+t)^{-\frac{9}{8}}\|(u-\bar u, \theta-\bar\theta)\|_{L^2(\bbr)}^2+C\deltar^{\frac{1}{4}}(1+t)^{-\frac{9}{8}},
    \end{align*}
    where the good terms $\mathcal{G}^S(U)$ and $\mathcal{G}^R(U)$ are as in \eqref{Gs} and \eqref{Gr}.
    \end{lemma}
    \begin{proof}
    In the proof, we need to estimate every remaining term on the right-hand side of $\eqref{d-a-2}$. 
            	       
    \noindent$\bullet$ {\bf Estimates on \(\mathbf B_{i}(i=3,4) \):}\\
    For $\mb B_3$, using Young's inequality and \eqref{a-prime}, we have the following estimates
    \begin{align*} 
            \mb B_{3} & \leq C \int_{\mathbb{R}}|a_{\xi}^{-\mathbf{X}}||u-\bar{u}|\left(\left|(u-\bar{u})_{\xi}\right|+\left|\bar{u}_{\xi}\right||v-\bar{v}|\right) d \xi \\ 				 
            & \leq \frac{1}{80} \mathbf D+C \lambda^2 \int_{\mathbb{R}}|(v^{S})_{\xi}^{-\mathbf{X}}|\left|(u-\bar{u})\right|^{2} d \xi+C \int_{\mathbb{R}}|\bar{u}_{\xi}|^{2}|v-\bar{v}|^{2} d \xi \\ 
            & \leq \frac{1}{80} \mathbf D+C \delta_0 \mathcal{G}^{S}+C \int_{\mathbb{R}}\left|\bar{u}_{\xi}\right|^{2} |v-\bar{v}|^{2} d \xi\\
            & \leq \frac{1}{80} \mathbf D+C \delta_0 \mathcal{G}^{S}+C \delta_0 \mathcal{G}^{R}. 
    \end{align*}
     As for $\mb B_4$, similarly, we obtain
	\begin{align*} 
    \mathbf B_{4} &\leq C \int_{\mathbb{R}}\left|\bar{u}_{\xi}\right||v-\bar{v}|\left|(u-\bar{u})_{\xi}\right|d \xi +C \int_{\mathbb{R}}|\theta-\bar{\theta}|^{2}\left[\left|\left(\bar{v}_{\xi}, \bar{\theta}_{\xi}, \bar{u}_{\xi}\right)\right|^{2}\right] d \xi\\ 
    & \quad +C \int_{\mathbb{R}}|\theta-\bar{\theta}|\left(\left|(u-\bar{u})_{\xi}\right|^{2}+\left|\bar{u}_{\xi}\right|\left|(u-\bar{u})_{\xi}\right|+\left|\bar{u}_{\xi}\right|^{2}|v-\bar{v}|\right) d \xi \\ 
    & \leq  \frac{1}{40} \mathbf D+C \delta_0 \mathcal{G}^{S}+C \delta_0 \mathcal{G}^{R},
	\end{align*}
  where we use the smallness of $\delta_0$ and $\varepsilon_{1}$.
  By combining the above two estimates, we obtain
  \begin{equation} \label{sum-B}         		
    \sum_{i=3}^{4} \mathbf B_{i} \leq \frac{1}{20} \mathbf D+C\delta_0(\mathcal{G}^{S}+\mathcal{G}^{R}).      
    \end{equation}
 \noindent$\bullet$ {\bf Estimates on \(\mathbf B_i(i=5,6) \):}\\
    We can handle $\mb B_5(U)$ as follows. Set 
    \begin{align*}
     J_2 =  R\left((\theta^R_t-\s \theta^R_\x)- \sigma (\theta^S)^{-\mathbf{X}} _{\xi}\right)\Phi\left( \frac{v}{\bar{v}} \right) - \frac{\bar p^{th}\bar u_\x}{v \bar{v}} (v - \bar{v})^2.
    \end{align*}
    Using \eqref{shock-base}, \eqref{shock-vu}, \eqref{theta-s}, \eqref{phi-v} and \eqref{bar-p}, we have
    \begin{align*}
    J_{2} &\leq \frac{(\gamma + 1)p_m\sigma_m}{2(v_m)^2} (v^S)_\xi^{-\mathbf{X}} (v - \bar{v})^2 + C(\delta_0 + \varepsilon_1)|(v^S)_{\xi}^{-\mathbf{X}}||v - \bar{v}|^2\\
    &\quad -\frac{(\gamma+1)\bar p^{th}}{2\bar v^2}u^R_\xi (v-\bar v)^2 + C(\delta_0 +\eps_1 ) |u^R_\xi| |v-\bar v|^2,		
    \end{align*}
    which means
    \begin{equation}\label{J2}
        \begin{aligned}
        \mb B_5=&\int_{\mathbb{R}} a^{-\mathbf{X}}J_{2} d \xi-\int_{\mathbb{R}} a^{-\mathbf{X}} \frac{(\gamma + 1)p_m\sigma_m}{2(v_m)^2} (v^S)_\xi^{-\mathbf{X}} (v - \bar{v})^2 d\x \\
        \leq
        &C\left(\delta_0+\varepsilon_{1}\right) \int_{\mathbb{R}} a^{-\mathbf{X}}\left( |(v^{S})_{\xi}^{-\mathbf{X}}|+|u^R_\xi|  \right)|v-\bar{v}|^{2} d \xi- \int_{\mathbb{R}} a^{-\mathbf{X}} \frac{(\gamma+1)\bar p^{th}}{2\bar v^2}u^R_\xi (v-\bar v)^2 d \xi.
        \end{aligned}
    \end{equation}
    Similarly, for $\mb B_6$, set 
    \begin{align*}
        J_3=\frac{R}{\gamma-1}\Big((\theta_t^R-\sigma \theta_\x^R)-\sigma (\theta^S)^{-\mb X}_\xi \Big)\Phi\left(\frac{\theta}{\bar \theta}\right) - \frac{\bar u_\xi}{\theta}(\theta-\bar \theta)(p-\bar p^{th}) +\bar p^{th} \bar u_\xi \frac{(\theta-\bar\theta)^2}{\theta\bar\theta}.
    \end{align*}
    Then, we have 
    \begin{equation}\label{J3}
        \begin{aligned}
        \mb B_6&=	\int_{\mathbb{R}} a^{-\mathbf{X}}J_{3} d \xi -\int_{\mathbb{R}} a^{-\mathbf{X}} \frac{\sigma_m\left(v^{S}\right)_{\xi}^{-\mathbf{X}}}{2 v_m \theta_m}\left(R(\theta-\bar{\theta})-2 p_m(v-\bar{v})\right)(\theta-\bar{\theta}) d\x \\
        &\leq C\left(\delta_0+\varepsilon_{1}\right) \int_{\mathbb{R}} a^{-\mathbf{X}}\left( |(v^{S})_{\xi}^{-\mathbf{X}}|+|u^R_\xi|  \right)|(v-\bar{v}, \theta-\bar{\theta})|^{2} d \xi\\
        &\quad + \int_{\mathbb{R}} a^{-\mathbf{X}}u^R_\x \Big[ -\frac{R}{2\bar v\bar\theta}(\theta-\bar \theta) + \frac{\bar p^{th}}{\bar v \bar \theta} (v-\bar v)\Big] (\theta-\bar \theta) d \xi.
        \end{aligned}
    \end{equation}
    Therefore, by combining the estimates \eqref{J2} and \eqref{J3} together with $\bar p^{th}=\frac{R \bar{\theta}}{\bar{v}}$, 
    \begin{align*}
        \mb B_5(U)+\mb B_6(U) &=
        \int_{\mathbb{R}} a^{-\mathbf{X}}\left(J_{2}+J_{3}\right) d \xi -\mb B_1(U)\\ 
        &\leq C\left(\delta_0+\varepsilon_{1}\right) \int_{\mathbb{R}} a^{-\mathbf{X}}\left( |(v^{S})_{\xi}^{-\mathbf{X}}|+|u^R_\xi|  \right)|(v-\bar{v}, \theta-\bar{\theta})|^{2} d \xi\\
        &\quad -  \int_\bbr \!a^{-\mb X}   \frac{R u^R_\xi }{2\bar v^2}\Big[\frac{\gamma \bar\theta}{\bar v}(v-\bar v)^2+\Big(\sqrt{\frac{\bar\theta}{\bar v}}(v-\bar v) - \sqrt{\frac{\bar v}{\bar\theta}}(\theta-\bar\theta)\Big)^2\Big] d\xi\\
        &:= C\left(\delta_0+\varepsilon_{1}\right) \int_{\mathbb{R}} a^{-\mathbf{X}}\left( |(v^{S})_{\xi}^{-\mathbf{X}}|+|u^R_\xi|  \right)|(v-\bar{v}, \theta-\bar{\theta})|^{2} d \xi + Z.
    \end{align*}
    We now derive the simpler form $\mathcal{G}^R$ from the above good term $Z$ as follows. Using
    \[
    2(v-\bar v)(\theta-\bar\theta) \le  \frac{ 2\bar\theta}{\bar v} (v-\bar v)^2 + \frac{\bar v}{2\bar\theta}(\theta-\bar\theta)^2,
    \]
    and $u^R_\xi \sim v^R_\xi >0$ by Lemma \ref{lem:rare-v} (1), 
    we have
    \[
    -Z  \ge \int_\bbr \!a^{-\mb X}   \frac{R u^R_\xi }{2\bar v^2}\Big[\frac{(\gamma-1) \bar\theta}{\bar v}(v-\bar v)^2 + \frac{\bar v}{2\bar\theta}(\theta-\bar\theta)^2 \Big] d\xi  \ge C \mathcal{G}^R(U).
    \]
    Therefore, using the smallness of $\delta_{0}$ and $\eps_1$, we have
    \begin{equation}\label{R}
    \mb B_5(U)+\mb B_6(U) \leq C (\delta_0+\varepsilon_{1}) \mathcal{G}^S(U)-C\mathcal{G}^R(U).
    \end{equation}
            			
    \noindent$\bullet$ {\bf Estimates on \(\mathbf Y_{i}(i=4,5,6) \):}\\
    Using \eqref{theta-s} and $\eqref{phi-v}$, we have
     \begin{equation*}
         \left|\left(\mathbf Y_{4}, \mathbf Y_{5}\right)\right| \leq C \int_{\mathbb{R}}|(v^{S})_{\xi}^{-\mathbf{X}}| |(v-\bar{v}, \theta-\bar{\theta})|^{2} d \xi.
     \end{equation*}           	
    In addition, Lemmas \ref{lem:shock-v} and \ref{lem:variables-equivalence} yield
    \begin{equation*}
    \left|\left(\mathbf Y_{4}, \mathbf  Y_{5}\right)\right| \leq C\delta_{S}^2 \int_{\mathbb{R}}|(v-\bar{v}, \theta-\bar{\theta})|^{2} d \xi \leq C\delta_{S}^{2} (\varepsilon_{1}+\delta_0)^{2},
    \end{equation*}
    which means             			
    \begin{equation*}
        \frac{C}{\delta_{S}}\left|\left(\mathbf  Y_{4}, \mathbf  Y_{5}\right)\right|^{2} \leq C \delta_{S} (\varepsilon_{1}+\delta_0)^{2} \mathcal{G}^{S}.
    \end{equation*}            			
    For $\mb Y_6$, similarly, we have
    \begin{align*} 
        \frac{C}{\delta_{S}}\left|\mathbf  Y_{6}\right|^{2} & \leq \frac{C}{\delta_{S}}\left(\int_{\mathbb{R}}|a_{\xi}^{-\mathbf{X}}||(v-\bar{v}, u-\bar{u}, \theta-\bar{\theta})|^{2} d \xi\right)^{2} \\ 
        & \leq \frac{C \lambda^{2}}{\delta_{S}^{3}}\left(\int_{\mathbb{R}}|(v^{S})_{\xi}^{-\mathbf{X}}||(v-\bar{v}, u-\bar{u}, \theta-\bar{\theta})|^{2} d \xi\right)^{2} \\ 
        & \leq \frac{C \lambda^{2}}{\delta_{S}}\left\| (v-\bar{v}, u-\bar{u}, \theta-\bar{\theta})\right\| _{L^{2}(\mathbb{R})}^{2} \int_{\mathbb{R}}|(v^{S})_{\xi}^{-\mathbf{X}}||(v-\bar{v}, u-\bar{u}, \theta-\bar{\theta})|^{2} d \xi \\ 
        & \leq C (\varepsilon_{1}+\delta_0)^{2} \int_{\mathbb{R}}|(v^{S})_{\xi}^{-\mathbf{X}}||(v-\bar{v}, u-\bar{u}, \theta-\bar{\theta})|^{2} d \xi \leq C (\varepsilon_{1}+\delta_0)^{2} \mathcal{G}^{S}. 
    \end{align*}
    Combining the above estimates, we obtain
    \begin{equation}\label{es-Y}
        \frac{C}{\delta_{S}} \sum_{i=4}^{6}\left|\mathbf  Y_{i}\right|^{2} \leq C (\varepsilon_{1}+\delta_0) \mathcal{G}^{S}.
    \end{equation}
    \noindent$\bullet$ {\bf Estimates on the terms $\mb S_i~(i=1,2)$:}\\
    First, by \eqref{QR}, for $i=1,2,$ 
    \begin{align*}
     |Q^R_i| &\le C \Big[|u^R_{\xi\xi}|+|u^R_{\xi}||v^R_{\xi}|+|u^R_{\xi}|^2\Big]\\
     &\le C  \Big[\big|u^R_{\xi\xi}\big|+\big|\big(v^R_{\xi},u^R_{\xi}\big)\big|^2\Big],
    \end{align*}
    and then, using Sobolev's inequality and Lemma \ref{lem:variables-equivalence}, we have
    \begin{align*}
        \mb S_1+\mb S_2&=-\int_\bbr a^{\mb {-X}}Q_1(u-\bar u) d\x- \int_\bbr  a^{\mb {-X}}Q_2\left[ (1-\frac{\bar \theta}{\theta})+\frac{(\theta-\bar\theta)^2}{\theta\bar\theta} \right]d\x \\
        &\leq C\int_\bbr \big(|Q^I_1|+|Q_1^R|\big)|u-\bar u| d\x+ C\int_\bbr  \big(|Q^I_2|+|Q_2^R|\big)|\theta-\bar \theta|d\x\\
        &\le C\|Q_1^I\|_{L^2(\bbr)} \|u-\bar u\|_{L^2(\bbr)}+ C \|Q_2^I\|_{L^2(\bbr)} \|\theta-\bar \theta\|_{L^2(\bbr)}\\
        &\quad+C\int_{\bbr} \big|u^R_{\xi\xi}\big| |(u-\bar u, \theta-\bar\theta)|d\x+C\int_{\bbr} \big|\big(v^R_{\xi},u^R_{\xi}\big)\big|^2 |(u-\bar u, \theta-\bar\theta)|d\x\\
        &\le C\|Q_1^I\|_{L^2(\bbr)} \|u-\bar u\|_{L^2(\bbr)}+ C \|Q_2^I\|_{L^2(\bbr)} \|\theta-\bar \theta\|_{L^2(\bbr)}\\   &\quad+C\|u^R_{\xi\xi}\|_{L^2(\bbr)}+C\|u^R_{\xi\xi}\|_{L^2(\bbr)}\|(u-\bar u, \theta-\bar\theta)\|_{L^2(\bbr)}^2\\           	&\quad+C\|\big(v^R_{\xi},u^R_{\xi}\big)\|_{L^{\frac{5}{2}}(\bbr)}^{\frac{5}{2}}+C\|\big(v^R_{\xi},u^R_{\xi}\big)\|_{L^{\infty}(\bbr)}^{\frac{3}{2}}\|(u-\bar u, \theta-\bar\theta)\|_{L^2(\bbr)}^2,
    \end{align*}
    which together with Lemma \ref{lem:wave-inter} implies
    \begin{align}\label{S12-e}
        \begin{aligned}
        \mb S_1+\mb S_2&\leq C \deltas\deltar \left(e^{-C \deltas t}+|(1+\deltar t)|^{-\frac{5}{2}}\right)\big(\|u-\bar u\|_{L^2(\bbr)}+\|\theta-\bar \theta\|_{L^2(\bbr)}\big)\\
        &\quad + C(1+t)^{-\frac{9}{8}}\|(u-\bar u, \theta-\bar\theta)\|_{L^2(\bbr)}^2+C\deltar^{\frac{1}{4}}(1+t)^{-\frac{9}{8}}.
        \end{aligned}
    \end{align}     			
    Then, from \eqref{sum-B}, \eqref{R}, \eqref{es-Y} and \eqref{S12-e}, together with the smallness of $\delta_0$ and $\eps_1$, we complete the proof of the lemma.
    \end{proof}
    Then, we proceed to prove Lemma \ref{Re-en} as follows. First, applying Lemmas \ref{leading-order} and \ref{lem-sharp2} and Young's inequality to \eqref{d-a-2}, we obtain
    \begin{align*}
        \frac{d}{d t} \int_{\mathbb{R}} a^{-\mathbf{X}} \bar{\theta} \eta(U | \bar{U}) d \xi 
        &\leq -\frac{\delta_{S}}{2 M}|\dot{\mathbf{X}}(t)|^{2}+\mathbf B_{1}+\mathbf B_{2}-\mathbf G-\frac{3}{4} \mathbf D \\ 
        & \quad -\frac{\delta_{S}}{2 M}|\dot{\mathbf{X}}(t)|^{2}+\dot{\mathbf{X}}(t) \sum_{i=4}^{6} \mathbf Y_{i}  +\sum_{i=3}^{6} \mathbf B_{i}+\mb S_1+\mb S_2-\frac{1}{8} \mathbf D -\frac{1}{8} \mathbf D\\
        &\leq -\frac{\delta_{S}}{2 M}|\dot{\mathbf{X}}(t)|^{2}+\mathbf B_{1}+\mathbf B_{2}-\mathbf G-\frac{3}{4} \mathbf D \\ 
        & \quad +\frac{C}{\delta_{S}} \sum_{i=4}^{6}\left|\mathbf Y_{i}\right|^{2}+\sum_{i=3}^{6} \mathbf B_{i}+\mb S_1+\mb S_2-\frac{1}{8} \mathbf D -\frac{\delta_{S}}{4 M}|\dot{\mathbf{X}}(t)|^{2}-\frac{1}{8} \mathbf D\\
        &\leq- \frac{C}{2} \mathcal{G}^{S}(U)+C \deltas^{4/3} \deltar^{4/3}  \left( e^{-C \deltas t}+|(1+\deltar t)|^{-10/3}\right) -C \mathcal{G}^R(U)\\
        &\quad + C \deltas\deltar \left( e^{-C \deltas t}+|(1+\deltar t)|^{-\frac{5}{2}}\right)\big(\|u-\bar u\|_{L^2(\bbr)}+\|\theta-\bar \theta\|_{L^2(\bbr)}\big)\\
        &\quad + C(1+t)^{-\frac{9}{8}}\|(u-\bar u, \theta-\bar\theta)\|_{L^2(\bbr)}^2+C\deltar^{\frac{1}{4}}(1+t)^{-\frac{9}{8}}\\
        &\quad  -\frac{\delta_{S}}{4 M}|\dot{\mathbf{X}}(t)|^{2}-\frac{1}{8} \mathbf D,
	\end{align*}
     together with the smallness of $\varepsilon_1$ and $\delta_0$.
            		  
     Integrating the inequality over $[0,t]$, for all $0\leq t\leq T$, together with Gronwall's inequality, we obtain
    \begin{align*}
         & \int_{\mathbb{R}} \eta(U | \bar{U})(t,\cdot) d \xi+\delta_{S} \int_{0}^{t}|\dot{\mathbf{X}}(\tau)|^{2} d \tau+\int_{0}^{t}\left(\mathcal{G}^{R}+\mathcal{G}^{S}+\mathbf D\right) d \tau \\
         &\leq C \int_{\mathbb{R}} \eta\left(U_{0}(\cdot) | \bar{U}(0, \cdot)\right) d \xi+C(\deltas^{1/3}\deltar ^{4/3}+\deltas^{4/3}\deltar^{1/3}+\deltar^{1/4})\\
        &\quad+C( \deltas+\deltar)\big(\|u-\bar u\|_{L^{\infty}(0,T;L^2(\bbr))}+\|\theta-\bar \theta\|_{L^{\infty}(0,T;L^2(\bbr))}\big),
    \end{align*}
    where $U_0=(v_0,u_0,\theta_0)^t$. Using the fact
    \begin{equation*}
       \| U-\bar{U}\| _{L^{2}(\mathbb{R})}^{2} \sim \int_{\mathbb{R}} \eta(U | \bar{U}) d \xi, \quad \forall t \in[0, T],
    \end{equation*}
    together with Young's inequality, we have
    \begin{align*}
        &\| U-\bar{U}(t,\cdot)\| _{L^{2}(\mathbb{R})}^{2}+\delta_{S} \int_{0}^{t}|\dot{\mathbf{X}}(\tau)|^{2} d \tau+\int_{0}^{t}\left(\mathcal{G}^{S}+\mathcal{G}^{R}+\mathbf D\right) d \tau \\
        & \leq C\left\| U_{0}(\cdot)-\bar{U}(0, \cdot)\right\| _{L^{2}(\mathbb{R})}^{2}+C\delta_0^{\frac{1}{4}},
    \end{align*}
    which implies
    \begin{align*}
        & \| (v - \bar{v}, u - \bar{u}, \theta- \bar{\theta})(t, \cdot)\| _{L^{2}(\mathbb{R})}^{2}+\delta_{S} \int_{0}^{t}|\dot{\mathbf{X}}(\tau)|^{2} d \tau+\int_{0}^{t}\left(\mathcal{G}^{S}+\mathcal{G}^{R}+\mathbf D\right) d \tau \\
        & \leq C\left\| (v_0(\cdot) - \bar{v}(0, \cdot), u_0(\cdot) - \bar{u}(0, \cdot), \theta_0(\cdot) - \bar{\theta}(0, \cdot))\right\| _{L^{2}(\mathbb{R})}^{2}+C\delta_0^{\frac{1}{4}}.
    \end{align*}            		        
    Then, using Lemma \ref{lem:variables-equivalence} and above inequality, we complete the proof of Lemma \ref{Re-en}.     		         		
    \section{Proof of Proposition \ref{prop:priori}}
            			
    In order to close the {\it a priori} estimates, we need to prove the following Lemma \ref{1-d} for $L^2$-estimates of the first-order derivatives and Lemma \ref{weight-s} for the weighted entropy derivative to complete the proof of Proposition \ref{prop:priori}. Throughout this section, the computations are carried out for smooth solutions; under the regularity assumed in Proposition \ref{prop:priori}, they are justified by a standard mollification argument, since all the resulting estimates involve only norms that are finite in that class.
    For simplicity, we first introduce the following notations
    \begin{align*}
        W(t, \xi) &:= p(t, \xi) - \bar{p}(t, \xi), \quad \psi(t, \xi) := u(t, \xi) - \bar{u}(t, \xi), \\
        \zeta(t, \xi) &:= s(t, \xi) - \bar{s}(t, \xi), \quad \phi(t,\xi) := v(t, \xi) - \bar{v}(t, \xi),\\
        \vartheta(t,\xi)&:= \theta(t, \xi) - \bar{\theta}(t, \xi).
    \end{align*}            		           				            
    Then, using the above definitions and \eqref{p-u-s-1}, together with equations \eqref{shock-p-u-s}, we obtain the following perturbation equations
        \begin{equation}\label{PUS}
            \begin{cases}
            W_t - \sigma W_{\xi} - \dot{\mathbf{X}}(t) (p^S)_{\xi}^{-\mathbf{X}}  +\frac{\gamma p}{v} \psi_{\xi}+ \left(\frac{\gamma p}{v}-\frac{\gamma \bar p}{\bar v} \right){\bar u}_{\xi}  = (\gamma-1)\mu \left( \frac{(u_{\xi})^2}{v^2} - \frac{(\bar{u}_{\xi})^2}{\bar{v}^2}  \right)-R_1,\\[1mm]
            \psi_t - \sigma \psi_{\xi} - \dot{\mathbf{X}}(t) (u^S)_{\xi}^{-\mathbf{X}} + W_{\xi} = \mu \left( \frac{u_{\xi}}{v} - \frac{\bar{u}_{\xi}}{\bar{v}} \right)_{\xi}-R_2, \\[1mm]
            \zeta_t -  \sigma\zeta_\xi -  \dot{\mathbf{X}}(t) (s^S)_{\xi}^{-\mathbf{X}} = \mu \left( \frac{(u_{\xi})^2}{v\theta} - \frac{(\bar{u}_{\xi})^2}{\bar{v}\bar{\theta}} \right)-R_3.
            \end{cases}
        \end{equation}

\begin{lemma}\label{lem:variables-weighted}
Assume that all the hypotheses of Proposition~\ref{prop:priori} hold. 
Let $\omega=\omega(t,\x)$ denote a generic function associated with the
wave profiles, such as $|(v^R)_{\x}|$ or $|(v^S)^{-\mb X}_{\x}|$, satisfying
\[
0\leq\omega(t,\x)\leq C.
\]
Then, we have
\begin{equation}
\label{eq:weighted-zero-forward}
\begin{aligned}
    \int_0^t\int_{\mathbb R}
    \omega |(v-\bar v,\theta-\bar\theta)|^2 d\x d\tau
    \leq
    C\int_0^t\int_{\mathbb R}
    \omega |(p-\bar p,s-\bar s)|^2 d\x d\tau
    +
    C\delta_0,
\end{aligned}
\end{equation}
and
\begin{equation}
\label{eq:weighted-zero-backward}
\begin{aligned}
       \int_0^t\int_{\mathbb R}
    \omega |(p-\bar p,s-\bar s)|^2 d\x d\tau
    \leq
    C\int_0^t\int_{\mathbb R}
    \omega |(v-\bar v,\theta-\bar\theta)|^2 d\x d\tau
    +
    C\delta_0.
\end{aligned}
\end{equation}

Furthermore, one has
\begin{equation}
\label{eq:weighted-first-forward}
\begin{aligned}
    &\int_0^t\int_{\mathbb R}
    \omega |(v-\bar v,\theta-\bar\theta)_\x|^2 d\x d\tau
    \leq
    C\int_0^t\int_{\mathbb R}
    \omega |(p-\bar p,s-\bar s)_\x|^2 d\x d\tau
    \\
    &\qquad\quad+
    C\int_0^t\int_{\mathbb R}
    \omega
    |(\bar v,\bar\theta)_\x|^2 |(v-\bar v,\theta-\bar \theta)|^2 d\x d\tau
    +
    C\delta_0.
\end{aligned}
\end{equation}
\end{lemma}
\begin{proof}
By \eqref{eq:raw-pointwise-forward}, we have
\begin{align}
\label{eq:weighted-pointwise-forward-proof}
    |(v-\bar v,\theta-\bar\theta)|^2
    \leq
    C|(p-\bar p,s-\bar s)|^2+
    C|(\mathcal Q_p,\mathcal Q_s)|^2.
\end{align}
Multiplying \eqref{eq:weighted-pointwise-forward-proof} by $\omega$
and integrating over $[0,t]\times\mathbb R$, we obtain
\begin{equation}
\begin{aligned}\label{eq:weighted-zero-forward-proof}
    \int_0^t\int_{\mathbb R}
    \omega
    |(v-\bar v,\theta-\bar\theta)|^2
    d\x d\tau
    \leq
    C\int_0^t\int_{\mathbb R}
    \omega
    |(p-\bar p,s-\bar s)|^2
    d\x d\tau+
    C\int_0^t\int_{\mathbb R}
    \omega
        |(\mathcal Q_p,\mathcal Q_s)|^2
    d\x d\tau.
\end{aligned}
\end{equation}
Since
\[
    0\leq\omega(t,\x)\leq C,
\]
the estimate \eqref{eq:Q-L2-proof} obtained in the proof of
Lemma~\ref{lem:variables-equivalence} gives
\begin{align}
\label{eq:weighted-Q-zero-proof}
    \int_0^t\int_{\mathbb R}
    \omega
        |(\mathcal Q_p,\mathcal Q_s)|^2
    d\x d\tau\leq
    C\int_0^t
    \|(\mathcal Q_p,\mathcal Q_s)(\tau,\cdot)\|_{L^2}^2 d\tau
    \leq
    C\delta_0.
\end{align}
Substituting \eqref{eq:weighted-Q-zero-proof} into
\eqref{eq:weighted-zero-forward-proof}, we obtain
\eqref{eq:weighted-zero-forward}. Similarly, we can obtain \eqref{eq:weighted-zero-backward}.

It remains to prove \eqref{eq:weighted-first-forward}. From
\eqref{eq:corrected-derivative-backward-proof}, we have
\begin{align}
\label{eq:weighted-first-pointwise-proof}
    |(v-\bar v,\theta-\bar\theta)_\x|^2
    \leq{}&
    C|(p-\bar p,s-\bar s)_\x|^2+
    C|(\bar v,\bar\theta)_\x|^2
      |(p-\bar p,s-\bar s)|^2
    \nonumber\\
    &+
    C|(\mathcal Q_p,\mathcal Q_s)_\x|^2+
    C|(\bar v,\bar\theta)_\x|^2
      |(\mathcal Q_p,\mathcal Q_s)|^2.
\end{align}
Multiplying \eqref{eq:weighted-first-pointwise-proof} by $\omega$ and
integrating over $[0,t]\times\mathbb R$, together with \eqref{eq:Q-L2-proof} and \eqref{eq:Qx-sup-proof}, we obtain
\begin{align*}
    &\int_0^t\int_{\mathbb R}
    \omega
    |(v-\bar v,\theta-\bar\theta)_\x|^2
    d\x d\tau
    \\
    &\leq
    C\int_0^t\int_{\mathbb R}
    \omega
    |(p-\bar p,s-\bar s)_\x|^2
    d\x d\tau+
    C\int_0^t\int_{\mathbb R}
    \omega
    |(\bar v,\bar\theta)_\x|^2
    |(p-\bar p,s-\bar s)|^2
    d\x d\tau
    +
    C\delta_0.
\end{align*}
Then, from Lemma \ref{lem:variables-equivalence}, we prove \eqref{eq:weighted-first-forward}.
\end{proof}

    \subsection{Estimates for the first-order derivatives}\label{subsec:first}
    We first carry out the estimates of the effective viscous flux in Lemma \ref{flux} and use it to establish the space-time estimates for the first-order derivative of pressure in Lemma \ref{pres}. Subsequently, we present the first-order estimate of the entropy in Lemma \ref{S}. On this basis, we eventually derive the estimates for all the first-order derivatives in Lemma \ref{1-d}.
    
    We introduce the effective viscous flux $F=\mu\frac{u_\x}{v}-p$, then define $\bar F=\mu\frac{\bar u_\x}{\bar v}-\bar p$ and $\mathcal{F}:=F-\bar F$. 

\medskip
\noindent\textbf{Roadmap of the first-order estimates.}
Since the argument closes through a loop, we first describe how the lemmas of
this section depend on each other. Set
\[
    \mathcal W(t):=\int_0^t\int_{\bbr}|(v^S)_\xi^{-\mathbf X}|\,|\zeta_\xi|^2\,d\xi d\tau .
\]
\begin{enumerate}[label=\textup{(\arabic*)},leftmargin=*]
\item \emph{Lemma \ref{flux} (velocity and effective viscous flux).}
Testing the velocity equation with $\psi_t$ controls $\|\psi_\xi(t)\|_{L^2}^2$
and $\int_0^t\|\mathcal F_\xi\|_{L^2}^2d\tau$ by a \emph{small} multiple
$(C\delta_0+\chi)$ of $\int_0^t\|W_\xi\|_{L^2}^2d\tau$, small multiples of
$\int_0^t\|\psi_{\xi\xi}\|_{L^2}^2d\tau$ and of $\mathcal W(t)$, and a term
$\int_0^t(1+\tau)^{-2}\|\zeta_\xi\|_{L^2}^2d\tau$ with a time-integrable
coefficient.
\item \emph{Lemma \ref{pres} (pressure).}
Written in terms of $\mathcal F$, the pressure equation \eqref{P} contains
the damping term $\frac{\gamma p}{\mu}W$. This gives
$\|W_\xi(t)\|_{L^2}^2+\int_0^t\|W_\xi\|_{L^2}^2d\tau
\le C\int_0^t\|\mathcal F_\xi\|_{L^2}^2d\tau+\cdots$.
Choosing $\chi$ and $\delta_0$ small, steps (1) and (2) close with each
other.
\item \emph{Lemma \ref{S} (entropy, pointwise in time).}
The differentiated entropy equation has neither dissipation nor damping; it
yields $\|\zeta_\xi(t)\|_{L^2}^2\le C\,\mathcal W(t)+\cdots$, where the
constant in front of $\mathcal W(t)$ is \emph{not} small.
\item \emph{Lemma \ref{1-d}.}
Combining (1)-(3), recovering
$\int_0^t\|\psi_{\xi\xi}\|_{L^2}^2d\tau$ from the momentum equation and
applying Gronwall's inequality, all first-order quantities are controlled by
the initial data, $\delta_0^{1/4}$ and $C\,\mathcal W(t)$.
\item \emph{Lemmas \ref{trans}-\ref{trans-1} (transport lemmas).}
These three lemmas convert source terms that are weighted by the shock
profile, weighted by the rarefaction wave, or merely integrable in time with
values in $L^2$, into bounds for shock-weighted integrals along the entropy
characteristics. Each of them rests on the mechanism \eqref{key-mech}, and
each gains a factor $\delta_S$.
\item \emph{Lemma \ref{weight-s} (weighted entropy estimate).}
Representing $\zeta_\xi$ along the characteristics and applying step (5),
$\mathcal W(t)$ is bounded by the initial data, $\delta_0^{1/4}$, \emph{small}
multiples of $\int_0^t\|W_\xi\|_{L^2}^2d\tau$ and
$\int_0^t\|\psi_{\xi\xi}\|_{L^2}^2d\tau$, and a rarefaction-weighted term
with a time-integrable coefficient. The quadratic term
$\delta_S\bigl(\int_0^t\|\psi_\xi\|_{H^1}^2d\tau\bigr)^2$ is where the
$L^2(0,T;H^2)$-part of the {\it a priori} assumption \eqref{apri-ass} is used.
\item \emph{Lemma \ref{one-order} (conclusion).}
Inserting (6) into (4), the small multiples are absorbed into the left-hand
side and Gronwall's inequality closes the estimate.
\end{enumerate}
The smallness that closes the loop between (4) and (6) comes from the factor
$\delta_S$ gained in step (5) and from $\delta_0$ and $\varepsilon_1$.
\medskip
    \begin{lemma}\label{flux}
        Under the hypotheses of Proposition \ref{prop:priori}, there exists $C > 0$ such that for all $t \in (0, T]$,
         \begin{align*}
          &\|\psi_\x(t,\cdot)\|_{L^2(\mathbb{R})}^2 +\int_0^t   \|\mathcal{F}_\x\|^2_{L^2(\bbr)} d\tau\\
          &\le C\|(W,\psi,\zeta)(0,\cdot)\|_{H^1(\mathbb{R})}^2+(C\delta_0+\chi)\int_0^t\int_{\bbr} |W_\x|^2 d\x d\tau+C\int_0^t (1+\tau)^{-2}\int_{\bbr} |\zeta_\x|^2d\x d\tau\\
          &\quad +C\delta_0\int_0^t\int_{\bbr} |(v^S)_\xi^{-\mathbf{X}}| |\zeta_\x|^2d\x d\tau+C(\varepsilon_1+\delta_0) \int_0^t\int_{\bbr}  |\psi_{\x\x}|^2d\x d\tau+C\delta_0^{\frac 14},
        \end{align*}
        where $\chi$ is a sufficiently small positive constant independent of $\delta_0, \deltas, \varepsilon_1$ and $T$.
    \end{lemma}
        \begin{proof}
        Multiplying $\eqref{PUS}_2$ by $\psi_t$, we get
            \begin{equation*}
            	\psi^2_t - \sigma \psi_{\xi} \psi_t - \dot{\mathbf{X}}(t) (u^S)_{\xi}^{-\mathbf{X}}\psi_t + W_{\xi} \psi_t= \mu \left( \frac{\psi_{\xi}}{v}\right)_\x \psi_t +\mu\left[\bar u_\x\left(\frac{1}{v}- \frac{1}{\bar{v}} \right)\right]_{\xi}\psi_t-R_2\psi_t,
            \end{equation*}
        which implies
        \begin{align*}
            &\mu\left(\frac{1}{2v}\psi_\x^2\right)_t+\psi^2_t - (W\psi_\x)_t-\mu \left( \frac{\psi_{\xi}}{v}\psi_t\right)_\x -\mu \left( \frac{1}{2v}\right)_t\psi^2_\x+(W\psi_t)_{\x} \\
            &- \sigma \psi_{\xi} \psi_t - \dot{\mathbf{X}}(t) (u^S)_{\xi}^{-\mathbf{X}}\psi_t + W_{t} \psi_\x-\mu\left[\bar u_\x\left(\frac{1}{v}- \frac{1}{\bar{v}} \right)\right]_{\xi}\psi_t+R_2\psi_t= 0.
        \end{align*}
    Integrating the equation, we obtain
        \begin{equation*}
            \|\psi_\x(t,\cdot)\|_{L^2(\mathbb{R})}^2 +\int_0^t\int_{\bbr}\psi^2_td\x d\tau \le C\|\psi_\x(0,\cdot)\|_{L^2(\mathbb{R})}^2+ \sum_{i=1}^{5} L_i,
        \end{equation*}
    where 
        \begin{align*}            				                   &L_1:=C\left(\int_{\bbr}|W\psi_\x|d\x \right)\Big|_{\tau=0}+C\left(\int_{\bbr}|W\psi_\x|d\x \right)\Big|_{\tau=t},\\
    	   &L_2:=C\int_0^t\int_{\bbr}\left|\left( \frac{1}{v}\right)_{\tau}\psi^2_\x\right|d\x d\tau,\\
          &L_3:=C\int_0^t\int_{\bbr}|\psi_t| \left(|\psi_{\xi} |+ |\dot{\mathbf{X}}(\tau) (u^S)_{\xi}^{-\mathbf{X}}| +\left|\left(\bar u_\x\left(\frac{1}{v}- \frac{1}{\bar{v}} \right)\right)_{\xi}\right|\right)d\x d\tau, \\
          &L_4:=C\int_0^t\int_{\bbr}|\psi_\x| |W_t|d\x d\tau,\\
          &L_5:=C\int_0^t\int_{\bbr}|\psi_t| |R_2|d\x d\tau.
        \end{align*}
    For $L_1$, by direct computation and Young's inequality, we can get
    \begin{equation}\label{L1}
       L_1\leq C(\|W(t,\cdot)\|_{L^2(\mathbb{R})}^2+\|W(0,\cdot)\|_{L^2(\mathbb{R})}^2+\|\psi_\x(0,\cdot)\|_{L^2(\mathbb{R})}^2)+\frac{1}{2}\|\psi_\x(t,\cdot)\|_{L^2(\mathbb{R})}^2.
    \end{equation}
    For $L_2$, using the fact 
        \begin{align*}
         |v_t| &\leq C(|v_\x|+|u_\x|)\\
         &\leq C (|(\phi_\x,\psi_\x)|+|(\bar v_\x,\bar u_\x)|),
        \end{align*}
     and Lemma \ref{lem:variables-boundedness}, together with Sobolev's inequality, Lemmas \ref{lem:shock-v} and \ref{lem:rare-v}, we have 
    \begin{equation}\label{L2}
       \begin{aligned}
          L_2&\leq C\int_0^t\int_{\bbr}|(\phi_\x,\psi_\x)| \psi^2_\x d\x d\tau  +C(\deltas+\deltar)\int_0^t\int_{\bbr}\psi^2_\x d\x d\tau\\
        &\leq C\int_0^t
        \|\psi_\x\|_{L^\infty(\mathbb{R})}
        \|(\phi_\x,\psi_\x)\|_{L^2(\mathbb{R})} \|\psi_\x \|_{L^2(\mathbb{R})} d\tau
        +C\delta_0\int_0^t\int_{\bbr}\psi^2_\x d\x d\tau\\
        &\leq C(\varepsilon_1+\delta_0)\int_0^t\int_{\bbr}\psi^2_{\x\x} d\x d\tau +C(\varepsilon_1+\delta_0)\int_0^t\int_{\bbr}\psi^2_\x d\x d\tau.
        \end{aligned}
    \end{equation}
     For $L_3$, from Lemmas \ref{lem:shock-v}, \ref{lem:variables-weighted} and \ref{lem:rare-v}, using Young's inequality, we obtain
    \begin{equation}\label{L3}
    	\begin{aligned}
            L_3&\leq C\int_0^t\int_{\bbr}|\psi_t| (|\psi_{\xi} |+ |\dot{\mathbf{X}}(\tau) (u^S)_{\xi}^{-\mathbf{X}}| +|\bar{v}_\x||(\phi_\x,\phi)|)d\x d\tau\\
            &\leq \frac{1}{4}\int_0^t\int_{\bbr}\psi^2_td\x d\tau
            +C\deltas\int_0^t |\dot{\mathbf{X}}(\tau)|^2  d\tau+ C\int_0^t\int_{\bbr}\psi^2_\x d\x d\tau\\
            &\quad +C\int_0^t\int_{\bbr} |\bar{v}_\x|^2 |(\phi_\x,\phi)|^2d\x d\tau\\
            &\leq \frac{1}{4}\int_0^t\int_{\bbr}\psi^2_td\x d\tau
            +C\deltas\int_0^t |\dot{\mathbf{X}}(\tau)|^2  d\tau+C\int_0^t\int_{\bbr}\psi^2_\x d\x d\tau+C\delta_0 \int_0^t(\mathcal{G}^S+\mathcal{G}^R)d\tau\\
            &\quad +C\int_0^t\int_{\bbr} |v^R_\x|^2 |\phi_\x|^2d\x d\tau+ C\delta_0\int_0^t\int_{\bbr} |(v^S)^{-\mb X}_\x| |\phi_\x|^2d\x d\tau\\
            &\leq \frac{1}{4}\int_0^t\int_{\bbr}\psi^2_td\x d\tau
            +C\deltas\int_0^t |\dot{\mathbf{X}}(\tau)|^2  d\tau+C\int_0^t\int_{\bbr}\psi^2_\x d\x d\tau+C\delta_0 \int_0^t(\mathcal{G}^S+\mathcal{G}^R)d\tau+C\delta_0\\
            &\quad +C\int_0^t(1+\tau)^{-2}\int_{\bbr}  |\zeta_\x|^2d\x d\tau+ C\delta_0\int_0^t\int_{\bbr} |(v^S)^{-\mb X}_\x| |\zeta_\x|^2d\x d\tau+ C\delta_0\int_0^t\int_{\bbr} |W_\x|^2d\x d\tau.
        \end{aligned}
    \end{equation}
     For $L_4$, we use $\eqref{PUS}_1$, Young's inequality and Sobolev's inequality, together with Lemmas \ref{lem:shock-p}, \ref{lem:wave-inter} and \ref{lem:variables-equivalence} to obtain
    \begin{equation}\label{L4}
        \begin{aligned}
        L_4&\le C \int_0^t\int_{\bbr}  |\psi_\x| (|W_\x| + |\dot{\mathbf{X}} (p^S)_{\xi}^{-\mathbf{X}}|  + |\psi_{\xi}| +|\bar{v}_\x| |(W,\phi)| + 
        |\psi_\x|^2 +|R_1| ) d\x d\tau\\
        &\le C \int_0^t\int_{\bbr}  |\psi_\x|^2d\x d\tau
        +C\deltas\int_0^t |\dot{\mathbf{X}}|^2  d\tau
        +C\int_0^t\int_{\bbr} |\bar{v}_\x|^2 |(W,\phi)|^2d\x d\tau\\
        &\quad+\chi  \int_0^t\int_{\bbr}  |W_\x|^2d\x d\tau
        +C\int_0^t \| \psi_\x\|_{L^2(\mathbb{R})}^2\| \psi_\x\|_{L^\infty(\mathbb{R})} d\tau+C\int_0^t\int_{\bbr}  |R_1|^2d\x d\tau\\\
        &\le C \int_0^t\int_{\bbr}  |\psi_\x|^2d\x d\tau
        +C\deltas\int_0^t |\dot{\mathbf{X}}(\tau)|^2  d\tau
        +C\deltas\int_0^t\int_{\bbr} |(v^S)_\xi^{-\mathbf{X}}| |(W,\phi)|^2d\x d\tau\\
        &\quad +C\int_0^t\int_{\bbr} |v^R_\xi|^2 |(W,\phi)|^2d\x d\tau+\chi  \int_0^t\int_{\bbr}  |W_\x|^2d\x d\tau
        +C\varepsilon_1 \int_0^t\int_{\bbr}  |\psi_{\x\x}|^2d\x d\tau+C\delta_0,
        \end{aligned}
    \end{equation}
    where $\chi$ is a small positive constant.\\
    For $L_5$, using Lemmas \ref{lem:rare-v}, \ref{lem:wave-inter} and Young's inequality, we obtain
    \begin{equation}\label{L5}
    	\begin{aligned}
            L_5 &\leq \frac{1}{4}\int_0^t\int_{\bbr}\psi^2_t d\x d\tau+C \int_0^t\int_{\bbr}|R_2|^2 d\x d\tau\\
            &\leq \frac{1}{4}\int_0^t\int_{\bbr} \psi^2_td\x d\tau+C\delta_0.
        \end{aligned}
    \end{equation}
     Then, by Lemmas \ref{Re-en}, \ref{lem:variables-weighted} and \ref{lem:rare-v}, together with \eqref{L1}-\eqref{L5}, we have
      \begin{align*}
        &\|\psi_\x(t,\cdot)\|_{L^2(\mathbb{R})}^2 +\int_0^t\int_{\bbr}\psi^2_td\x d\tau \\
        &\le C\|(W,\psi,\zeta)(0,\cdot)\|_{H^1(\mathbb{R})}^2+(C\delta_0+\chi)\int_0^t\int_{\bbr} |W_\x|^2 d\x d\tau+C\int_0^t (1+\tau)^{-2} \int_{\bbr}  |\zeta_\x|^2d\x d\tau\\
    	&\quad +C\delta_0\int_0^t\int_{\bbr} |(v^S)_\xi^{-\mathbf{X}}| |\zeta_\x|^2d\x d\tau+C(\varepsilon_1+\delta_0) \int_0^t\int_{\bbr}  |\psi_{\x\x}|^2d\x d\tau+C\delta_0^{\frac 14}.
     \end{align*}
    By the definition of $\mathcal{F}$ and $\eqref{PUS}_2$, we can get the equation for $\mathcal{F}$
    \begin{equation*}
        \psi_t - \sigma \psi_{\xi} - \dot{\mathbf{X}}(t) (u^S)_{\xi}^{-\mathbf{X}}  = \mathcal{F}_\x-R_2,
    \end{equation*}
 which implies 
  \begin{align*}
    \int_0^t   \|\mathcal{F}_\x\|^2_{L^2(\bbr)} d\tau
    \le&
	C\int_0^t   \|\psi_t\|^2_{L^2(\bbr)} d\tau+C \int_0^t   \|\psi_\x\|^2_{L^2(\bbr)} d\tau \\
    &+C\deltas\int_0^t |\dot{\mathbf{X}}(\tau)|^2  d\tau+	C\int_0^t   \|R_2\|^2_{L^2(\bbr)} d\tau.
    \end{align*}
   Thus, we can complete the proof of Lemma \ref{flux} together with Lemmas \ref{lem:wave-inter} and \ref{Re-en}.
   \end{proof}
 \begin{remark}
    In the proof of Lemma \ref{flux}, we also obtain the space-time estimate for $\psi_t$ as follows
 \begin{equation}\label{psi-t}
    \begin{aligned}
        \int_0^t\int_{\bbr}\psi^2_td\x d\tau  
        &\le C\|(W,\psi,\zeta)(0,\cdot)\|_{H^1(\mathbb{R})}^2+(C\delta_0+\chi)\int_0^t\int_{\bbr} |W_\x|^2 d\x d\tau\\
        &\quad +C\int_0^t (1+\tau)^{-2} \int_{\bbr}  |\zeta_\x|^2d\x d\tau+C\delta_0\int_0^t\int_{\bbr} |(v^S)_\xi^{-\mathbf{X}}| |\zeta_\x|^2d\x d\tau\\
        &\quad +C(\varepsilon_1+\delta_0) \int_0^t\int_{\bbr}  |\psi_{\x\x}|^2d\x d\tau+C\delta_0^{\frac 14},
    \end{aligned}
    \end{equation}
    which we will use in the proof of Lemma \ref{1-d}.
  \end{remark}
            	
  \begin{lemma}\label{pres}
   Under the hypotheses of Proposition \ref{prop:priori}, there exists $C > 0$ such that for all $t \in (0, T]$,
           	 \begin{align*}
           			\|W_\x(t,\cdot)\|_{L^2(\mathbb{R})}^2 +\int_0^t   \|W_\x\|^2_{L^2(\bbr)} d\tau
           			&\le C\|(W,\psi,\zeta)(0,\cdot)\|_{H^1(\mathbb{R})}^2
           			+C\int_0^t\int_{\bbr}|\mathcal{F}_\x|^2d\x d\tau\\
           			&\quad +C\varepsilon_1 \int_0^t\int_{\bbr}  |\psi_{\x\x}|^2d\x d\tau+C\delta_0^{\frac 14}.
           	\end{align*}
           \end{lemma}
            				
            \begin{proof}
            	To employ the effective viscous flux $\mathcal{F}$, we rewrite the equation $\eqref{PUS}_1$ as
            	\begin{equation}\label{P}
            		\begin{aligned}
            		&W_t - \sigma W_{\xi} - \dot{\mathbf{X}}(t) (p^S)_{\xi}^{-\mathbf{X}}  +\frac{\gamma p}{\mu} W+ \frac{\gamma p}{\mu} \mathcal{F}+\gamma\frac{\bar u_\x}{\bar v}W\\            		
            		&= (\gamma-1)(\mathcal{F}+W) \left( \frac{\mathcal{F}+W}{\mu} +2\frac{\bar{u}_{\xi}}{\bar{v}}  \right)-R_1.
            		\end{aligned}
            	\end{equation}
            	Differentiating both sides of \eqref{P} with respect to $\x$
            	and multiplying by $W_\x$, we obtain
            	\begin{equation}\label{P-d}
            		\begin{aligned}
            			&\quad \frac{1}{2}\frac{d}{dt} |W_\x|^2
            			-\frac{1}{2}\sigma \left(W^2_\x\right)_\x
            			+\frac{\gamma p}{\mu} W^2_\x
            			+\frac{\gamma}{\mu} W W^2_\x
            			+\frac{\gamma}{\mu} \bar p_\x W W_\x\\
            			&=\dot{\mathbf{X}}(t) (p^S)_{\xi\xi}^{-\mathbf{X}} W_\x
            			-\frac{\gamma p}{\mu} \mathcal{F_\x} W_\xi
            			-\frac{\gamma}{\mu} \bar p_\x \mathcal{F} W_\x
            			-\gamma\frac{\bar u_\x}{\bar v} W^2_\x 
            			-\gamma\left(\frac{\bar u_\x}{\bar v}\right)_\x W W_\x\\
            			&\quad -\frac{\gamma}{\mu} \mathcal{F} W^2_\xi+ (\gamma-1)\left[(\mathcal{F}+W) \left( \frac{\mathcal{F}+W}{\mu} +2\frac{\bar{u}_{\xi}}{\bar{v}}  \right)\right]_\x W_\x-W_\x R_{1,\x}.
            		\end{aligned}
            	\end{equation}
            	Integrating the equation \eqref{P-d}, we get 
            	\begin{equation*}
            		\|W_\x(t,\cdot)\|_{L^2(\mathbb{R})}^2 +\int_0^t   \|W_\x\|^2_{L^2(\bbr)} d\tau
            		\le C\|W_\x(0,\cdot)\|_{L^2(\mathbb{R})}^2 +
            		\sum_{i=1}^{4} K_i,          			
            	\end{equation*}
            	where 
             \begin{align*}
            			&K_1:=C\int_0^t\int_{\bbr} |W_\x| 
            			\left( |W W_\x|
            			+| \bar p_\x W| 
            			+|\dot{\mathbf{X}}(\tau) (p^S)_{\xi\xi}^{-\mathbf{X}}|
            			+|\mathcal{F_\x}| \right.\\
                        &\qquad\qquad\qquad\qquad\qquad\left.
            			+|\bar p_\x \mathcal{F}| 
            			+|\bar u_\x W_\x| 
            			+\left|\left(\frac{\bar u_\x}{\bar v}\right)_\x W\right|  
            			\right) d\x d\tau,\\
            			&K_2:=C\int_0^t\int_{\bbr} 
            			|\mathcal{F} W^2_\x|            		
            			d\x d\tau,\\
            			&K_3:=C\int_0^t\int_{\bbr} 
            			\left|\left[(\mathcal{F}+W) \left( \frac{\mathcal{F}+W}{\mu} +2\frac{\bar{u}_{\xi}}{\bar{v}}  \right)\right]_\x\right| |W_\x|
            			d\x d\tau\\
            			&K_4:=C\int_0^t\int_{\bbr} 
            			|W_\x| |R_{1,\x}|            		
            			d\x d\tau.
            	\end{align*}
            	For $K_1$, using the fact 
            		\begin{align*}
            			|\mathcal{F}|\le C(|\psi_\x|+|\bar u_\x \phi|+|W|),
            		\end{align*}
            	 and applying Young's inequality and Sobolev's inequality, together with Lemmas \ref{lem:shock-v}, \ref{lem:shock-p} and \ref{lem:variables-weighted}, we have
            	 \begin{equation}\label{K1}
            		\begin{aligned}
            			K_1
            			&\le \frac{1}{8} \int_0^t\int_{\bbr}  |W_\x|^2d\x d\tau
            			+C\deltas\int_0^t |\dot{\mathbf{X}}(\tau)|^2  d\tau
            			+C\deltas\int_0^t\int_{\bbr} |(v^S)_\xi^{-\mathbf{X}}| |(W,\phi)|^2d\x d\tau\\
            			&\quad +C\int_0^t\int_{\bbr} |(p^R_\xi,u^R_\x)|^2 |(W,\phi)|^2d\x d\tau+C \int_0^t\int_{\bbr}  (|\mathcal{F}_\x|^2+|\psi_\x|^2) d\x d\tau\\
            			&\le \frac{1}{8} \int_0^t\int_{\bbr}  |W_\x|^2d\x d\tau
            			+C\deltas\int_0^t |\dot{\mathbf{X}}(\tau)|^2  d\tau
            			+C\delta_0 \int_0^t(\mathcal{G}^R+\mathcal{G}^S)d\tau\\
            			&\quad+ C \int_0^t\int_{\bbr}  (|\mathcal{F}_\x|^2+|\psi_\x|^2) d\x d\tau+C\delta_0,
            		\end{aligned}
            	\end{equation}
            	where we use the smallness of $\varepsilon_1$ and $\delta_0$.\\
            	For $K_2$,  using the fact 
            	\begin{equation*}
            			|\mathcal{F}|\le C(|\psi_\x|+|\bar u_\x \phi|+|W|),
            	\end{equation*}
            	and Young's inequality, together with Sobolev's inequality and the {\it a priori} assumption \eqref{apri-ass}, we have
            	\begin{equation}\label{K2}
            		\begin{aligned}
            			K_{2}
            			&\le   C(\varepsilon_1+\delta_0) \int_0^t\int_{\bbr}  |W_\x|^2d\x d\tau 	
            			 +C \int_0^t\int_{\bbr}  |\psi_\x||W_\x|^2d\x d\tau\\
            			&\le   C(\varepsilon_1+\delta_0)  \int_0^t\int_{\bbr}  |W_\x|^2d\x d\tau 	
            			 +C \int_0^t  \|\psi_\x\|_{L^\infty(\bbr)} \|W_\x\|^2_{L^2(\bbr)} d\tau\\
            			 &\le   C(\varepsilon_1+\delta_0)  \int_0^t\int_{\bbr}  |W_\x|^2d\x d\tau 	
            			 +C\sup_{0\le \tau\le t} \|W_\x(\tau, \cdot)\|_{L^2(\bbr)}  \int_0^t  \|\psi_\x\|_{H^1(\bbr)} \|W_\x\|_{L^2(\bbr)} d\tau\\
            			 &\le   C(\varepsilon_1+\delta_0)  \int_0^t\int_{\bbr}  |W_\x|^2d\x d\tau 	
            			 +C\varepsilon_1  \int_0^t  \|\psi_\x\|_{H^1(\bbr)}^2 d\tau.
            		\end{aligned}
            	\end{equation}
            	
            	As for  $K_3$, it can be estimated similarly as for 
            	$K_1$ and $K_{2}$, together with the smallness of $\varepsilon_1$ and $\delta_0$,  so we obtain
            \begin{equation}\label{K3}
            		\begin{aligned}
            			K_3
            			&\le \frac{1}{8} \int_0^t\int_{\bbr}  |W_\x|^2d\x d\tau
            			+C\delta_0 \int_0^t(\mathcal{G}^R+\mathcal{G}^S)d\tau+C \int_0^t\int_{\bbr}  (|\mathcal{F}_\x|^2+|\psi_\x|^2) d\x d\tau\\
            			&\quad 
            			+C\varepsilon_1 \int_0^t\int_{\bbr}  |\psi_{\x\x}|^2 d\x d\tau+C\delta_0.
            		\end{aligned}
            	\end{equation}
            	For $K_4$, using Young's inequality, Lemmas \ref{lem:rare-v} and \ref{lem:wave-inter}, we get
            	 \begin{equation}\label{K4}
            		\begin{aligned}
            			K_4
            			&\le \frac{1}{8} \int_0^t\int_{\bbr}  |W_\x|^2d\x d\tau
            			+C \int_0^t\int_{\bbr}  |R_{1,\x}|^2d\x d\tau\\
            			&\le \frac{1}{8} \int_0^t\int_{\bbr}  |W_\x|^2d\x d\tau
            			+C\delta_0.
            		\end{aligned}
            	\end{equation}
            	Combining Lemma \ref{Re-en} with \eqref{K1}-\eqref{K4}, we have 
            	 \begin{align*}
            			\|W_\x(t,\cdot)\|_{L^2(\mathbb{R})}^2 +\int_0^t   \|W_\x\|^2_{L^2(\bbr)} d\tau
            			&\le C\|(W,\psi,\zeta)(0,\cdot)\|_{H^1(\mathbb{R})}^2
            			+C\int_0^t\int_{\bbr}|\mathcal{F}_\x|^2d\x d\tau\\
            			&\quad +C\varepsilon_1 \int_0^t\int_{\bbr}  |\psi_{\x\x}|^2d\x d\tau+C\delta_0^{\frac 14}. 
            	\end{align*}
            	Thus, we complete the proof of the lemma.
            \end{proof}
            
            \begin{lemma}\label{S}
            	Under the hypotheses of Proposition \ref{prop:priori}, there exists $C > 0$ such that for all $t \in (0, T]$,
            	 \begin{align*}
            			\|\zeta_\x(t,\cdot)\|_{L^2(\mathbb{R})}^2 
            			&\le 
            			C\|(W,\psi,\zeta)(0,\cdot)\|_{H^1(\mathbb{R})}^2+
            			C\int_0^t\int_{\bbr}|(v^S)_{\xi}^{-\mathbf{X}}| |\zeta_\x|^2  
            			d\x d\tau  \\
            			&\quad+C\delta_0\int_0^t\int_{\bbr} |W_\x|^2  
            			d\x d\tau  
            			+C(\varepsilon_1+\delta_0)
            			\int_0^t \|\psi_{\x\x}\|^2_{L^2(\bbr)} d\tau+C\delta_0^{\frac 14}.
            		\end{align*}
            \end{lemma}
            \begin{proof}
            	 For simplicity, we introduce the following notation
            	\begin{equation}\label{H}
            		\begin{aligned}
            			H(t, \xi):= \mu \left( \frac{(u_{\xi})^2}{v\theta} - \frac{(\bar{u}_{\xi})^2}{\bar{v}\bar{\theta}} \right).
            		\end{aligned}            			
            	\end{equation}
            	Differentiating both sides of $\eqref{PUS}_3$ with respect to $\x$
            	and multiplying by $\zeta_\x$, we obtain
            	\begin{equation}\label{s-d}
            		\frac{1}{2}\frac{d}{dt}|\zeta_\x|^2-\frac{1}{2}\sigma\left(\zeta^2_\x\right)_\x- \dot{\mathbf{X}}(t) (s^S)_{\xi\x}^{-\mathbf{X}} \zeta_\x=H_\x \zeta_\x-R_{3,\x}\zeta_\x.
            	\end{equation}
            	Integrating the equation \eqref{s-d}, we get 
            \begin{align*}
            			\|\zeta_\x(t,\cdot)\|_{L^2(\mathbb{R})}^2 \le \|\zeta_\x(0,\cdot)\|_{L^2(\mathbb{R})}^2+ \sum_{i=1}^3 N_i,
            \end{align*}
            	where 
    \begin{align*}
        &N_1:=C\int_0^t\int_{\bbr} |\zeta_\x|                 	
        |\dot{\mathbf{X}}(\tau) (s^S)_{\xi\xi}^{-\mathbf{X}}|
        d\x d\tau,\\
        &N_2:=C\int_0^t\int_{\bbr} 
        |H_\x||\zeta_\x| d\x d\tau,\\
        &N_3:=C\int_0^t\int_{\bbr} 
        |R_{3,\x}||\zeta_\x| d\x d\tau.
    \end{align*}
    For $N_1$, using Lemma \ref{lem:shock-p} and Young's inequality, we have 
               \begin{equation}\label{N1}
               	\begin{aligned}
               		N_1
               		&\le 
               		C\deltas\int_0^t\int_{\bbr} |\zeta_\x|                 	
               		|\dot{\mathbf{X}}(\tau)| |(s^S)_{\xi}^{-\mathbf{X}}|
               		d\x d\tau\\
               		&\le
               		C\deltas\int_0^t\int_{\bbr}|(v^S)_{\xi}^{-\mathbf{X}}| |\zeta_\x|^2  
               		d\x d\tau
               		+C\deltas\int_0^t |\dot{\mathbf{X}}(\tau)|^2  d\tau.
               	\end{aligned}
               \end{equation}
               For $N_2$,  we first write all terms of $H_{\x}$  as follows                 
               \begin{align*}
               		H_\xi =\ & \frac{2\mu}{v\theta} \bar{u}_\xi \psi_{\xi\xi}+ \frac{2\mu}{v\theta} \bar{u}_{\xi\xi} \psi_\xi - 2\mu \left( \frac{\phi \bar{\theta} + \bar{v} \vartheta + \phi \vartheta}{v\theta \bar{v}\bar{\theta}} \right) \bar{u}_\xi \bar{u}_{\xi\xi}\\
               		&+ \mu \left( \frac{\bar{v}_\xi \bar{\theta} + \bar{v} \bar{\theta}_\xi}{\bar{v}^2 \bar{\theta}^2} - \frac{(\bar{v}_\xi + \phi_\xi)\theta + v(\bar{\theta}_\xi + \vartheta_\xi)}{v^2 \theta^2} \right) \bar{u}_\xi^2 \\
               		&- 2\mu \frac{\bar{v}_\xi \theta + v\bar{\theta}_\xi }{v^2 \theta^2} \bar{u}_\xi \psi_\xi  +\left[ \frac{2\mu}{v\theta} \psi_\xi \psi_{\xi\xi} - \mu \frac{ \phi_\xi\theta + v\vartheta_\xi}{v^2 \theta^2} \psi_\xi^2\right]\\
               		&- \left[2\mu \frac{\phi_\xi\theta + v\vartheta_\xi}{v^2 \theta^2} \bar{u}_\xi \psi_\xi + \mu \frac{\bar{v}_\xi\theta + v\bar{\theta}_\xi}{v^2 \theta^2} \psi_\xi^2\right],
               	\end{align*}
              and then using Lemmas \ref{lem:rare-v} and \ref{lem:shock-v}, we have 
               	\begin{align*}
               		|H_\x|&\leq 
               		C(|\bar{v}_\xi|  |\psi_{\xi\xi}| + |\bar{v}_\xi| |\psi_{\xi}| +|\bar{v}_\xi|^2 \big( |\phi| + |\vartheta| \big)+|\bar{v}_\xi|^2 ( |\phi_\xi| + |\vartheta_\xi|) )\\
               		&+C(|\psi_\xi|  |\psi_{\xi\xi}| + |(\phi_\xi,\vartheta_\xi)|  |\psi_{\xi}|^2)
               		+C(|\bar{v}_\xi|  |\psi_\xi|  |(\phi_{\xi},\vartheta_{\x})| + |\bar{v}_\xi|  |\psi_{\xi}|^2).
               	\end{align*}
               Therefore, we decompose $N_2$ into three parts, defined as follows:
               	\begin{align*}
               		&N_{21}:=C\int_0^t\int_{\bbr} |\zeta_\x| |\bar{v}_\x|         	
               		( |\psi_{\xi\xi}| +   |\psi_{\xi}| +|\bar{v}_\x|\big( |\phi| + |\vartheta| \big)+|\bar{v}_\x|( |\phi_\xi| + |\vartheta_\xi|) )
               		d\x d\tau,\\
               		&N_{22}:=C\int_0^t\int_{\bbr} 
               	|\zeta_\x| (|\psi_\xi|  |\psi_{\xi\xi}| + |(\phi_\xi,\vartheta_\xi)|  |\psi_{\xi}|^2)d\x d\tau,\\
               	&N_{23}:=C\int_0^t\int_{\bbr} 
               	|\zeta_\x| (|\bar{v}_\x|
               	 |\psi_\xi|  |(\phi_{\xi},\vartheta_{\x})| + |\bar{v}_\x|  |\psi_{\xi}|^2)d\x d\tau.           	
               	\end{align*}
               For $N_{21}$, from Lemmas \ref{lem:shock-v}, \ref{lem:rare-v} and \ref{lem:variables-weighted}, using Young's inequality, we have 
               \begin{equation}\label{N21}
               	\begin{aligned}
               		N_{21}
               		&\leq 
               		C\int_0^t\int_{\bbr}|(v^S)_{\xi}^{-\mathbf{X}}| |(\phi,\vartheta)|^2  
               		d\x d\tau          	
               		+C\int_0^t\int_{\bbr}|(v^S)_{\xi}^{-\mathbf{X}}| |(\phi,\vartheta,\zeta)_\x|^2  
               		d\x d\tau\\
               		&\quad + 
               		C\int_0^t\int_{\bbr}|v^R_{\xi}| |(\phi,\vartheta)|^2  
               		d\x d\tau          	
               		+C\int_0^t\int_{\bbr}|v^R_{\xi}|^{\frac 32} |(\phi,\vartheta,\zeta)_\x|^2  
               		d\x d\tau\\
               		&\quad +C\delta_0\int_0^t\int_{\bbr} |(\psi_\x,\psi_{\x\x})|^2  
               		d\x d\tau \\
               		&\leq C\int_0^t(\mathcal{G}^S+\mathcal{G}^R)d\tau      	
               		+C\int_0^t\int_{\bbr}|(v^S)_{\xi}^{-\mathbf{X}}| |\zeta_\x|^2  
               		d\x d\tau
               		+C\delta_0\int_0^t\int_{\bbr}|W_\x|^2  
               		d\x d\tau\\
               		&\quad +C\int_0^t(1+\tau)^{-\frac 32} \|\zeta_\x\|^2 _{L^2(\bbr)} 
               		 d\tau+C\delta_0\int_0^t\int_{\bbr} |(\psi_\x,\psi_{\x\x})|^2  
               		d\x d\tau+C\delta_0.
               	\end{aligned}
               \end{equation}               
               For $N_{22}$, using Lemma \ref{lem:variables-equivalence}, Young's inequality and Sobolev's inequality, we obtain
               \begin{equation}\label{N22}
               	\begin{aligned}
               		N_{22}&\leq 
               		C\int_0^t (\|\psi_\x\|_{L^\infty(\bbr)} \|\zeta_\x\|_{L^2(\bbr)} \|\psi_{\x\x}\|_{L^2(\bbr)} 
               		+\|\psi_\x\|^2_{L^\infty(\bbr)} \|(\phi_\xi,\vartheta_\xi)\|_{L^2(\bbr)} \|\zeta_{\x}\|_{L^2(\bbr)} )d\tau  \\
               		&\leq C(\varepsilon_1+\delta_0)
               		\int_0^t (\|\psi_\x\|_{L^\infty(\bbr)} \|\psi_{\x\x}\|_{L^2(\bbr)} 
               		+\|\psi_\x\|_{L^\infty(\bbr)}^2)d\tau  \\
               		&\leq C(\varepsilon_1+\delta_0)
               		\int_0^t \|\psi_\x\|^2_{H^1(\bbr)} d\tau.
               	\end{aligned}
               \end{equation}
               Similarly, for $N_{23}$, we can get
              \begin{equation}\label{N23}
               	\begin{aligned}
               		N_{23}&\leq 
               		C\int_0^t (\|\psi_\x\|_{L^\infty(\bbr)} \||\bar{v}_\x|\zeta_\x\|_{L^2(\bbr)} \|(\phi_\xi,\vartheta_\xi)\|_{L^2(\bbr)}+\|\psi_\x\|_{L^\infty(\bbr)} \|\psi_\xi\|_{L^2(\bbr)} \|\zeta_{\x}\|_{L^2(\bbr)} )d\tau  \\
               		&\leq C(\varepsilon_1+\delta_0)
               		\int_0^t (\|\psi_\x\|_{L^\infty(\bbr)} \||\bar{v}_\x|\zeta_\x\|_{L^2(\bbr)} 
               		+\|\psi_\x\|_{L^\infty(\bbr)} \|\psi_\xi\|_{L^2(\bbr)}  )d\tau  \\
               		&\leq C(\varepsilon_1+\delta_0)
               		\int_0^t \|\psi_\x\|^2_{H^1(\bbr)} d\tau
               		+ C(\varepsilon_1+\delta_0)\int_0^t\int_{\bbr}|(v^S)_{\xi}^{-\mathbf{X}}| |\zeta_\x|^2  
               		d\x d\tau  \\
               		&\quad +C(\varepsilon_1+\delta_0)\int_0^t(1+\tau)^{-2} \|\zeta_\x\|^2 _{L^2(\bbr)} 
               		d\tau .
               	\end{aligned}
               \end{equation}
               Combining \eqref{N21}-\eqref{N23}, we obtain the estimate for $N_2$ as follows 
               \begin{equation}\label{N2}
               	\begin{aligned}
               	    N_2
               	    &\leq C\int_0^t(\mathcal{G}^S+\mathcal{G}^R)  d\tau      	
               	    +C\int_0^t\int_{\bbr}|(v^S)_{\xi}^{-\mathbf{X}}| |\zeta_\x|^2  
               	    d\x d\tau
               	    +C\delta_0\int_0^t\int_{\bbr}|W_\x|^2  
               	    d\x d\tau\\
               	    &\quad +C\int_0^t(1+\tau)^{-\frac 32} \|\zeta_\x\|^2 _{L^2(\bbr)} 
               	    d\tau+C(\eps_1+\delta_0)\int_0^t\int_{\bbr} |(\psi_\x,\psi_{\x\x})|^2  
               	    d\x d\tau+C\delta_0.
               	    \end{aligned}
               \end{equation}
               For $N_3$, using Young's inequality, together with Lemmas \ref{lem:rare-v} and \ref{lem:wave-inter}, we have
               \begin{equation}\label{N3}
               	\begin{aligned}
               		N_3
               		&=C\int_0^t\int_{\bbr} 
               		|R_{3,\x}||\zeta_\x| d\x d\tau\\
               		&\leq C\int_0^t
               		\|R_{3,\x}\|_{L^2(\bbr)}  \|\zeta_\x\|^2_{L^2(\bbr)}  d\tau+C\int_0^t
               		\|R_{3,\x}\|_{L^2(\bbr)}  d\tau\\
               		&\leq C\int_0^t \left[\deltar^{\frac 12}(1+\tau)^{-\frac 32} + \deltas \deltar \left(e^{-C \deltas \tau}+|(1+\deltar \tau)|^{-\frac{5}{2}}\right)\right]\|\zeta_\x\|^2 _{L^2(\bbr)} 
               		d\tau\\
               		&\quad +C\int_0^t \left[\deltar^{\frac 12}(1+\tau)^{-\frac 32} + \deltas \deltar \left(e^{-C \deltas \tau}+|(1+\deltar \tau)|^{-\frac{5}{2}}\right)\right]
               		d\tau\\
               		&\leq C\int_0^t \left[\deltar^{\frac 12}(1+\tau)^{-\frac 32} + \deltas \deltar \left(e^{-C \deltas \tau}+|(1+\deltar \tau)|^{-\frac{5}{2}}\right)\right]\|\zeta_\x\|^2 _{L^2(\bbr)} 
               		d\tau+C\delta_0^{\frac 12}.
               	\end{aligned}
               \end{equation}
               Therefore, by Lemma \ref{Re-en}, using Gronwall's inequality, together with \eqref{N1}, \eqref{N2} and \eqref{N3}, we obtain
               \begin{align*}
               		\|\zeta_\x(t,\cdot)\|_{L^2(\mathbb{R})}^2 
               		&\le 
               		C\|(W,\psi,\zeta)(0,\cdot)\|_{H^1(\mathbb{R})}^2+
               		C\int_0^t\int_{\bbr}|(v^S)_{\xi}^{-\mathbf{X}}| |\zeta_\x|^2  
               		d\x d\tau  \\
               		&\quad+C\delta_0\int_0^t\int_{\bbr} |W_\x|^2  
               		d\x d\tau  
               		+C(\varepsilon_1+\delta_0)
               		\int_0^t \|\psi_{\x\x}\|^2_{L^2(\bbr)} d\tau+C\delta_0^{\frac 14}.
               	\end{align*}
            \end{proof}
            
               \begin{lemma}\label{1-d}
               		Under the hypotheses of Proposition \ref{prop:priori}, there exists $C > 0$ such that for all $t \in (0, T]$,
               		\begin{align*}
               			&\|(W_\x,\psi_\x,\zeta_\x)(t,\cdot)\|_{L^2(\mathbb{R})}^2 +\int_0^t   \|W_\x\|^2_{L^2(\bbr)} d\tau+
               			\int_0^t \|\psi_{\x\x}\|_{L^2(\mathbb{R})}^2  d\tau\\
               			&\le
               			C\|(W,\psi,\zeta)(0,\cdot)\|_{H^1(\mathbb{R})}^2
               			+C\int_0^t\int_{\bbr}|(v^S)_{\xi}^{-\mathbf{X}}| |\zeta_\x|^2  
               			d\x d\tau+C\delta_0^{\frac 14},            		
               		\end{align*}
               \end{lemma}
               \begin{proof}
               	By Lemmas \ref{flux}-\ref{S}, using Gronwall's inequality, together with the smallness of $\delta_0$ and $\chi$, we have 
               	\begin{equation}\label{eng-1}
               		\begin{aligned}
               			&\|(W_\x,\psi_\x,\zeta_\x)(t,\cdot)\|_{L^2(\mathbb{R})}^2 +\int_0^t   \|W_\x\|^2_{L^2(\bbr)} d\tau
               			+\int_0^t   \|\mathcal{F}_\x\|^2_{L^2(\bbr)} d\tau\\
               			&\le C\|(W,\psi,\zeta)(0,\cdot)\|_{H^1(\mathbb{R})}^2
               			+C\int_0^t\int_{\bbr}|(v^S)_{\xi}^{-\mathbf{X}}| |\zeta_\x|^2  
               			d\x d\tau \\
               			&\quad +C (\varepsilon_1+\delta_0) \int_0^t\int_{\bbr}  |\psi_{\x\x}|^2d\x d\tau+C\delta_0^{\frac 14}. 
               		\end{aligned}
               	\end{equation}
               	  Next, we aim to obtain the space-time estimate for 
               	  $\psi_{\x\x}$. From $\eqref{PUS}_2$, we can get the expression for $\psi_{\x\x}$:
               	   \begin{equation*}
               	   \mu\frac{\psi_{\x\x}}{v}= \psi_t - \sigma \psi_{\xi} - \dot{\mathbf{X}}(t) (u^S)_{\xi}^{-\mathbf{X}} + W_{\xi} + \mu\frac{v_\x}{v^2}\psi_\x- \mu\left[ \bar u_\x \left( \frac{1}{v} - \frac{1}{\bar{v}} \right)\right]_{\xi}+R_2.
               	   \end{equation*}
                Then, we have 
               	   	\begin{align*}
               	   	\|\psi_{\x\x}\|_{L^2(\mathbb{R})}^2   
               	   	\le  &C(
               	   	\|\psi_t \|_{L^2(\mathbb{R})}^2   
               	   	+\|\psi_{\x}\|_{L^2(\mathbb{R})}^2   
               	   	+\deltas|\dot{\mathbf{X}}(t)|^2
               	   	+\|W_{\x}\|_{L^2(\mathbb{R})}^2
                    +\|\bar u_\xi\bar v_\xi\phi\|_{L^2(\mathbb{R})}^2\\  
               	   	&+\|v_{\x}\psi_\x\|_{L^2(\mathbb{R})}^2  
               	   	+\|\bar u_{\x\x} \phi\|_{L^2(\mathbb{R})}^2 
               	   	+\|\bar u_\x\phi_{\x}\|_{L^2(\mathbb{R})}^2
               	   	+\|R_2\|_{L^2(\mathbb{R})}^2),
               	   	\end{align*}
               	   which implies
               	   \begin{align*}
               	   	\int_0^t \|\psi_{\x\x}\|_{L^2(\mathbb{R})}^2  d\tau
               	   	\le&  C \int_0^t (
               	   	\|\psi_t \|_{L^2(\mathbb{R})}^2   
               	   	+\|\psi_{\x}\|_{L^2(\mathbb{R})}^2   
               	   	+\deltas|\dot{\mathbf{X}}|^2
               	   	+\|W_{\x}\|_{L^2(\mathbb{R})}^2
                    +\|\bar u_\xi\bar v_\xi\phi\|_{L^2(\mathbb{R})}^2) d\tau\\
               	   	&+C\int_0^t (\|v_{\x}\psi_\x\|_{L^2(\mathbb{R})}^2  
               	   	+ \|\bar u_{\x\x} \phi\|_{L^2(\mathbb{R})}^2 
               	   	+\|\bar u_\x\phi_{\x}\|_{L^2(\mathbb{R})}^2+\|R_2\|_{L^2(\mathbb{R})}^2) d\tau.
               	   	\end{align*}
               	   By Lemmas \ref{Re-en}, \ref{lem:rare-v}, \ref{lem:shock-v}, \ref{lem:variables-equivalence}, \ref{lem:wave-inter} and \ref{lem:variables-weighted}, using \eqref{psi-t} and Sobolev's inequality, we get 
               	   	\begin{align*}
               	   		\int_0^t \|\psi_{\x\x}\|_{L^2(\mathbb{R})}^2  d\tau
               	   		\le
               	   		&C\|(W,\psi,\zeta)(0,\cdot)\|_{H^1(\mathbb{R})}^2
               	   		+C\int_0^t\int_{\bbr}|(v^S)_{\xi}^{-\mathbf{X}}| |\zeta_\x|^2  
               	   		d\x d\tau \\
               	   		& +C \int_0^t (
               	   		\|\psi_t \|_{L^2(\mathbb{R})}^2
               	   		+\|v_{\x}\psi_\x\|_{L^2(\mathbb{R})}^2    
               	   		+\|W_{\x}\|_{L^2(\mathbb{R})}^2 ) d\tau\\
               	   		&+C\int_0^t (1+\tau)^{-2}\|\zeta_\x\|^2 _{L^2(\bbr)} 
               	   		d\tau+C\delta_0^{\frac 14}\\
               	   		\le
               	   		&C\|(W,\psi,\zeta)(0,\cdot)\|_{H^1(\mathbb{R})}^2
               	   		+C\int_0^t\int_{\bbr}|(v^S)_{\xi}^{-\mathbf{X}}| |\zeta_\x|^2  
               	   		d\x d\tau \\
               	   		& +C \int_0^t (
               	   		\|\psi_t \|_{L^2(\mathbb{R})}^2
               	   		+\|\phi_{\x}\|_{L^2(\mathbb{R})}^2  \|\psi_{\x}\|_{L^\infty(\mathbb{R})}^2 
               	   		+\|W_{\x}\|_{L^2(\mathbb{R})}^2 ) d\tau\\
               	   		&+C\int_0^t (1+\tau)^{-2}\|\zeta_\x\|^2 _{L^2(\bbr)} 
               	   		d\tau+C\delta_0^{\frac 14}\\
               	   		\le
               	   		&C\|(W,\psi,\zeta)(0,\cdot)\|_{H^1(\mathbb{R})}^2
               	   		+C\int_0^t\int_{\bbr}|(v^S)_{\xi}^{-\mathbf{X}}| |\zeta_\x|^2  
               	   		d\x d\tau \\
               	   		& +C \int_0^t (
               	   		\|\psi_t \|_{L^2(\mathbb{R})}^2
               	   		+(\varepsilon_1+\delta_0) \|\psi_{\x}\|_{H^1(\mathbb{R})}^2 
               	   		+\|W_{\x}\|_{L^2(\mathbb{R})}^2 ) d\tau\\
               	   		&+C\int_0^t (1+\tau)^{-2}\|\zeta_\x\|^2 _{L^2(\bbr)} 
               	   		d\tau+C\delta_0^{\frac 14}.
               	   	\end{align*}
               	Using the smallness of $\varepsilon_1$ and $\delta_0$, we have 
               	\begin{align*}
               			\int_0^t \|\psi_{\x\x}\|_{L^2(\mathbb{R})}^2  d\tau
               			\le
               			&C\|(W,\psi,\zeta)(0,\cdot)\|_{H^1(\mathbb{R})}^2
               			+C\int_0^t\int_{\bbr}|(v^S)_{\xi}^{-\mathbf{X}}| |\zeta_\x|^2  
               			d\x d\tau \\
               			& +C \int_0^t (
               			\|\psi_t \|_{L^2(\mathbb{R})}^2
               			+\|W_{\x}\|_{L^2(\mathbb{R})}^2 ) d\tau\\
               			&+C\int_0^t (1+\tau)^{-2}\|\zeta_\x\|^2 _{L^2(\bbr)} 
               			d\tau+C\delta_0^{\frac 14}.
               	\end{align*}               	
               	Using Gronwall's inequality, together with \eqref{eng-1} and \eqref{psi-t}, we obtain
               		\begin{align*}
               			&\|(W_\x,\psi_\x,\zeta_\x)(t,\cdot)\|_{L^2(\mathbb{R})}^2 +\int_0^t   \|W_\x\|^2_{L^2(\bbr)} d\tau+
               			\int_0^t \|\psi_{\x\x}\|_{L^2(\mathbb{R})}^2  d\tau\\
               			&\le
               			C\|(W,\psi,\zeta)(0,\cdot)\|_{H^1(\mathbb{R})}^2
               			+C\int_0^t\int_{\bbr}|(v^S)_{\xi}^{-\mathbf{X}}| |\zeta_\x|^2  
               			d\x d\tau+C\delta_0^{\frac 14},            		
               		\end{align*}
               	where we use the smallness of $\varepsilon_1$ and $\delta_0$.
               \end{proof}

 \subsection{Weighted space-time estimates for the entropy derivative}\label{subsec:weighted}
As indicated by Lemma \ref{1-d}, closing the energy estimate
requires weighted space-time estimates for the first-order derivative of the entropy. Before proceeding to the formal proof, we first introduce the following three auxiliary lemmas.
    \begin{lemma}\label{trans}
    For all $t \in (0, T]$, let $r(t,\x)$ be an arbitrary function satisfying 
    $$
    \int_0^t  \int_{\mathbb{R}} | (v^S)_\xi^{-\mathbf{X}}| r^2(\tau,\xi) d\xi d\tau < +\infty,
    $$
    with $(v^S)_\xi^{-\mathbf{X}}$ as defined in Lemma \ref{lem:shock-v}. If the function $h(t,\x)$ satisfies
     $$0\leq h(t,\xi)\leq |(v^S)_\xi^{-\mathbf{X}}|| r(t,\xi)|,$$ 	 
    then there exists $C > 0$ such that for all $t \in (0, T]$,
    \begin{equation*}
      I[h]\leq  C\deltas\int_0^t  \int_{\mathbb{R}} | (v^S)_\xi^{-\mathbf{X}}| r^2(\tau,\xi) d\xi d\tau,
    \end{equation*}
    where 
     $$I[h]:=\int_0^t  \int_{\mathbb{R}} | (v^S)_\xi^{-\mathbf{X}} | \left|\int_0^{\tau}h(\eta, \xi+\sigma (\tau-\eta)) d\eta\right|^2 d\xi d\tau.$$
    \end{lemma}
    \begin{proof}
        Set $$ \tilde{I}[h] :=\left|\int_0^{\tau}h(\eta, \xi+\sigma (\tau-\eta)) d\eta\right|^2.$$
        Using the smallness of $\eps_1$ and $\delta_0$, we have $|\dot{\mathbf X}(t)|
    \le C_0(\varepsilon_1+\delta_0)
    \le \frac{\sigma}{2}$.
Then, the change of variables
$y=\xi+\sigma(\tau-\eta)-\mathbf X(\eta)$,
together with Lemma \ref{lem:shock-v}, gives
    \begin{align*}
        &\int_0^\tau
        \left|
            (v^S)'
            \bigl(
                \xi+\sigma(\tau-\eta)-\mathbf X(\eta)
            \bigr)
        \right|d\eta\\
        &\quad\le
        \frac{2}{\sigma}
        \int_0^\tau
        \left|
            (v^S)'
            \bigl(
                \xi+\sigma(\tau-\eta)-\mathbf X(\eta)
            \bigr)
        \right|
        \bigl(\sigma+\dot{\mathbf X}(\eta)\bigr)d\eta\\
        &\quad=
        \frac{2}{\sigma}
        \int_{\xi-\mathbf X(\tau)}^{
            \xi+\sigma\tau-\mathbf X(0)
        }
        |(v^S)'(y)|d y\le
        \frac{2}{\sigma}
        \int_{\mathbb{R}}|(v^S)'(y)|d y
        \le C\delta_S.
    \end{align*}
        Then, using H\"older's inequality, we have
        \begin{align*}
           \tilde{I}[h] &\leq \left(\int_0^{\tau} |h(\eta, \xi+\sigma (\tau-\eta))| d\eta\right)^2\\
            &\leq  \left(\int_0^{\tau} |(v^S)'(\xi+\sigma (\tau-\eta)-\mathbf{X}(\eta))| d\eta\right) \\
            &\cdot\left(\int_0^{\tau} |(v^S)'(\xi+\sigma (\tau-\eta)-\mathbf{X}(\eta))| r^2(\eta, \xi+\sigma (\tau-\eta))d\eta\right)\\
            &\leq C\deltas\int_0^{\tau} |(v^S)'(\xi+\sigma (\tau-\eta)-\mathbf{X}(\eta))| r^2(\eta, \xi+\sigma (\tau-\eta))d\eta.
        \end{align*}
       Applying the Fubini-Tonelli theorem and the change of variables $y=\xi+\sigma(\tau-\eta)$, we can get the following estimate
        \begin{align*}
            I[h]&\leq C\deltas \int_0^t\int_{\eta}^t\int_{\mathbb{R}}|(v^S)'(\xi-\mathbf{X}(\tau))| |(v^S)'(\xi+\sigma (\tau-\eta)-\mathbf{X}(\eta))| \\
            &\qquad\qquad\qquad \cdot r^2(\eta, \xi+\sigma (\tau-\eta)) d\xi d\tau d\eta\\
            &\leq C\deltas \int_0^t\int_{\eta}^t\int_{\mathbb{R}}|(v^S)'(y-\sigma(\tau-\eta)-\mathbf{X}(\tau))| |(v^S)'(y-\mathbf{X}(\eta))| r^2(\eta, y) dy d\tau d\eta\\
            &\leq C\deltas \int_0^t\int_{\mathbb{R}}|(v^S)'(y-\mathbf{X}(\eta))| r^2(\eta, y)\int_{\eta}^t  |(v^S)'(y-\sigma(\tau-\eta)-\mathbf{X}(\tau))|  d\tau dy d\eta\\      
            &\leq C\deltas \int_0^t\int_{\mathbb{R}}|(v^S)'(y-\mathbf{X}(\eta))| r^2(\eta, y) dy d\eta,          	
        \end{align*}
        so that we can complete the proof.
    \end{proof}
            		
	\begin{lemma}\label{trans-R}
        For all $t \in (0, T]$, let $r(t,\x)$ be an arbitrary function satisfying 
        $$
        \int_0^t  \int_{\mathbb{R}}  r^2(\tau,\xi) d\xi d\tau < +\infty.
        $$ 			
        If the function $h(t,\x)$ satisfies
        $$0\leq h(t,\xi)\leq |v^R_\x| | r(t,\xi)|,$$ 	 
        with $(v^R)_\xi$ as defined in Lemma \ref{lem:rare-v}, then there exists $C > 0$ such that for all $t \in (0, T]$,
        \begin{equation*}
        I[h]\leq  C\deltas\int_0^t  \int_{\mathbb{R}} r^2(\tau,\xi) d\xi d\tau,
        \end{equation*}
        where 
        $$I[h]:=\int_0^t  \int_{\mathbb{R}} | (v^S)_\xi^{-\mathbf{X}} | \left|\int_0^{\tau}h(\eta, \xi+\sigma (\tau-\eta)) d\eta\right|^2 d\xi d\tau.$$
        \end{lemma}
    \begin{proof}
    Let $$ \tilde{I}[h] :=\left|\int_0^{\tau}h(\eta, \xi+\sigma (\tau-\eta)) d\eta\right|^2,$$
    and then using H\"older's inequality and Lemma \ref{lem:rare-v}, we have
    \begin{align*}
        \tilde{I}[h] &\leq \left(\int_0^{\tau} |h(\eta, \xi+\sigma (\tau-\eta))| d\eta\right)^2\\
        &\leq  \left(\int_0^{\tau} |u^R_\x(\eta, \xi+\sigma (\tau-\eta))|^2 d\eta\right) \left(\int_0^{\tau}  r^2(\eta, \xi+\sigma (\tau-\eta))d\eta\right)\\
        &\leq  \left(\int_0^{\tau} \|u^R_\x(\eta,\cdot)\|^2_{L^\infty(\bbr)} d\eta\right) \left(\int_0^{\tau}  r^2(\eta, \xi+\sigma (\tau-\eta))d\eta\right)\\
        &\leq C\int_0^{\tau}  r^2(\eta, \xi+\sigma (\tau-\eta))d\eta.
    \end{align*}
    Using the Fubini-Tonelli theorem, \eqref{dxbound}, and
Lemma \ref{lem:shock-v}, together with the change of variable
$y=\xi+\sigma(\tau-\eta)$, we obtain the following estimate:
    \begin{align*}
        I[h]&\leq C \int_0^t\int_{\eta}^t\int_{\mathbb{R}}|(v^S)'(\xi-\mathbf{X}(\tau))| r^2(\eta, \xi+\sigma (\tau-\eta)) d\xi d\tau d\eta\\
        &\leq C \int_0^t\int_{\eta}^t\int_{\mathbb{R}}|(v^S)'(y-\sigma(\tau-\eta)-\mathbf{X}(\tau))|r^2(\eta, y) dy d\tau d\eta\\
        &\leq C \int_0^t\int_{\mathbb{R}}r^2(\eta, y)\int_{\eta}^t  |(v^S)'(y-\sigma(\tau-\eta)-\mathbf{X}(\tau))|  d\tau dy d\eta\\ &\leq C\deltas \int_0^t\int_{\mathbb{R}}r^2(\eta, y) dy d\eta,
    \end{align*}
    so that we can complete the proof.
    \end{proof}
            		
    \begin{lemma}\label{trans-1}
    For all $t \in (0, T]$, let $h(t,\x)$ be an arbitrary function satisfying 
    $$
    \int_0^t\|h(\tau,\cdot)\|_{L^2(\bbr)}\,d\tau < +\infty,
    $$
    with $(v^S)_\xi^{-\mathbf{X}}$ as defined in Lemma \ref{lem:shock-v}.
    Then, there exists $C > 0$ such that for all $t \in (0, T]$,
    $$
    I[h]
    \le C\deltas\left(\int_0^t\|h(\tau,\cdot)\|_{L^2(\bbr)}\,d\tau\right)^2,$$
    where 
    $$I[h]=\int_0^t  \int_{\mathbb{R}} | (v^S)_\xi^{-\mathbf{X}} | \left|\int_0^{\tau}h(\eta, \xi+\sigma (\tau-\eta)) d\eta\right|^2 d\xi d\tau,$$ 
    \end{lemma}
    \begin{proof}
    We first make the change of variable $y=\xi+\sigma \tau.$
    Then,
    \[
    I[h]
    =
    \int_0^t\int_{\mathbb R}
    |(v^S)'(y-\sigma \tau-\mathbf{X}(\tau))|
    \left|
    \int_0^{\tau} h(\eta, y-\sigma\eta)\,d\eta
    \right|^2
    dy d\tau.
    \]
    Since
    \[
    \left|
    \int_0^{\tau} h(\eta, y-\sigma\eta)\,d\eta
    \right|
    \le
    \int_0^t |h(\eta, y-\sigma\eta)|\,d\eta,
    \]
    we get
    \[
    I[h]
    \le
    \int_0^t \int_{\mathbb R}
    |(v^S)'(y-\sigma \tau-\mathbf{X}(\tau))|
    \left(
    \int_0^t |h(\eta, y-\sigma\eta)|\,d\eta
    \right)^2
    dy d\tau.
    \]
    By the Fubini-Tonelli theorem, we may exchange the \(\tau\)- and \(y\)-integrations:
    \[
    I[h]
    \le
    \int_{\mathbb R}
    \left(
    \int_0^t |h(\eta, y-\sigma\eta)|\,d\eta
    \right)^2
    \left(
    \int_0^t  |(v^S)'(y-\sigma \tau-\mathbf{X}(\tau))|\,d\tau
    \right)
    dy.
    \]
    Consequently, using Lemma \ref{lem:shock-v} and \eqref{dxbound}, we have
    \[
    I[h]
    \le C\deltas
    \int_{\mathbb R}
    \left(
    \int_0^t |h(\eta, y-\sigma\eta)|\,d\eta
    \right)^2
    dy.
    \]
    Using Minkowski's integral inequality, we get
    \[
    \left\|
    \int_0^t |h(\eta, \cdot-\s\eta)|\,d\eta
    \right\|_{L^2_y}
    \le
    \int_0^t
    \|h(\eta, \cdot-\sigma\eta)\|_{L^2_y}\,d\eta.
    \]
    Since spatial translations preserve the \(L^2\)-norm, we obtain
    \[
    I[h] \le
    C\deltas\left(\int_0^t\|h(\eta,\cdot)\|_{L^2(\bbr)}\,d\eta\right)^2.   \]
    \end{proof}
            		
    Applying Lemmas \ref{trans}-\ref{trans-1}, we can derive the following weighted space-time estimates for the first derivative of the entropy.
     \begin{lemma}\label{weight-s}
        Under the hypotheses of Proposition \ref{prop:priori}, there exists $C > 0$ such that for all $t \in (0, T]$,
        \begin{align*}
         & \int_0^t  \int_{\mathbb{R}} | (v^S)_\xi^{-\mathbf{X}} |  \zeta^2_{\x} (\tau,\x)d\xi d\tau\\
        &\leq  C\|(W,\psi,\zeta)(0,\cdot)\|_{H^1(\mathbb{R})}^2 
         +C\deltas\int_0^t\int_{\bbr}  |v^R_\xi|^2(|W_\xi|^2 + |\zeta_\xi|^2) d\x d\tau\\
        &\quad+C\delta_0 \int_0^t  \int_{\mathbb{R}} |W_\x|^2d\xi d\tau+
        C(\delta_0+\eps_1)
        \left(
        \int_0^t
        \|\psi_{\x\x}\|^2_{L^2(\bbr)}\,d\tau
        \right)
        +
    	C\delta_0^{\frac 14}.
    \end{align*}
     \end{lemma}  
     \begin{proof}
    For convenience, define the initial entropy perturbation by
    \begin{equation*}
        \zeta_0(\xi):= s(0, \xi) - \bar{s}(0, \xi).
     \end{equation*}
    Differentiating the equation $\eqref{PUS}_3$ with respect to $\xi$ and using \eqref{H}, we have 
    \begin{equation}\label{en-d}
        \zeta_{\x t} -  \sigma\zeta_{\x\x}  = \dot{\mathbf{X}}(t) (s^S)_{\xi\xi}^{-\mathbf{X}}+H_{\x}-R_{3,\x}.
    \end{equation}
     Using the method of characteristics, we can obtain the following representation of the solution to \eqref{en-d}:
     \begin{align*}
        \zeta_{\x} (t,\x)=&  \zeta_{0\xi}(\xi+\sigma t) +\int_0^t \dot{\mathbf{X}}(\tau) (s^S)_{\xi\xi}^{-\mathbf{X}}(\xi+\sigma (t-\tau)) d\tau\\
        &+\int_0^t H_{\x}(\tau, \xi+\sigma (t-\tau)) d\tau-\int_0^t R_{3,\x}(\tau, \xi+\sigma (t-\tau)) d\tau.
    \end{align*}
                    It follows that
                    \begin{equation}\label{en-i}
                    	\begin{aligned}
                    	& \int_0^t  \int_{\mathbb{R}} | (v^S)_\xi^{-\mathbf{X}} |  \zeta^2_{\x} (\tau,\x)d\xi d\tau\\
                    	 \leq &  C\int_0^t  \int_{\mathbb{R}} | (v^S)_\xi^{-\mathbf{X}} |  \zeta^2_{0\xi}(\xi+\sigma \tau) d\xi d\tau\\
                    	 &+C\int_0^t  \int_{\mathbb{R}} | (v^S)_\xi^{-\mathbf{X}} | \left|\int_0^{\tau} \dot{\mathbf{X}}(\eta) (s^S)_{\xi\xi}^{-\mathbf{X}}(\xi+\sigma (\tau-\eta)) d\eta\right|^2 d\xi d\tau\\
                    	 & +C\int_0^t  \int_{\mathbb{R}} | (v^S)_\xi^{-\mathbf{X}} | \left|\int_0^{\tau} H_{\x}(\eta, \xi+\sigma (\tau-\eta)) d\eta\right|^2 d\xi d\tau \\
                    	 &+C\int_0^t  \int_{\mathbb{R}} | (v^S)_\xi^{-\mathbf{X}} | \left|\int_0^{\tau} R_{3,\x}(\eta, \xi+\sigma (\tau-\eta)) d\eta\right|^2 d\xi d\tau\\
                    	 &:= I_1 + I_2 + I_3+I_4.
                    \end{aligned}
                    \end{equation}
                    Next, we estimate each term of \eqref{en-i} one by one. For $I_1$, applying the Fubini-Tonelli theorem and the change of variables $y=\xi+\sigma \tau$, we have
                   \begin{equation}\label{I1}
                   	  \begin{aligned}
                    	 I_1&=C\int_0^t  \int_{\mathbb{R}} | (v^S)'(y-\sigma \tau-\mathbf{X}(\tau))|  \zeta^2_{0y}(y) dy d\tau\\
                    	 &\leq C\int_{\mathbb{R}}  \left(\int_0^t | (v^S)'(y-\sigma \tau-\mathbf{X}(\tau))| d\tau \right)\zeta^2_{0y}(y) dy\\
                    	 &\leq C \deltas \int_{\mathbb{R}}\zeta^2_{0y}(y) dy=C\deltas \| \zeta_{0\xi} \|_{L^2(\bbr)}^2.
                       \end{aligned}
                    \end{equation}
                  For $I_2$, applying Lemmas \ref{lem:shock-v} and \ref{trans} with $r(t,\xi)=\dot{\mathbf{X}}(t)$, we obtain
                   \begin{equation}\label{I2}
                   	\begin{aligned}
                   		I_2 &\leq C\deltas\int_0^t  \int_{\mathbb{R}} | (v^S)_\xi^{-\mathbf{X}}| | \dot{\mathbf{X}}(\tau)|^2d\xi d\tau\\
                   		&\leq C\deltas\int_0^t  \int_{\mathbb{R}} | (v^S)_\xi^{-\mathbf{X}}|d\xi | \dot{\mathbf{X}}(\tau)|^2d\tau\\
                   		&\leq C\deltas^2\int_0^t  | \dot{\mathbf{X}}(\tau)|^2d\tau.
                   	\end{aligned}
                   \end{equation}
                   For $I_3$,  we first partition the terms in $H_{\x}$ into three groups as follows                 
                   	\begin{align*}
                   		H_\xi =\ &\left[ \frac{2\mu}{v\theta} \bar{u}_\xi \psi_{\xi\xi}+ \frac{2\mu}{v\theta} \bar{u}_{\xi\xi} \psi_\xi - 2\mu \left( \frac{\phi \bar{\theta} + \bar{v} \vartheta + \phi \vartheta}{v\theta \bar{v}\bar{\theta}} \right) \bar{u}_\xi \bar{u}_{\xi\xi}\right.\\
                   		&\left.+ \mu \left( \frac{\bar{v}_\xi \bar{\theta} + \bar{v} \bar{\theta}_\xi}{\bar{v}^2 \bar{\theta}^2} - \frac{(\bar{v}_\xi + \phi_\xi)\theta + v(\bar{\theta}_\xi + \vartheta_\xi)}{v^2 \theta^2} \right) \bar{u}_\xi^2 - 2\mu \frac{\bar{v}_\xi \theta + v\bar{\theta}_\xi }{v^2 \theta^2} \bar{u}_\xi \psi_\xi\right]\\
                   		& +\left[ \frac{2\mu}{v\theta} \psi_\xi \psi_{\xi\xi} - \mu \frac{ \phi_\xi\theta + v\vartheta_\xi}{v^2 \theta^2} \psi_\xi^2\right]- \left[2\mu \frac{\phi_\xi\theta + v\vartheta_\xi}{v^2 \theta^2} \bar{u}_\xi \psi_\xi + \mu \frac{\bar{v}_\xi\theta + v\bar{\theta}_\xi}{v^2 \theta^2} \psi_\xi^2\right]\\
                   		&:= \tilde I_{31}+\tilde I_{32}+\tilde I_{33}.
                \end{align*}    
                 For simplicity, we denote the integral terms corresponding to $\tilde{I}_{31}, \tilde{I}_{32}, \tilde{I}_{33}$ by $I_{31}, I_{32}, I_{33}$, respectively.      
                   In order to estimate $I_{31}$, we start by estimating $\tilde I_{31}$. For $\tilde I_{31}$, using Lemmas \ref{lem:shock-v} and \ref{lem:rare-v}, we have 
                   \begin{equation*}
                   	|\tilde I_{31}|
                   	\leq C(|\bar{v}_\xi|  |\psi_{\xi\xi}| + |\bar{v}_\xi|  |\psi_{\xi}| +|\bar{v}_\xi|^2 \big( |\phi| + |\vartheta| \big)+|\bar{v}_\xi|^2 ( |\phi_\xi| + |\vartheta_\xi|) ).                   
                   \end{equation*}  
                   As the composite wave consists of a rarefaction wave and a viscous shock wave, using Lemma \ref{lem:shock-v}, we divide $\tilde I_{31}$ into two parts as follows
                   	\begin{align*}
                   		|\tilde I_{31}|
                   		&\leq C(|v^R_\xi|  |\psi_{\xi\xi}| + |v^R_\xi|  |\psi_{\xi}| +|v^R_\xi|^2 \big( |\phi| + |\vartheta| \big)+|v^R_\xi|^2 ( |\phi_\xi| + |\vartheta_\xi|) )\\
                   		&\quad+ C|(v^S)^{-{\mb X}}_\x|  (|\psi_{\xi\xi}| +   |\psi_{\xi}| +\big( |\phi| + |\vartheta| \big)+ ( |\phi_\xi| + |\vartheta_\xi|) )\\
                   		&\quad := \tilde I^R_{31}+\tilde I^S_{31}.     
                   	\end{align*}  
                   Then, we divide $I_{31}$ into two corresponding components, designating them as $I^R_{31}$ and $I^S_{31}$. For  $I^R_{31}$, from Lemmas \ref{lem:rare-v} and \ref{trans-R}, we get 
                   \begin{equation}\label{IR31}
                   	\begin{aligned}
                   		I^R_{31} 
                   		&\leq C \int_0^t  \int_{\mathbb{R}} | (v^S)_\xi^{-\mathbf{X}} | \left|\int_0^{\tau} |v^R_\xi| ( |\psi_{\xi\xi}| + |\psi_{\xi}| + |v^R_\xi|(|\phi| + |\vartheta| +|\phi_\xi| + |\vartheta_\xi|)) d\eta\right|^2 d\xi d\tau \\
                   		&\leq C\deltas\int_0^t\int_{\bbr}  ( |\psi_{\xi\xi}| ^2+ |\psi_{\xi}| ^2+ |v^R_\xi|(|\phi|^2 + |\vartheta|^2) +|v^R_\xi|^2(|\phi_\xi|^2 + |\vartheta_\xi|^2)) d\x d\tau.                   		
                   	\end{aligned}
                   \end{equation}
                    For $I^S_{31}$,  by Lemmas \ref{lem:shock-v} and \ref{trans}, we have 
                   \begin{equation}\label{IS31}
                    	\begin{aligned}
                    		I^S_{31} 
                    		&\leq C \int_0^t  \int_{\mathbb{R}} | (v^S)_\xi^{-\mathbf{X}} | \left|\int_0^{\tau}  | (v^S)_\xi^{-\mathbf{X}} | ( |\psi_{\xi\xi}| + |\psi_{\xi}| + |\phi| + |\vartheta| +|\phi_\xi| + |\vartheta_\xi|) d\eta\right|^2 d\xi d\tau \\
                    		&\leq C\deltas\int_0^t  \int_{\mathbb{R}} | (v^S)_\xi^{-\mathbf{X}}| (| \psi_\x|^2+| \psi_{\x\x}|^2+|\phi|^2+|\vartheta|^2+|\phi_\x|^2+|\vartheta_\x|^2)d\xi d\tau\\
                    		& \leq C\deltas\int_0^t  \int_{\mathbb{R}} (| \psi_\x|^2+| \psi_{\x\x}|^2)d\xi d\tau+C\deltas\int_0^t  \int_{\mathbb{R}} | (v^S)_\xi^{-\mathbf{X}}| (|\phi|^2+|\vartheta|^2)d\xi d\tau\\
                    		&\quad +C\deltas\int_0^t  \int_{\mathbb{R}} | (v^S)_\xi^{-\mathbf{X}}| (|\phi_\x|^2+|\vartheta_\x|^2)d\xi d\tau.
                    	\end{aligned}
                    \end{equation}
                    From \eqref{IR31} and \eqref{IS31}, we have
                    \begin{equation}\label{I31}
                    	\begin{aligned}
                    		I_{31}& \leq I^R_{31}+I^S_{31}\\
                    		& \leq C\deltas\int_0^t  \int_{\mathbb{R}} (| \psi_\x|^2+| \psi_{\x\x}|^2)d\xi d\tau+C\deltas\int_0^t(\mathcal{G}^S+\mathcal G^R)d\tau\\
                    		&\quad +C\deltas\int_0^t  \int_{\mathbb{R}} | (v^S)_\xi^{-\mathbf{X}}| (|\phi_\x|^2+|\vartheta_\x|^2)d\xi d\tau\\
                    		&\quad +C\deltas\int_0^t\int_{\bbr}  |v^R_\xi|^2(|\phi_\xi|^2 + |\vartheta_\xi|^2) d\x d\tau.             
                    	\end{aligned}
                    \end{equation}
                   Similarly, for $\tilde I_{32}$, we have 
                    \begin{align*}
                   		|\tilde I_{32}|\leq C(|\psi_\xi|  |\psi_{\xi\xi}| + |(\phi_\xi,\vartheta_\xi)|  |\psi_{\xi}|^2).
                   \end{align*}  
                   Then, by Lemma \ref{trans-1}, we get
                    \begin{align*}
                    		I_{32}
                    		&\le C \int_0^t  \int_{\mathbb{R}} | (v^S)_\xi^{-\mathbf{X}} | \left|\int_0^{\tau} (|\psi_\xi|  |\psi_{\xi\xi}| + |(\phi_\xi,\vartheta_\xi)|  |\psi_{\xi}|^2)(\eta, \xi+\sigma (\tau-\eta)) d\eta\right|^2 d\xi d\tau \\
                    		&\le
                    		C\deltas
                    		\left(
                    		\int_0^t
                    		\|\psi_\xi\psi_{\xi\xi}\|_{L^2(\bbr)}+ \|(\phi_\xi,\vartheta_\xi)  |\psi_{\xi}|^2\|_{L^2(\bbr)}\,d\tau
                    		\right)^2 .
                   \end{align*}
                    Using
                    \[
                    \|\psi_\xi\psi_{\xi\xi}\|_{L^2(\bbr)}
                    \le
                    \|\psi_\xi\|_{L^\infty(\bbr)}\|\psi_{\xi\xi}\|_{L^2(\bbr)}
                    \]
                    and 
                     \[
                   \|(\phi_\xi,\vartheta_\xi)  |\psi_{\xi}|^2\|_{L^2(\bbr)}
                    \le
                    \|\psi_\xi\|^2_{L^\infty(\bbr)}\|(\phi_{\xi},\vartheta_{\x})\|_{L^2(\bbr)},
                    \]
                    together with Sobolev's inequality and Lemma \ref{lem:variables-equivalence}, we have
                     \begin{equation}\label{I32}
                     	\begin{aligned}
                         I_{32}
                     		&\le
                     		C\deltas
                     		\left(
                     		\int_0^t
                     		\|\psi_\xi\|^2_{L^2(\bbr)}+ \| \psi_{\xi\x}\|^2_{L^2(\bbr)}\,d\tau
                     		\right)^2\\
                     		&\le C\deltas
                     		\left(
                     		\int_0^t
                     		\|\psi_\xi\|^2_{L^2(\bbr)}d\tau\right)^2+
                     		C\deltas\left(\int_0^t \| \psi_{\xi\x}\|^2_{L^2(\bbr)}\,d\tau
                     		\right)^2.
                     	\end{aligned}
                      \end{equation}
            As for $I_{33}$, by Lemmas \ref{trans} and \ref{trans-R}, together with 
              \begin{equation*}
            		|\tilde I_{33}|\leq C(|\bar{v}_\xi|  |\psi_\xi|  |(\phi_{\xi},\vartheta_{\x})| + |\bar{v}_\xi|  |\psi_{\xi}|^2),
            \end{equation*}  
            and 
            $$|\bar{v}_\xi| \le C|v^R_\x|+|(v^S)^{-\mb X}_\x|,$$
                  we get 
            \begin{align*}
                    		I_{33} 
                    		&\le C \int_0^t  \int_{\mathbb{R}} | (v^S)_\xi^{-\mathbf{X}} | \left|\int_0^{\tau} |\bar{v}_\xi| (  |\psi_\xi|  |(\phi_{\xi},\vartheta_{\x})| +   |\psi_{\xi}|^2)(\eta, \xi+\sigma (\tau-\eta)) d\eta\right|^2 d\xi d\tau \\
                    	    &\leq C\deltas\int_0^t  \int_{\mathbb{R}}  ( |(\phi_{\xi},\vartheta_{\x})|^2| \psi_\x|^2+| \psi_{\x}|^4)d\xi d\tau.
            \end{align*}
               Similarly, using
                    \[
                    \|(\phi_{\xi},\vartheta_{\x})\psi_\xi\|_{L^2(\bbr)}
                    \le
                    \|\psi_\xi\|_{L^\infty(\bbr)} \|(\phi_{\xi},\vartheta_{\x})\|_{L^2(\bbr)}
                    \]
                    and 
                    \[
                    \| \psi_{\xi}^2\|_{L^2(\bbr)}
                    \le
                    \|\psi_\xi\|_{L^\infty(\bbr)}\|\psi_{\xi}\|_{L^2(\bbr)},
                    \]
                    together with Sobolev's inequality and Lemma \ref{lem:variables-equivalence}, we have
                    \begin{equation}\label{I33}
                    	\begin{aligned}
                    		I_{33}
                    		&\le
                    		C(\delta_0+\eps_1)
                    		\left(
                    		\int_0^t
                    		(\|\psi_\xi\|^2_{L^2(\bbr)}+ \| \psi_{\xi\x}\|^2_{L^2(\bbr)})\,d\tau
                    		\right).
                    	\end{aligned}
                    \end{equation}
                    Thus, combined with \eqref{I31}-\eqref{I33}, the estimate for $I_{3}$ can be obtained as follows
                    \begin{equation}\label{I3}
                     	\begin{aligned}
                     		I_{3}
                     		\le
                     		& C\deltas\int_0^t(\mathcal{G}^S+\mathcal G^R)d\tau+C\deltas\int_0^t  \int_{\mathbb{R}} | (v^S)_\xi^{-\mathbf{X}}| (|\phi_\x|^2+|\vartheta_\x|^2)d\xi d\tau\\
                    		& +C\deltas\int_0^t\int_{\bbr}  |v^R_\xi|^2(|\phi_\xi|^2 + |\vartheta_\xi|^2) d\x d\tau\\
                     		&+
                     		C(\delta_0+\eps_1)
                     		\left(
                     		\int_0^t
                     		\|\psi_\xi\|^2_{H^1(\bbr)}\,d\tau
                     		\right)
                     		+
                     		C\deltas
                     		\left(
                     		\int_0^t
                     		\|\psi_\xi\|^2_{H^1(\bbr)}\,d\tau
                     		\right)^2.
                     	\end{aligned}
                     \end{equation}
                   For $I_4$, using Lemmas \ref{lem:rare-v}, \ref{lem:wave-inter} and \ref{trans-1}, we have 
                   \begin{equation}\label{I4}
                   	\begin{aligned}
                   		I_4&=\int_0^t  \int_{\mathbb{R}} | (v^S)_\xi^{-\mathbf{X}} | \left|\int_0^{\tau} R_{3,\x}(\eta, \xi+\sigma (\tau-\eta)) d\eta\right|^2 d\xi d\tau\\
                   		&\leq C\deltas\left(\int_0^t\|R_{3,\x}(\tau,\cdot)\|_{L^2(\bbr)}\,d\tau\right)^2\leq C\delta_0.   		
                   	\end{aligned}
                   \end{equation}
                  Subsequently, combining \eqref{I1}, \eqref{I2}, \eqref{I3} and \eqref{I4}, we can obtain
                   	\begin{align*}
                   		& \int_0^t  \int_{\mathbb{R}} | (v^S)_\xi^{-\mathbf{X}} |  \zeta^2_{\x} (\tau,\x)d\xi d\tau\\
                   		&\leq  C\deltas \| \zeta_{0\xi} \|_{L^2(\bbr)}^2+C\deltas^2\int_0^t  | \dot{\mathbf{X}}(\tau)|^2d\tau+C\deltas\int_0^t(\mathcal{G}^S+\mathcal G^R)d\tau\\
                     	& +C\deltas\int_0^t  \int_{\mathbb{R}} | (v^S)_\xi^{-\mathbf{X}}| (|\phi_\x|^2+|\vartheta_\x|^2)d\xi d\tau
                        +C\deltas\int_0^t\int_{\bbr}  |v^R_\xi|^2(|\phi_\xi|^2 + |\vartheta_\xi|^2) d\x d\tau\\
                     	&+
                     	C(\delta_0+\eps_1)
                     	\left(
                     	\int_0^t
                     	\|\psi_\xi\|^2_{H^1(\bbr)}\,d\tau
                     	\right)
                     	+
                     	C\deltas
                     	\left(
                     	\int_0^t
                     	\|\psi_\xi\|^2_{H^1(\bbr)}\,d\tau
                     	\right)^2+ C\delta_0,
                   	\end{align*}
                   which together with Lemmas \ref{Re-en}, \ref{lem:variables-weighted}, the {\it a priori} assumption \eqref{apri-ass} and the smallness of $\deltas$, implies
                   	\begin{align*}
                   		& \int_0^t  \int_{\mathbb{R}} | (v^S)_\xi^{-\mathbf{X}} |  \zeta^2_{\x} (\tau,\x)d\xi d\tau\\
                   		&\leq  C\|(W,\psi,\zeta)(0,\cdot)\|_{H^1(\mathbb{R})}^2 
                        +C\deltas\int_0^t\int_{\bbr}  |v^R_\xi|^2(|W_\xi|^2 + |\zeta_\xi|^2) d\x d\tau\\
                     	&\quad+C\delta_0 \int_0^t  \int_{\mathbb{R}} |W_\x|^2d\xi d\tau+
                     	C(\delta_0+\eps_1)
                     	\left(
                     	\int_0^t
                     	\|\psi_{\x\x}\|^2_{L^2(\bbr)}\,d\tau
                     	\right)
                     	+
                     	C\delta_0^{\frac 14}.
                   	\end{align*}
                   
                   Therefore, we complete the proof of the lemma.
                   \end{proof}
                   \subsection{Conclusion}
                   The following lemma provides the complete first-order derivative estimates. Together with Lemma \ref{Re-en}, we can complete the energy estimates.
                   \begin{lemma}\label{one-order}
                   	Under the hypotheses of Proposition \ref{prop:priori}, there exists $C > 0$ such that for all $t \in (0, T]$,
                   		\begin{align*}
                   			&\|(W_\x,\psi_\x,\zeta_\x)(t,\cdot)\|_{L^2(\mathbb{R})}^2 +\int_0^t   \|W_\x\|^2_{L^2(\bbr)} d\tau+
                   			\int_0^t \|\psi_{\x\x}\|_{L^2(\mathbb{R})}^2  d\tau\\
                   			&\leq C\|(W,\psi,\zeta)(0,\cdot)\|_{H^1(\mathbb{R})}^2+C\delta_0^{\frac 14}.	
                   		\end{align*}
                   \end{lemma}
           \begin{proof}
            By Lemmas \ref{Re-en}, \ref{1-d} and \ref{weight-s}, together with the smallness of $\deltas, \delta_0$ and $\eps_1$, we can obtain
               \begin{align*}
               			&\|(W_\x,\psi_\x,\zeta_\x)(t,\cdot)\|_{L^2(\mathbb{R})}^2 +\int_0^t   \|W_\x\|^2_{L^2(\bbr)} d\tau+
               			\int_0^t \|\psi_{\x\x}\|_{L^2(\mathbb{R})}^2  d\tau\\
               			&\le
               			C\|(W,\psi,\zeta)(0,\cdot)\|_{H^1(\mathbb{R})}^2
               			+C\int_0^t\int_{\bbr}|(v^S)_{\xi}^{-\mathbf{X}}| |\zeta_\x|^2  
               			d\x d\tau+C\delta_0^{\frac 14}\\
                        &\leq C\|(W,\psi,\zeta)(0,\cdot)\|_{H^1(\mathbb{R})}^2
               			+C\deltas\int_0^t\int_{\bbr}  |v^R_\xi|^2(|W_\xi|^2 + |\zeta_\xi|^2) d\x d\tau+C\delta_0^{\frac 14}.
               		\end{align*}
            Then, using Lemma \ref{lem:rare-v} and Gronwall's inequality, we have
            \begin{align*}
               	&\|(W_\x,\psi_\x,\zeta_\x)(t,\cdot)\|_{L^2(\mathbb{R})}^2 +\int_0^t   \|W_\x\|^2_{L^2(\bbr)} d\tau+
               	\int_0^t \|\psi_{\x\x}\|_{L^2(\mathbb{R})}^2  d\tau\\
                &\leq C\|(W,\psi,\zeta)(0,\cdot)\|_{H^1(\mathbb{R})}^2
               	+C\delta_0^{\frac 14}.
            \end{align*}
           \end{proof}     
    Combining Lemmas \ref{Re-en} and \ref{one-order},  we have
    \begin{align*}
        &\|(p-\bar p,u-\bar u,s-\bar s)(t,\cdot)\|_{H^1(\mathbb{R})}^2 +\deltas\int_0^t|\dot{\mb{\mathbf{X}}}(\tau)|^2 d\tau+\int_0^t   \|(u-\bar u)_\x\|^2_{H^1(\bbr)} d\tau\\
        &+\int_0^t   \|(p-\bar p)_\x\|^2_{L^2(\bbr)}d\tau+\int_0^t\mathcal G^R(U)d\tau+\int_0^t\mathcal G^S(U)d\tau\\
        &\le C \|(p-\bar p,u-\bar u,s-\bar s)(0,\cdot) \|_{H^1(\bbr)}^2+C\delta_0^{\frac 14}.
    \end{align*}
     Therefore, we complete the proof of Proposition \ref{prop:priori}.

\section*{Acknowledgment}
L.-A. Li was partially supported by the National Natural Science Foundation of China (No. 12501295), Beijing Natural Science Foundation (No. 1254045) and the Fundamental Research Funds for the Central Universities (No. 2243100008). J. Wu was partially supported by the National Science Foundation of the United States under Grants DMS-2104682 and DMS-2309748. X. Xu was partially supported by the National Natural Science Foundation of China (grants 12571244, 12171040) and the National Key R\&D Program of China (grant 2020YFA0712900).

\bibliographystyle{plain}
\bibliography{reference}

@article{KVW25,
	author = {Kang, M.-J. and Vasseur, A. and Wang, Y.},
	date = {2025/06/25},
	doi = {10.1007/s00205-025-02116-w},
	id = {Kang2025},
	isbn = {1432-0673},
	journal = {Arch. Ration. Mech. Anal.},
	number = {4},
	pages = {42},
	title = {Time-Asymptotic Stability of Generic {R}iemann Solutions for Compressible {N}avier--{S}tokes--{F}ourier Equations},
	url = {https://doi.org/10.1007/s00205-025-02116-w},
	volume = {249},
	year = {2025}}

@article {KVW23,
	AUTHOR = {Kang, M.-J. and Vasseur, A. and Wang, Y.},
	TITLE = {Time-asymptotic stability of composite waves of viscous shock
	and rarefaction for barotropic {N}avier-{S}tokes equations},
	JOURNAL = {Adv. Math.},
	FJOURNAL = {Adv. Math.},
	VOLUME = {419},
	YEAR = {2023},
	PAGES = {Paper No. 108963, 66},
	ISSN = {0001-8708,1090-2082},
	MRCLASS = {35Q30 (76N06 76P05)},
	MRNUMBER = {4560999},
	MRREVIEWER = {Ond\v rej\ Kreml},
	DOI = {10.1016/j.aim.2023.108963},
	URL = {https://doi.org/10.1016/j.aim.2023.108963},
}

@article{HLM,
	author = {Huang, F. and Li, J. and Matsumura, A.},
	date = {2010/07/01},
	doi = {10.1007/s00205-009-0267-0},
	id = {Huang2010},
	isbn = {1432-0673},
	journal = {Arch. Ration. Mech. Anal.},
	number = {1},
	pages = {89--116},
	title = {Asymptotic Stability of Combination of Viscous Contact Wave with Rarefaction Waves for One-Dimensional Compressible {N}avier-{S}tokes System},
	url = {https://doi.org/10.1007/s00205-009-0267-0},
	volume = {197},
	year = {2010}}

@article{Xin-c,
	author  = {T.-P. Liu and Z. Xin},
	title   = {Pointwise decay to contact discontinuities for systems of viscous conservation laws},
	journal = {Asian J. Math.},
	year    = {1997},
	volume  = {1},
	pages   = {34--84},
}

@article{MN86,
	author = {Matsumura, A. and Nishihara, K.},
	date = {1986/06/01},
	doi = {10.1007/BF03167088},
	id = {Matsumura1986},
	isbn = {0910-2043},
	journal = {Japan J. Appl. Math.},
	number = {1},
	pages = {1--13},
	title = {Asymptotics toward the rarefaction waves of the solutions of a one-dimensional model system for compressible viscous gas},
	url = {https://doi.org/10.1007/BF03167088},
	volume = {3},
	year = {1986}}

@article{Nash,
    author  = {Nash, J.},
    title   = {{Le probl{\`e}me de {Cauchy} pour les {\'e}quations diff{\'e}rentielles d'un fluide g{\'e}n{\'e}ral}},
    journal = {Bull. Soc. Math. France},
    volume  = {90},
    pages   = {487--497},
    year    = {1962},
    doi     = {10.24033/bsmf.1586},
}

@article{MN85,
	author = {Matsumura, A. and Nishihara, K.},
	title   = {On the stability of travelling wave solutions of a one-dimensional model system for compressible viscous gas},
	journal = {Japan J. Appl. Math.},
	year    = {1985},
    number = {1},
	volume  = {2},
	pages   = {17--25},
    }

@article{Good,
		author = {Goodman, J.},
		title   = {Nonlinear asymptotic stability of viscous shock profiles for conservation laws},
		journal = {Arch. Ration. Mech. Anal.},
		year    = {1986},
		volume  = {95},
        number  = {4},
		pages   = {325--344},
	}

@article{Liu85,
		author = {Liu, T.-P.},
		title   = {Nonlinear stability of shock waves for viscous conservation laws},
		journal = {Mem. Amer. Math. Soc.},
		year    = {1985},
		volume  = {56},
        number  = {328},
		pages   = {v+108},
	}

@article{LZ15,
		author = {Liu, T.-P. and Zeng, Y.},
		title   = {Shock waves in conservation laws with physical viscosity},
		journal = {Mem. Amer. Math. Soc.},
		year    = {2015},
		volume  = {234},
        number  = {1105},        
		pages   = {vi+168},
	}

@article{SX93,
	author = {Szepessy, A. and Xin, Z.},
	title   = {Nonlinear stability of viscous shock waves},
	journal = {Arch. Ration. Mech. Anal.},
	year    = {1993},
	volume  = {122},
    number  = {1},
	pages   = {53--103},
}

@article{Ma13,
  author  = {S. Ma},
  title   = {Asymptotic Stability of Contact Waves for the Compressible {N}avier-{S}tokes Equations},
  journal = {Appl. Anal.},
  volume  = {92},
  number  = {11},
  pages   = {2320--2334},
  year    = {2013}
}

@article{HM,
   	author = {Huang, F. and Matsumura, A.},
   	title   = {Stability of a composite wave of two viscous shock waves for the full compressible {N}avier-{S}tokes equation},
   	journal = {Comm. Math. Phys. },
   	year    = {2009},
   	volume  = {289},
    number  = {3},
   	pages   = {841--861},
   }

@article{HKK,
	author = {Han, S. and  Kang, M.-J. and  Kim, J.},
	title   = {Large-time behavior of composite waves of viscous shocks for the barotropic {N}avier-{S}tokes equations},
	journal = {SIAM J. Math. Anal.},
	year    = {2023},
	volume  = {55},
    number  = {5},
	pages   = {5526--5574},
}

@article{HXY,
	author = {Huang, F. and  Xin, Z. and Yang, T.},
	title   = {Contact discontinuity with general perturbations for gas motions},
	journal = {Adv. Math. },
	year    = {2008},
	volume  = {219},
	pages   = {1246--1297},
}

@article{HMX,
	author = {F. Huang and A. Matsumura and Z. Xin},
	title   = {Stability of contact discontinuities for the 1-{D} 
	compressible {N}avier-{S}tokes equations},
	journal = {Arch. Ration. Mech. Anal.},
	year    = {2006},
	volume  = {179},
    number  = {1},
	pages   = {55--77},
}

@article{MN92,
	author = {Matsumura, A. and Nishihara, K.},
	title   = {Global stability of the rarefaction wave of a one-dimensional model system for compressible viscous gas},
	journal = {Comm. Math. Phys.},
	year    = { 1992},
	volume  = {144},
    number  = {2},
	pages   = {325--335},
}

@article{NYZ04,
	author = {Nishihara, K. and  Yang, T. and Zhao, H.},
	title   = {Nonlinear stability of strong rarefaction waves for compressible {N}avier-{S}tokes equations},
	journal = {SIAM J. Math. Anal.},
	year    = {2004},
	volume  = {35},
    number  = {6},
	pages   = {1561--1597},
}

@article{LX88,
	author = {Liu, T.-P. and Xin, Z.},
	title   = {Nonlinear stability of rarefaction waves for compressible {N}avier-{S}tokes equations},
	journal = {Comm. Math. Phys.  },
	year    = {1988},
	volume  = {118},
    number  = {3},
	pages   = {451--465},
}

@article{Riemann,
    author  = {Riemann, B.},
    title   = {{\"Uber die Fortpflanzung ebener Luftwellen von endlicher Schwingungsweite}},
    journal = {Abh. K{\"o}nigl. Ges. Wiss. G{\"o}ttingen},
    volume  = {8},
    pages   = {43--65},
    year    = {1860},
}

@article{Lax,
	author = {Lax, P. D.},
	title   = {Hyperbolic systems of conservation laws {II}},
	journal = {Comm. Pure Appl. Math.},
	year    = {1957},
	volume  = {10},
	pages   = {537--566},
}

@article{FM,
	author = { L. Fan and A. Matsumura},
	title   = {Asymptotic stability of a composite wave of two viscous  shock waves for a one-dimensional system of  non-viscous and heat-conductive ideal gas},
	journal = {J. Differential Equations},
	year    = {2015},
	volume  = {258},
	pages   = {1129--1157},
}

@article{FGT,
	author = {Fan, L. and Gong, G. and Tang, S.},
	title = {Asymptotic stability of viscous contact wave and rarefaction waves for the system of heat-conductive ideal gas without viscosity},
	journal = {Anal. Appl. (Singap.)},
	volume = {17},
	number = {02},
	pages = {211--234},
	year = {2019},
}

@phdthesis{SK2,
    author = {S. Kawashima},
    title = {Systems of a Hyperbolic-Parabolic Composite Type, with Applications to the Equations of Magnetohydrodynamics},
    school = {Kyoto University},
    year = {1984},
    url = {http://hdl.handle.net/2433/97887},
}

@article{SK1,
	author = {S. Kawashima},
	title = {Large-time behaviour of solutions to hyperbolic-parabolic systems of conservation laws and applications},
	journal = {Proc. Roy. Soc. Edinburgh Sect. A},
	volume = {106},
	pages = {169--194},
	year = {1987},
}

@article{MW,
	author = {Ma, S. and Wang, J.},
	title = {Decay rates to viscous contact waves for the compressible {N}avier-{S}tokes equations},
	journal = {J. Math. Phys.},
	volume = {57},
	number = {2},
	pages = {021501},
	year = {2016},
	issn = {0022-2488},
	doi = {10.1063/1.4938574},
	url = {https://doi.org/10.1063/1.4938574},
}

@article{HW16,
	author = {F. Huang and T. Wang},
	title = {Stability of superposition of viscous contact wave and 
	rarefaction waves for compressible {N}avier-{S}tokes system},
	journal = {Indiana Univ. Math. J.},
	volume  = {65},
    number  = {6},
	pages = {1833--1875},
	year = {2016},
}

@article{Dafermos1,
	author = {Dafermos, C.M.},
	title = {The second law of thermodynamics and stability},
	journal = {Arch. Ration. Mech. Anal.},
	volume  = {70},
    number  = {2},
	pages = {167--179},
	year = {1979},
}

@article{KVJEMS21,
	author = {Kang, M.-J. and Vasseur, A.},
	title = {Contraction property for large perturbations of shocks of the barotropic {N}avier-{S}tokes system},
	journal = {J. Eur. Math. Soc.},
	volume  = {23},
    number  = {2},
	pages = {585--638},
	year = {2021},
}

@article{Hou,
  author  = {M. Hou},
  title   = {Asymptotic stability of the rarefaction wave for the non-viscous and heat-conductive ideal gas in half space},
  journal = {Acta Math. Sci.},
  volume  = {39},
  number  = {4},
  pages   = {1195--1212},
  year    = {2019},
  doi     = {10.1007/s10473-019-0421-1}
}

@article{LZ99,
  author  = {T.-P. Liu and Y. Zeng},
  title   = {Compressible {N}avier-{S}tokes Equations with Zero Heat Conductivity},
  journal = {J. Differential Equations},
  volume  = {153},
  number  = {2},
  pages   = {225--291},
  year    = {1999}
}

@article{HY18,
  author  = {J. Hu and H. Yin},
  title   = {Nonlinear Stability of Rarefaction Waves for the Compressible {N}avier-{S}tokes Equations with Zero Heat Conductivity},
  journal = {Nonlinear Anal.},
  volume  = {174},
  pages   = {242--277},
  year    = {2018}
}

@book{Smoller,
	author    = {J. Smoller},
	title     = {Shock Waves and Reaction-Diffusion Equations},
	edition = {second},
	publisher = {Springer-Verlag},
	year      = {1994},
}

@book{LiuShockWaves,
  author    = {T.-P. Liu},
  title     = {Shock Waves},
  series    = {Grad. Stud. Math.},
  volume    = {215},
  publisher = {Amer. Math. Soc.},
  address   = {Providence, RI},
  year      = {2021}
}

@book{VK65,
  author    =  {Vincenti, W. G. and Kruger, Jr., C. H.},
  title     = {Introduction to Physical Gas Dynamics},
  publisher = {John Wiley \& Sons},
  address   = {New York},
  year      = {1965}
}

@misc{HL25,
  author = {X. Huang and H. Lee},
  title  = {Time-asymptotic stability of composite wave of viscous shocks
            and viscous contact wave for {N}avier--{S}tokes--{F}ourier equations},
  year   = {2025},
  note   = {arXiv:2504.04014},
}

@misc{GLWX,
  author = {Y. Guo and L.-A. Li and J. Wu and X. Xu},
  title  = {Time-asymptotic stability of generic {R}iemann solutions
            for the system of heat-conductive ideal gas without viscosity},
  year   = {2026},
  note   = {arXiv:2609.02119},
}

@article{ShK85,
  author  = {Shizuta, Y. and Kawashima, S.},
  title   = {Systems of equations of hyperbolic-parabolic type with applications to the discrete {B}oltzmann equation},
  journal = {Hokkaido Math. J.},
  volume  = {14},
  number  = {2},
  pages   = {249--275},
  year    = {1985},
}

\end{document}